\documentclass[onecolumn,twoside,a4paper,11pt,smallheadings,headinclude,headsepline,cleardoubleempty,DIV11,DIVclassic,BCOR2mm,final]{scrbook}
\usepackage{scrpage2}
\usepackage[centertags]{amsmath}
\usepackage{amssymb,amsthm}
\usepackage[all]{xy}
\usepackage[latin1]{inputenc}
\usepackage[ngerman]{babel}
\usepackage{titlesec}
\usepackage{graphicx}

\CompileMatrices

\begin{document}
\newtheoremstyle{myplain}
    {\topsep}{\topsep}
    {\itshape}
    {}{\bfseries}{.\, ---}{.5em}
    {\thmname{#1} \thmnumber {#2}\thmnote{ (#3)}}
\theoremstyle{myplain}
\newtheorem{thm}{Theorem}
\newtheorem{kor}[thm]{Korollar}
\newtheorem{lem}[thm]{Lemma}
\newtheorem{prop}[thm]{Proposition}
\newtheorem{satz}[thm]{Satz}

\newtheoremstyle{mydefinition}
    {\topsep}{\topsep}
    {}
    {}{\bfseries}{.\, ---}{.7em}
    {\thmname{#1}\thmnumber{ #2}\thmnote{ (#3)}}
\theoremstyle{mydefinition}

\newtheorem{defn}[thm]{Definition}
\newtheorem{bem}[thm]{Bemerkung}


\newcommand{\lto}{\longrightarrow}
\newcommand{\opn}{\operatorname}
\newcommand{\id}{\opn{id}}
\newcommand{\dirlim}{\opn{dirlim}}
\newcommand{\Proj}{\opn{Proj}}
\newcommand{\Cyl}{\opn{Cyl}}
\newcommand{\Wh}{\opn{Wh}}
\newcommand{\Whp}[1]{\Wh(\pi #1)}
\newcommand{\Torus}{\opn{Torus}}
\newcommand{\copi}{\coprod_{i\in I_n}}
\newcommand{\copii}{\coprod_{i\in I'_n}}
\newcommand{\arincl}{\ar@{^{(}->}}
\newcommand{\arinclinv}{\ar@{_{(}->}}
\newcommand{\norm}[1]{\Vert #1\Vert}
\newcommand{\nosep}{\setlength{\topsep}{0pt}\setlength{\partopsep}{0pt}\setlength{\itemsep}{0pt}\setlength{\parskip}{0pt}}
\newcommand{\smallsep}{\setlength{\itemsep}{0pt}}

\newcommand{\RR}{\mathbb{R}}
\newcommand{\ZZ}{\mathbb{Z}}
\newcommand{\NN}{\mathbb{N}}


\bibliographystyle{plain}
\setcounter{secnumdepth}{1}
\numberwithin{thm}{chapter}

\SelectTips{cm}{}
\renewcommand{\theenumi}{\roman{enumi}}
\renewcommand{\labelenumi}{\textup{(}\theenumi\textup{)}}

\renewcommand{\theequation}{\arabic{equation}}
\renewcommand*{\captionformat}{}
\renewcommand{\figurename}{Fig.}

\setkomafont{pagehead}{\upshape\sffamily}
\clearscrheadings
\clearscrplain
\rehead[]{\leftmark}
\lohead[]{\rightmark}
\ohead{\thepage}


\frontmatter


\titlehead{
\begin{center}\scshape

\includegraphics[width=0.4\textwidth,height=2cm]{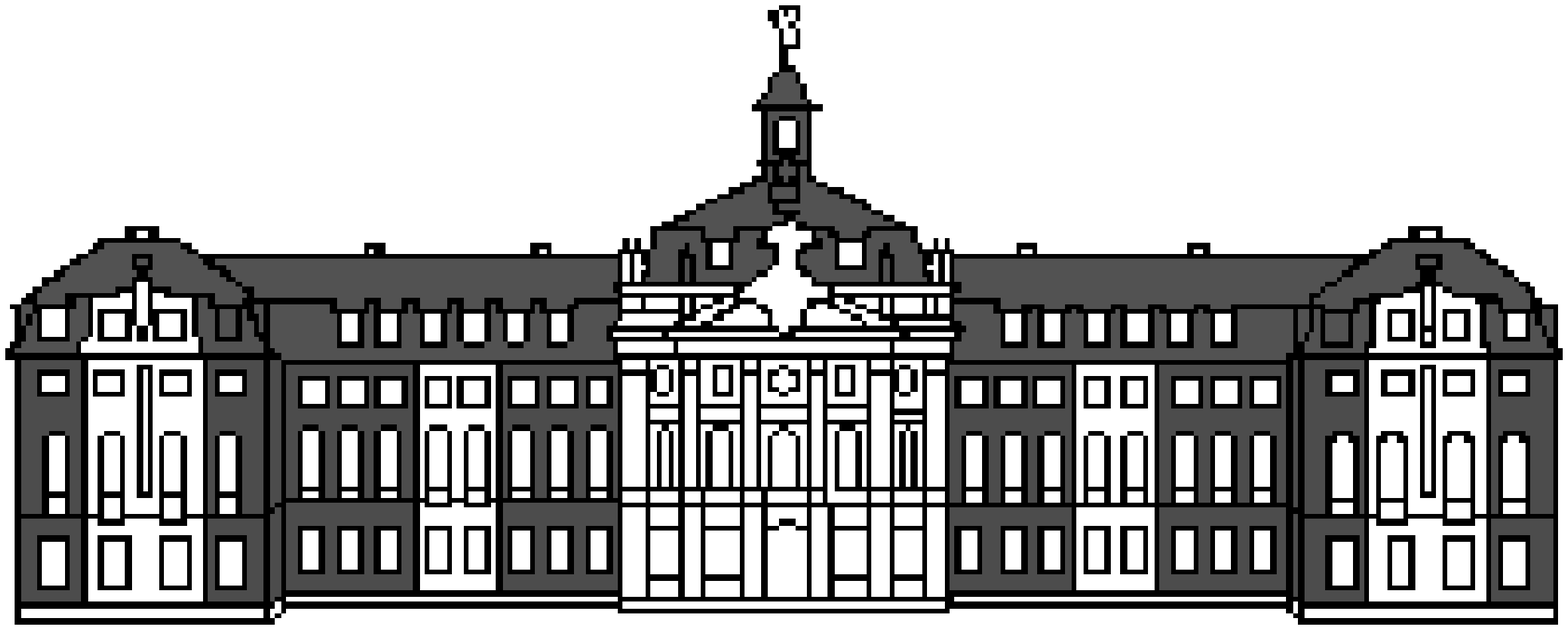}

\vspace{1.5ex}
{\large Westf\"{a}lische Wilhelms-Universit\"{a}t M\"{u}nster}

{Fachbereich Mathematik und Informatik}

\end{center}
\vspace{2cm}
}

\title{Whitehead-Torsion und Faserungen}

\author{Wolfgang Steimle}
\date{}

\publishers{\vspace{4cm}\large\vspace{1cm}\sffamily\bfseries Diplomarbeit

\vspace{1.5ex}
\rmfamily\mdseries zur Erlangung des Grades eines

Diplom-Mathematikers

\vspace{1.5ex}
Betreut durch Prof.~Dr.~Wolfgang Lück

Vorgelegt am 9.~Februar 2007 
}
\maketitle

\chapter*{Abstract}
This work treats on the question whether a given map $f:\, M \lto B$  of smooth closed manifolds is homotopic to a smooth fiber bundle. We define a first obstruction in $H^1(B;\Whp{E})$ and, provided that this obstruction vanishes and one additional condition is verified, a second obstruction in $\Whp{E}$. Both elements vanish if the answer to the above question is positive. In the case where $B$ is the 1-sphere and the dimension of $M$ exceeds five, we show that the converse is also true, using a relationship with two obstructions defined by Farrell in this particular situation.\\

\emph{Mathematics Subject Classification (2000):} 19J10, 55R15, 57R19.


\chapter*{Einleitung}

Im Jahr 1966 veröffentlichten W.~Browder und J.~Levine einen Artikel "`Fibering manifolds over a circle"' \cite{browder_levine}, in dem sie die folgende Frage untersuchten: Wann ist eine Abbildung $f:\,M\lto S^1$ von differenzierbaren kompakten Mannigfaltigkeiten homotop zu einem differenzierbaren Faserbündel? Dabei beschränkten sich die Autoren auf den Spezialfall, dass die zu erwartende Faser (deren Homotopietyp ja durch die lange exakte Faserungssequenz feststeht) 1-zusammenhängend ist. Unter diesen Voraussetzungen bewiesen sie unter Verwendung des 5 Jahre zuvor von Smale gefundenen $h$-Kobordismensatzes, dass in der Tat jede nicht nullhomotope Abbildung homotop zu einem differenzierbaren Faserbündel ist, sofern die Dimension von $M$ mindestens sechs ist.

Die Verallgemeinerung der Situation auf nicht einfach zusammenhängende Fasern wurde 1971 von F.~T.~Farrell in seiner Dissertation "`The Obstruction to Fibering an Manifold over a Circle"' gegeben, die in Artikelform in \cite{farrell} veröffentlicht ist. Tatsächlich ist das Resultat von Browder und Levine in diesem allgemeineren Kontext nicht mehr richtig. Dies legt bereits die Anwendung des $h$-Kobordismensatzes nahe, bei dessen Verallgemeinerung für nicht einfach zusammenhängende Mannigfaltigkeiten ein Hindernis in Form einer Whitehead-Torsion auftritt. Farrell definierte in seiner Arbeit zwei $K$-theoretische Hindernisse, die einer homotopen Abänderung der gegebenen Abbildung in ein Faserbündel entgegenstehen. Jedes der Hindernisse hat eine sehr konkrete geometrische Bedeutung; das Verschwinden beider impliziert bereits, dass die gesuchte homotope Abänderung möglich ist.

Die Argumentation von Browder/Levine und Farrell ist allerdings an die spezielle Situation gebunden, in der der Zielraum der Abbildung $f$ tatsächlich die 1-Sphäre ist, sodass eine weitere Verallgemeinerung in Richtung einer größeren Klasse von Bildräumen nicht ohne Weiteres möglich ist. Bereits die naheliegende Frage, welche Entsprechung Farrells Hindernisse für andere Bildräume $B$ hat, hat keine offensichtliche Antwort. 

In dieser Diplomarbeit wird eine Verallgemeinerung der von Farrell definierten Hindernisse auf eine größere Klasse von Abbildungen $f:\,M\lto B$ von geschlossenen zusammenhängenden differenzierbaren Mannigfaltigkeiten vorgeschlagen. Diese Klasse umfasst insbesondere den Fall, dass die Eulercharakteristik von $B$ verschwindet und beinhaltet damit die von Farrell betrachtete Situation. Konkret werden wir ein Element $\theta(f)\in H^1(B;\Whp{M})$ und, im Falle des Verschwindens von $\theta(f)$, ein Element $\tau(f)\in\Whp{M}$ definieren. Beide Elemente hängen nur von der Homotopieklasse der betrachteten Abbildung $f$ ab und verschwinden, falls $f$ ein differenzierbares Faserbündel ist. In der Tat ist das erste Hindernis sogar unter sehr allgemeinen Bedingungen definiert; die notwendigen Einschränkungen an die Klasse von betrachteten Abbildungen $f:\,M\lto B$ sind nötig, um die Wohldefiniertheit des zweiten Hindernisses sicherzustellen. 

Es sei nun kurz die Idee beschrieben, die der Definition unserer beiden Hindernisse zu Grunde liegt. Leitmotiv ist, statt der Abbildung $f$ deren "`assoziierte Faserung"' zu betrachten. Diese Standardkonstruktion der Homotopietheorie ordnet einer beliebigen Abbildung $f$ eine Faserung $p_f:\,M^f\lto M$ mit einem zu $M$ homotopie-äquivalenten Totalraum zu. Falls eine Homotopie von $f$ zu einem differenzierbaren Faserbündel existiert, ist zudem die Faser dieser Faserung vom Homotopietyp eines endlichen CW-Komplexes. Diese Endlichkeitsbedingung an die Faser ist die erste Voraussetzung für die Definition unserer beiden Elemente. 

Grob gesprochen misst das Element $\theta(f)$ die Einfachheit von Fasertransporten der assoziierten Faserung entlang von geschlossenen Wegen in $B$, während das Element $\tau(f)$ die Einfachheit von Trivialisierungen von $M^f$ über kontraktiblen Teilmengen von $B$ ausdrückt. Hierbei drückt "`Einfachheit"' jeweils das Verschwinden einer Whitehead-Torsion aus. Allerdings haben wir durch den Übergang zur assoziierten Faserung die Mannigfaltigkeits-Struktur von $M$ verloren; die Räume der assoziierten Faserung sind nur noch vom Homotopietyp eines endlichen CW-Komplexes. Daher ist a priori eine Whitehead-Torsion von Fasertransporten und Trivialisierungen nicht definiert. Für die Definition des ersten Elementes $\theta(f)$ ist dies allerdings, wie wir sehen werden, kein Problem. 

Anders ist dies beim zweiten Element: Einschränkungen von $M^f$ auf zusammenziehbare Teilmengen von $B$ haben im Allgemeinen keinen einfachen Homotopietyp, der mit entsprechenden Teilmengen von $M$ in Verbindung steht. Jedoch kann man geeignete lokale Trivialierungen des Totalraums der assoziierten Faserung zu einem endlichen CW-Modell $m:\,X\lto M^f$ zusammenfügen. Dieses kann nun mit dem durch $M$ selbst gegebenen endlichen CW-Modell $\lambda$ von $M^f$ verglichen werden. Die Torsion der "`Vergleichsabbildung"' $\lambda\circ x^{-1}$ definiert nun, nach Rücktransport in die Whitehead-Gruppe von $M$, das Element $\tau(f)$.  

Die exakte Form des CW-Modells $m$ spielt dabei offenbar überhaupt keine Rolle, da es ohnehin nur der Berechnung von Torsionen dient. Um diese Flexibilität zu nutzen, bietet es sich daher an, statt endlicher CW-Modelle von $M^f$ gleich "`einfache Strukturen"' auf $M^f$ zu betrachten. Grob gesprochen ist dies die Wahl eines endlichen CW-Modells $m$ auf $M^f$, wobei allerdings zwei endliche CW-Modelle definieren die gleiche einfache Struktur, falls die Torsion der Vergleichsabbildung verschwindet.

Für den von Farrell betrachteten Fall ist ein "`1:1-Vergleich"' unserer beiden Elemente mit den beiden Hindernissen von Farrell nicht möglich. Jedoch steht das hier definierte Element $\tau(f)$ mit Farrells Hindernissen in einem einfachen algebraischen Zusammenhang. Auf diese Weise erkennt man, dass das Element $\tau(f)$ die Information beider Hindernisse von Farrell zusammenfasst; insbesondere impliziert das Verschwinden von $\theta(f)$ und $\tau(f)$ auch das Verschwinden von Farrells Hindernissen.

Diese Arbeit ist wie folgt aufgebaut: Im ersten Kapitel dieser Arbeit werden wir kurz zwei für das Folgende wesentliche Elemente aus der Homotopietheorie von Faserungen in Erinnerung rufen: Zum einen betrachten wir den Fasertransport und gehen in diesem Kontext auch auf die Existenz und Eindeutigkeit von Trivialisierungen von Faserungen über kontraktiblen Basisräumen ein. Zum anderen wiederholen wir die Konstruktion und wesentliche Eigenschaften der assoziierten Faserung.

Im zweiten Kapitel definieren wir dann das erste Element $\theta(f)$. In diesem Zusammenhang bietet es sich an, bereits genauer auf den oben kurz erklärten Begriff der einfachen Struktur auf Räumen vom Homotopietyp eines endlichen CW-Komplexen einzugehen, auch wenn dieser im vollen Umfang erst im anschließenden dritten Kapitel zum Tragen kommt. Dieses dritte Kapitel ist dann der Definition und den Eigenschaften des Hindernisses $\tau(f)$ gewidmet. Hauptarbeit ist dabei die Nachweis der Existenz und Eindeutigkeit einer geeigneten einfachen Struktur auf $M^f$.

Das vierte Kapitel schließlich vergleicht für den Fall $B=S^1$ und $\dim M\geq 6$, also für Farrells Situation, unsere beiden Elemente mit den Hindernissen von Farrell. Das Ziel ist der Nachweis des bereits angesprochenen algebraischen Zusammenhangs. Als unmittelbare Folgerung erhält man, dass in dieser Situation das Verschwinden der in dieser Arbeit definierten Elemente bereits hinreichend dafür ist, die Frage nach der homotopen Abänderbarkeit positiv zu beantworten.

Abschließend möchte ich mich sehr herzlich bedanken bei Herrn Prof.~Wolfgang Lück für das interessante Diplomarbeits-Thema und die hervorragende Betreuung während dieser Arbeit. Mein weiterer Dank gilt der gesamten Arbeitsgruppe Topologie an der Universität Münster, die mich freundlich aufgenommen hat.


\tableofcontents

\mainmatter


\chapter{Faserungen}\label{fas}

In diesem Kapitel sollen einige wohlbekannte Tatsachen über Faserungen wiederholt werden. Das Konzept der Faserung kann man als eine Verallgemeinerung des Konzepts des Faserbündels auffassen, das zwar keine Anforderungen mehr an die lokale Form der Abbildung stellt, allerdings die Eigenschaft des Fasertransportes axiomatisch bewahrt. Dieses Konzept ist allgemein genug, um homotopietheoretisch eine beliebige Abbildung durch eine Faserung "`ersetzen"' zu können --- eine Tatsache, die für uns von Wichtigkeit sein wird.

In einem ersten Teil werden wir den Begriff des Fasertransportes wiederholen. Die Darstellung ist angelehnt an \cite{whitehead}, Seite 29ff.~und \cite{switzer}, Seite 342ff. Anschließend definieren wir die zu einer Abbildung assoziierte Faserung und leiten einige Eigenschaften her.

Es bezeichne stets $I := [0,1]$ das Einheitsintervall.

\section{Der Fasertransport}

Es sei $p:\,E \lto B$ eine Faserung. Wir werden zeigen, dass direkt aus der Homotopie-Hochhebungs-Eigenschaft (HHE) von Faserungen die Existenz eines Fasertransportes folgt. 

Wir betrachten dazu Pull-backs von Faserungen mit einer Basis-Abbildung:
$$\xymatrix{
f^*E \ar[rr]^{\bar{f}} \ar[d]^{p'} & & E \ar[d]^p\\
X \ar[rr]^f & & B
}$$
Dann ist $p'$ wieder eine Faserung. Es sei nun $H:\,X\times I\lto B$ eine Homotopie zwischen $f,g:\,X\lto B$. Wir betrachten das Homotopie-Hochhebungs-Problem (HHP)
\begin{equation}\begin{split}\xymatrix{
f^*E\times\{0\} \ar[rr]^{\bar{f}} \arinclinv[d] & & E\ar[d]^p\\
f^*E\times I \ar[r]_{p'\times\id_I} \ar@{.>}[rru]^{H'} & X\times I \ar[r]_H & B
}\label{fas:HHP_Fasertransport}\end{split}\end{equation}
Die Auswertung der Lösung $H'$ an $f^*E\times\{1\}$ ist eine Abbildung $f^*E\lto E$, die zusammen mit $p':\,f^*E\lto X$ eine Abbildung $\omega_H:\,f^*E\lto g^*E$ definiert.

\begin{satz}\label{fas:fasertransport}
Zu einer Homotopie $H:\,X\times I\lto B$ zwischen $f$ und $g$ sei $\omega_H:\,f^*E\lto g^*E$ wie oben. Dann gilt:
\begin{enumerate}\smallsep
\item Die Abbildungen $\bar{g}\circ\omega_H$ und $\bar{f}:\,f^*E\lto E$ sind homotop.
\item Ist $K:\,X\times I\lto B$ eine zweite Homotopie zwischen $f$ und $g$, die zu $H$ relativ $X\times\partial I$ homotop ist, so sind $\omega_H$ und $\omega_K$ faserhomotop.
\item $\omega_H:\, f^*E\lto g^*E$ ist eine Faserhomotopie-Äquivalenz von Faserungen über $X$.
\end{enumerate}
\end{satz}

\begin{proof}[Beweis]
(i) Nach Konstruktion stellt die Abbildung $H':\,f^*E\times I\lto E$ eine Homotopie dar.

(ii) Sei $L:\,X\times I\times I\lto B$ eine Homotopie zwischen $H$ und $K$ relativ $X\times\partial I$. Wir setzen 
$$\dot{I} := I\times\{0,1\} \cup \{0\}\times I \subset I\times I$$
und bemerken, dass das Paar $(I\times I,\dot{I})$ homöomorph zu $(I\times I, \{0\}\times I)$ ist. Dies erlaubt es uns, die HHE einer Faserung dahingehend zu interpretieren, dass eine Abbildung $\phi:\,A\times\dot{I}\lto E$ zu einer Abbildung $\Phi:\,A\times I\times I\lto E$ so fortgesetzt werden kann, dass $p\circ\Phi:\,A\times I\times I\lto B$ mit einer bereits gegebenen Fortsetzung von $p\circ\phi$ übereinstimmt.

Die bereits gegebene Fortsetzung sei hier die Abbildung $L\circ (p'\times\id_{I\times I}):\,f^*E\times I\times I\lto X\times I\times I\lto B$. Diese soll nach $E$ hochgehoben werden, mit der Maßgabe, auf $f^*E\times I\times \{0\}$ die Homotopie $H'$, auf $f^*E\times I\times \{1\}$ die Homotopie $K'$ und auf $f^*E\times \{0\}\times I$ konstant die Abbildung $\bar{f}:\,f^*E\lto E$ fortzusetzen. Wir bezeichnen mit $M:\,f^*E\times I\times I\lto E$ die gewünschte Fortsetzung. Die Auswertung von $M$ an $f^*E\times\{1\}\times I$ liefert zusammen mit der Abbildung $p'\circ\Proj:\, f^*E\times\{1\}\times f^*E \lto X$ eine Abbildung $f^*E\times\{1\}\times I\lto g^*E$. Es ist ersichtlich, dass es sich um die gesuchte Faserhomotopie handelt.

(iii) Es bezeichne $H^-:\,X\times I\lto B$ die Homotopie $H$ "`rückwärts durchlaufen"' (d.~h.\ $H^-(x,t):=H(x,1-t)$). Die dazu wie oben konstruierbare Abbildung $\omega_{H^-}:\,g^*E\lto f^*E$ ist ein Homotopie-Inverses zu $\omega_H$, wie ein ähnliches Fortsetzungsargument wie in (ii) zeigt: 

Wir betrachten dazu zunächst die stationäre Fortsetzung von $p\circ H':\,f^*E\times I\lto E\lto B$ zu einer Homotopie $f^*E\times I\times I \lto B$. Diese Abbildung soll nun nach $E$ geliftet werden, dieses Mal mit der Maßgabe, auf $f^*E\times I\times \{0\}$ die Abbildung $H'$ fortzusetzen, auf $f^*E\times \{0\}\times I$ konstant die Abbildung $\bar{g}\circ\omega_H:\,f^*E\lto E$ zu durchlaufen und auf $f^*E\times I\times \{1\}$ die Abbildung $\big((H^-)'\big)^-\circ(\omega_H\times\id_I):\,f^*E\times I\lto g^*E\times I\lto E$ fortzusetzen. Es bezeichne wieder $M:\,f^*E\times I\times I\lto E$ die gewünschte Fortsetzung. Wie in Teil (i) sieht man, dass die Auswertung von $M$ nun an $f^*E\times \{0\}\times I$ eine Faserhomotopie über $X$ zwischen $\id_{f^*E}$ und $\omega_{H^-}\circ\omega_H$ induziert.

Dasselbe Argument zeigt, dass auch $\omega_H\circ\omega_{H^-}$ über $X$ faserhomotop zu $\id_{g^*E}$ ist.
\end{proof}

Als Spezialfall dieses Sachverhalts erhält man den "`üblichen"' Fasertransport: Sei dazu $\gamma:\,I\lto B$ ein Weg in $B$. Indem man $\gamma$ als Homotopie zwischen den beiden Punktabbildungen $\star\lto B$ interpretiert, die den Basispunkt nach $\gamma(0)$ bzw.~$\gamma(1)$ schicken, folgt die Existenz einer bis auf Homotopie eindeutige Abbildung $\omega_\gamma:\,F_{\gamma(0)}\lto F_{\gamma(1)}.$ Diese Abbildung hängt nur von der Homotopieklasse relativ $\{0,1\}$ von $\gamma$ ab.

\begin{defn}
Seien $b,b'\in B$. Eine Abbildung $F_{b}\lto F_{b'}$ heißt \emph{Fasertransport,} wenn sie für einen geeigneten Weg $\gamma$ von $b$ nach $b'$ homotop zu einer Abbildung $\omega_\gamma$ wie oben ist. 
\end{defn}

Wir geben nun eine erste Anwendung von Satz \ref{fas:fasertransport} an.

\begin{prop}\label{fas:prop1}
Sei $p:\,E\lto B$ eine Faserung, sei $f:\,B'\lto B$ eine Homotopie-Äquivalenz. Dann ist auch $\bar{f}:\,f^*E\lto E$ eine Homotopie-Äquivalenz.
\end{prop}

\begin{proof}[Beweis]
Es seien $g:\,B\lto B'$ ein Homotopie-Inverses von $f$ sowie $H:\,B\times I\lto B$ bzw.~$K:\,B'\times I\lto B'$ Homotopien zwischen $\id_B$ und $f\circ g$ bzw.~zwischen $\id_{B'}$ und $g\circ f$. Wir zeigen, dass $\bar{g}\circ\omega_H:\,E\lto f^*E$ ein Homotopie-Inverses von $\bar{f}$ ist. Zunächst ist nach Satz \ref{fas:fasertransport}, Teil (i) die Komposition $\bar{f}\circ\bar{g}\circ\omega_H$ homotop zur Identität auf $E$; somit ist $\bar{f}\circ\bar{g}$ eine Homotopie-Inverses von $\omega_H$. Genauso gilt $\bar{g}\circ\bar{f}\circ\omega_K\simeq\id_{f^*E}$. Es folgt 
\begin{multline*}
\bar{g}\circ\omega_H \circ \bar{f} \simeq \bar{g}\circ\omega_H\circ\bar{f}\circ\bar{g}\circ\bar{f}\circ\omega_K \simeq \bar{g}\circ(\bar{f}\circ\bar{g})^{-1} \circ \bar{f}\circ\bar{g}\circ\bar{f}\circ\omega_K \\
\simeq \bar{g}\circ\bar{f}\circ\omega_K \simeq \id_{f^*E}.\tag*{\qedhere}
\end{multline*}
\end{proof}

Proposition \ref{fas:prop1} erlaubt uns, die Abbildung $\omega_H$ aus Satz \ref{fas:fasertransport} noch auf andere Weise zu beschreiben. Wir verwenden wieder die Bezeichnungen von oben und betrachten das Pull-back von $E$ mit der ganzen Homotopie $H:\,X\times I\lto B$. Auf Grund der Funktorialität des Pull-backs sind $f^*E$ und $g^*E$ die Einschränkungen von $H^*E$ auf $X\times\{0\}$ bzw.~$X\times\{1\}$, die wir als Pull-back von $H^*E$ mit den Inklusionen $i_0$ bzw.~$i_1:\,X\times\{i\}\lto X\times I$ interpretieren können. Da $i_0$ und $i_1$ Homotopie-Äquivalenzen sind, gilt gemäß Proposition \ref{fas:prop1} dasselbe für die Inklusionen von $f^*E$ und $g^*E$ nach $H^*E$. Betrachte nun das Diagramm
$$\xymatrix{
f^*E\times\{0\} \arincl[rr]^{\bar{i}_0} \arinclinv[d] && H^*E \ar[rr]^{\bar{H}} \ar[d] && E \ar[d] \\
f^*E\times I \ar[rr]^{p'\times\id_I} \ar@{.>}[rru]^{H''} && X\times I \ar[rr]^H && B
}$$
Das gesamte Diagramm stellt das HHP \eqref{fas:HHP_Fasertransport} dar, insbesondere löst $\bar{H}\circ H''$ das HHP \eqref{fas:HHP_Fasertransport}. Damit ist die Einschränkung von $H''$ auf $f^*E\times \{1\}$ gleich der oben definierte Abbildung $\omega_H:\,f^*E\lto g^*E$, wenn man den Bildbereich nur entsprechend einschränkt. Auch das linke Quadrat ist ein HHP vom Typ \eqref{fas:HHP_Fasertransport}, wenn wir in \eqref{fas:HHP_Fasertransport} als Homotopie $H$ die Identität auf $X\times I$ wählen; die Rolle der Abbildungen $f$ und $g$ übernehmen dann die beiden Inklusionen $i_0,i_1:\,X\lto X\times I$. Mit Satz \ref{fas:fasertransport}, Teil (i) gilt $\bar{i}_1\circ\omega_H\simeq \bar{i}_0$, oder mit anderen Worten:
\begin{equation}\label{fas:gl1}
\omega_H \simeq \bar{i}_1^{-1} \circ \bar{i}_0.
\end{equation}
Im Spezialfall, dass $p$ ein lokal-triviales Faserbündel über einem parakompakten Raum ist, kann die Abbildung $H''$ als ein Homöomorphismus gewählt werden (vgl.~\cite{husemoller}, Seite 28). Auch die Abbildung $\omega_H$ ist dann als Einschränkung eines Homöomorphismus wieder ein Homöomorphismus. Diese Tatsache wird in den folgenden Kapiteln eine Rolle spielen.

\begin{lem}\label{fas:lemma5}
Seien $p:\,E\lto B$ eine Faserung, $\varphi:\,B'\lto B$ eine Abbildung und $H:\,X\times I\lto B'$ eine Homotopie zwischen den Abbildungen $f,g:\,X\lto B'$. Bezeichne mit $p':\,\varphi^*E\lto B'$ die zurückgezogene Faserung und mit $\omega'_H:\,f^*\varphi^*E\lto g^*\varphi^*E$ eine mit $H$ und $p'$ wie oben konstruierte Abbildung. Dann ist $\omega'_H$ homotop zu $\omega_{\varphi\circ H}$.
\end{lem}

Dieses Lemma haben wir implizit bei der Herleitung von \eqref{fas:gl1} bereits bewiesen, jedoch erhält man es auch direkt aus der Charakterisierung \eqref{fas:gl1} des Fasertransports zurück, denn die Pull-backs $H^*\varphi^*E$ und $(\varphi\circ H)^*E$ sind kanonisch homöomorph. --- Als nächstes halten wir fest, dass die Konstruktion von $\omega_H$ mit Einschränkungen im folgenden Sinne verträglich ist. Sei $Y\subset X$ ein Teilraum. Dann ist die Einschränkung $K:=H\vert_{Y\times I}$ eine Homotopie zwischen $f\vert_Y$ und $g\vert_Y$, und man kann dazu wie oben eine Faserhomotopie-Äquivalenz $\omega_K:\,(f\vert_Y)^*E\lto (g\vert_Y)^*E$ konstruieren. 

\begin{lem}\label{fas:lemma4}
Unter den Identifikationen $(f\vert_Y)^*E=(f^*E)\vert_Y$ und $(g\vert_Y)^*E=(g^*E)\vert_Y$ sind $\omega_K$ und $\omega_H\vert_Y$ faserhomotop.
\end{lem}

\begin{proof}[Beweis]
Man überprüft unmittelbar, dass im HHP \eqref{fas:HHP_Fasertransport} die Einschränkung von $H'$ auf $(f\vert_Y)^*E\times I$ eine Lösung des entsprechenden HHP für die Homotopie $K$ ist. Satz \ref{fas:fasertransport}, Teil (ii) liefert die Behauptung.
\end{proof}

Wir betrachten nun den Fall einer nullhomotopen Abbildung $f:\,X\lto B$. Sei $K:\,X\times I\lto B$ eine Homotopie zwischen $K_0=f$ und der Kollapsabbildung, die ganz $X$ auf $b\in B$ abbildet. Als Anwendung von Satz \ref{fas:fasertransport} auf die Homotopie $K$ erhält man mit $(K_0)^*E=f^*E$ und $(K_1)^*E=F_b\times X$ eine Faserhomotopie-Äquivalenz $\omega_K:\, f^*E\lto F_b\times X$.

\begin{defn}
Für eine Abbildung $f:\,X\lto B$ heißt eine Faserhomotopie-Äquiva\-lenz $\varphi:\,f^*E\lto F_b\times X$ von Faserungen über $X$ eine \emph{Trivialisierung von $f^*E$}, wenn die Einschränkung $\varphi:\,F_{f(x)}\lto F_b\times \{x\}$ auf jede Faser ein Fasertransport (in $B$) ist.
\end{defn}

\begin{prop}\label{fas:trivialisierung}
Sei $X$ ein kontrahierbarer Raum, $x\in X$. Sei $f:\,X\lto B$ eine Abbildung.
\begin{enumerate}\setlength{\itemsep}{0cm}
\item Jede Nullhomotopie von $f$ liefert eine Trivialisierung $T:\,f^*E\lto F_b\times X$ von $f^*E$.
\item Sei $T':\,f^*E\lto F_{b'}\times X$ eine zweite Trivialisierung von $f^*E$. Dann existiert ein Fasertransport $g:\,F_b \lto F_{b'}$, der das Diagramm
$$\xymatrix{
& f^*E \ar[ld]_T \ar[rd]^{T'}\\
F_b\times X \ar[rr]^{g \times\id_X}&& F_{b'}\times X
}$$
bis auf Faserhomotopie kommutativ macht.
\end{enumerate}
\end{prop}

\begin{proof}[Beweis]
Die Behauptung (i) wurde bereits gezeigt. Für den Beweis von (ii) bezeichne mit $T^{-1}$ ein Faserhomotopie-Inverses von $T$. Die Abbildung $T'\circ T^{-1}:\,F_b\times X\lto F_{b'}\times X$ ist eine Faserhomotopie-Äquivalenz, deren Einschränkung auf jede Faser eine Fasertransport $F_b\lto F_{b'}$ ist. Bezeichnet nun $\theta:\, F_{b'}\times X\lto F_b$ die Komposition
\[F_{b'}\times X \xrightarrow{T'\circ T^{-1}} F_b\times X \xrightarrow{\Proj} F_b,\]
und ist $H:\,X\times I\lto X$ eine Kontraktion von $X$ auf $x_0$, so ist die Abbildung
\[F_{b'}\times X\times I\lto F_b\times X,\quad (y,x,t)\mapsto (\theta(y,H(x,t)),x)\]
eine Faserhomotopie zwischen $T'\circ T^{-1}$ und $\theta(\cdot,x_0)\times\id_X$. Hierbei ist die Abbildung $\theta(\cdot,x_0):\,F_{b'}\lto F_b$ ein Fasertransport.
\end{proof}

\section{Die assoziierte Faserung}

Wir bezeichnen wie üblich mit $X^I$ den Raum aller stetigen Wege $I\lto X$ in $X$ mit der Kompakt-offen-Topologie. Alle im Kontext der assoziierten Faserung betrachteten Räume seien stets aus der Kategorie der kompakt erzeugten Räume. Dies hat die technischen Vorteile, dass die "`natürlich auftretenden"' Abbildungen im Zusammenhang mit Abbildungsräumen automatisch stetig sind und das Exponentialgesetz gültig ist. Für Details vgl.~\cite{whitehead}, Seite 17ff. Im Kontext von Faserungen gilt insbesondere der folgende Satz, der \cite{whitehead}, Seite 31 entnommen ist:
\begin{satz}
Die Abbildung
$$p_{0,1}:\, X^I \lto X \times X,\quad \gamma\mapsto (\gamma(0),\gamma(1))$$
ist eine Faserung.
\end{satz}

\begin{kor}
Sei $i\in\{0,1\}$. Die Abbildung
$$p_i:\,X^I \lto X,\quad \gamma\mapsto\gamma(i)$$
ist eine Faserung, und ein Homotopie-Hochhebungs-Problem
$$\xymatrix{
Y\times\{0\} \ar[rr]^{H_0} \arinclinv[d] & & X^I \ar[d]^{p_i}\\
Y\times I \ar@{.>}[rru]^H \ar[rr]^h & & X }
$$
kann stets so gelöst werden, dass $p_{1-i}\circ H$ stationär ist.
\end{kor}

\begin{proof}[Beweis]
Weite die Abbildung $p_{1-i}\circ H_0 :\, Y\times\{0\} \lto X$ zu einer stationären Homotopie $g$ auf $Y\times I$ aus. Für $i=0$ leistet dann eine Lösung des HHP
$$\xymatrix{
Y\times\{0\} \ar[rr]^{H_0} \arinclinv[d] & & X^I \ar[d]^{p_{0,1}}\\
Y\times I \ar[rr]^{(h,g)} & & X\times X }
$$
das Verlangte. Für den Fall $i=1$ ersetze im Diagramm den unteren horizontalen Pfeil $(h,g)$ durch $(g,h)$.
\end{proof}

Sei zur Vereinfachung der Darstellung nun $i=0$. Eine wichtige Eigenschaft der Faserung $p_0$ ist, dass sie eine Homotopie-Äquivalenz ist, wie wir im Folgenden zeigen werden. Dazu bemerken wir, dass der "`triviale Schnitt"'
$$X \lto X^I,\quad x\mapsto \big[\,t\mapsto x\,\big]$$
ein Rechtsinverses von $p_0$ ist. Die Abbildung 
$$K:\, X^I \times I \lto X^I,\quad (\gamma,t)\mapsto \bigl[s\mapsto \gamma((1-t)\cdot s)\bigr]$$
ist eine Kontraktion von $X^I$ auf $X\hookrightarrow X^I$; diese hat die zusätzliche Eigenschaft, dass $p_0 \circ K$ stationär ist. Wir werden diese zusätzliche Eigenschaft im Folgenden mit \emph{"`faserweise auf $X$ kontraktibel"'} bezeichnen.

Zu einer Abbildung $f:\, X\lto Y$ bezeichnen wir mit $X^f$ den topologischen Raum, der als Pull-back von $p_0:\, Y^I \lto Y$ mit $f$ hervorgeht:
\[\xymatrix{
X^f \ar[rr]^{\bar{f}} \ar[d]_{p'_0} & &  Y^I \ar[d]^{p_0} \\
X \ar[rr]^f & & Y
}\]
Wir setzen
\[p := p_1 \circ \bar{f}:\, X^f \lto Y\]

\begin{prop}
Die Abbildung $p$ ist eine Faserung.
\end{prop}

\begin{proof}[Beweis]
Sei ein HHP
$$\xymatrix{
A\times \{0\} \ar[rr]^{H_0} \arinclinv[d] & & X^f \ar[d]^p \\
A\times I \ar[rr]^h & & Y
}$$
gegeben, dessen Lösung wir nun mit Hilfe der universellen Eigenschaft des Pull-backs aus zwei Homotopien $A\times I \lto Y^I$ und $A\times I \lto X$ konstruieren werden. Zunächst erhält man das folgende HHP für die Faserung $p_1$: 
$$\xymatrix{
A\times \{0\} \ar[rr]^{\bar{f}\circ H_0} \arinclinv[d] & & Y^I \ar[d]^{p_1} \\
A\times I \ar[rr]^h \ar@{.>}[rru]^{H'} & & Y
}$$
Dessen Lösung $H'$ kann so gewählt werden, dass $p_0 \circ H':\, A\times I \lto Y$ stationär ist. Die zweite gesuchte Abbildung ist dann die Ausweitung von $p'_0\circ H_0:\, A\times\{0\}\lto X$ zu einer stationären Homotopie auf $A\times I$. 
\end{proof}

\begin{defn}
Die Faserung $p:\, X^f \lto Y$ heißt die zur Abbildung $f:\, X\lto Y$ \emph{assoziierte Faserung}.
\end{defn}

Die Idee der assoziierten Faserung besteht darin, eine zunächst beliebige Abbildung $f:\,X\lto Y$ durch eine Faserung zu ersetzen, ohne deren homotopietheoretischen Eigenschaften zu verändern. Dazu bemerken wir, dass der triviale Schnitt $Y \lto Y^I$ von $p_0$ auf Grund der universellen Eigenschaft des Pull-backs einen Schnitt $\lambda: X \lto X^f$ von $p'_0$ induziert. Dabei gilt:

\begin{prop}
Die Abbildung $\lambda$ ist homotopie-invers zu $p'_0:\, X^f\lto X$.
\end{prop}

\begin{proof}[Beweis]
Wir haben bereits festgestellt, dass $Y^I$ faserweise auf $Y$ kontrahierbar ist. Interpretiert man diesen Sachverhalt derart, dass die Wegefaserung $p_0:\,Y^I\lto Y$ vermöge der Projektion $p_0$ und des trivialen Schnitts faserhomotopie-äquivalent zur trivialen Faserung $\id:\,Y\lto Y$ ist, folgt die Aussage aus dem nachfolgenden Lemma \ref{fas:lemma2}.
\end{proof}

\begin{lem}\label{fas:lemma2}
Seien $p:\,E\lto B$ und $p':\,E'\lto B$ Faserungen; sei $f:\,E\lto E'$ eine Faserhomotopie-Äquivalenz über $B$. Sei $q:\,X\lto E$ eine Abbildung. Dann ist auch die "`zurückgezogene Abbildung"'
$$q^*f:\,q^*E\lto q^*E',$$
die von $f\circ\bar{q}:\,q^*E\lto E\lto E'$ und der Faserung $q^*E\lto X$ induziert wird, eine Faserhomotopie-Äquivalenz von Faserungen (über $X$). Genauer gilt: Ist $g:\,E'\lto E$ ein Faserhomotopie-Inverses von $f$, so ist $q^*g$ ein Faserhomotopie-Inverses von $q^*f$.
\end{lem}

\begin{proof}[Beweis]
Sei $H:\,E\times I \lto E$ eine Faserhomotopie zwischen  $\id_E$ und $g\circ f$. Die natürliche Transformation von Pull-back-Diagrammen
$$\xymatrix{
X\times I \ar[rr]^{q\times\id_I} \ar[d]_{\Proj} & & B\times I \ar[d]_{\Proj} & & E\times I\ar[ll]_{p\times\id_I} \ar[d]^H\\
X \ar[rr]^q & & B & &  E\ar[ll]_p
}$$
induziert eine Abbildung $q^*E\times I\lto q^*E$, die eine Faserhomotopie zwischen $\id_{q^*E}$ und $q^*g\circ q^*f$ darstellt. Eine Faserhomotopie zwischen $\id_{q^*E'}$ und $q^*f\circ q^*g$ erhält man analog.
\end{proof}

Zusammenfassend haben wir also zu jeder Abbildung $f:\,X\lto Y$ ein (offenbar kommutatives) Diagramm
$$\xymatrix{
X \ar[rr]^f \ar[rd]_\lambda^\simeq & & Y\\
 & X^f \ar[ru]_p
}$$
zugeordnet. Dies rechtfertigt die Aussage, dass homotopietheoretisch jede Abbildung "`so gut wie eine Faserung"' ist.

Nachfolgend führen wir noch zwei Propositionen an, die im Folgenden von Wichtigkeit sein werden.

\begin{prop}\label{fas:prop3}
Sei $H:\,X\times I\lto Y$ eine Homotopie zwischen $f,g:\,X\lto Y$. Dann existiert eine Faserhomotopie-Äquivalenz
$\omega_H:\, X^f \lto X^g$ von Faserungen über $Y$, und das Diagramm
$$\xymatrix{
& X \ar[ld]_{\lambda} \ar[rd]^{\lambda} \\
X^f \ar[rr]^{\omega_H} && X^g
}$$
kommutiert bis auf Homotopie.
\end{prop}

\begin{proof}[Beweis] 
Aus Satz \ref{fas:fasertransport}, Teil (iii) folgt die Existenz einer Faserhomotopie-Äquivalenz $\omega_H:\, X^f \lto X^g$ von Faserungen über $X$. Tatsächlich liefert der Beweis jedoch schon die Existenz der Faserhomotopie-Äquivalenz über $Y$, wenn man nur bei jedem auftretenden HHP für $p_0:\,Y^I\lto Y$ zusätzlich fordert, dass $p_1:\, Y^I\lto Y$ unter der Lösung stationär ist. Man stellt so zunächst fest, dass dann die Abbildungen $\omega_H:\,X^f\lto Y^f$ und $\omega_{H^-}:\,Y^f\lto X^f$ Faserabbildungen über $Y$ sind. 

Die Homotopie zwischen $\omega_{H^-}\circ\omega_H$ und $\id_{X^f}$ erhält man ja als Forsetzung einer gewissen Abbildung $X^f\times \dot{I}\lto Y^I$ auf $X^f\times I\times I$. Die fortzusetzende Abbildung besteht dabei aus den Homotopien $H'$ und $(H^-)'$, die nach dem ersten Schritt bereits stationär über $p_1$ sind, sowie einer ohnehin stationären Abbildung. Damit ist ersichtlich, dass auch die Fortsetzung auf $X^f\times I\times I$ sowohl im zweiten als auch im dritten Argument stationär über $p_1$ gewählt werden kann.

Die zweite Aussage ist eine Konsequenz der Tatsache, dass $\omega_H$ eine Faserhomotopie-Äquivalenz über $X$ ist, denn die jeweiligen Inklusionen $\lambda$ sind Homotopie-Inverse der jeweiligen Projektionen $X^f\lto X$ bzw.~$X^g\lto X$.
\end{proof}

\begin{prop}\label{fas:prop2}
Ist $f:\, X\lto Y$ bereits eine Faserung, so ist die Inklusion $\lambda:\, X\lto X^f$ eine Faserhomotopie-Äquivalenz.
\end{prop}

\begin{proof}[Beweis]
Siehe \cite{whitehead}, Seite 43. Man beachte, dass $p_0$ natrlich kein Faserhomotopie-Inverses sein kann, da das Diagramm
$$\xymatrix{
X \ar[rr]^f & & Y\\
& X^f \ar[ul]^{p_0} \ar[ur]_p
}$$
im Allgemeinen \emph{nicht} kommutativ ist. Allerdings kann man nun die Faserungseigenschaft von $f$ benutzen, um "`die Punkte in $X$ entlang des Weges in $Y$ zu transportieren"'. Dazu betrachte man das HHP
$$\xymatrix{
X^f\times\{0\} \ar[rr]^{p'_0} \arinclinv[d] & & X \ar[d]^f\\
X^f \times I \ar[r]_{\bar{f}\times\id_I} \ar@{.>}[rru]^H & Y^I\times I \ar[r] & Y
}$$
wobei der rechte untere horizontale Pfeil die "`Auswertungsabbildung"' $(\gamma,t)\mapsto\gamma(t)$ ist. Die Auswertung von $H$ an $X^f\times\{1\}$ liefert dann ein Faserhomotopie-Inverses von $\lambda$.
\end{proof}


\chapter{Einfache Faserungen und das Element $\theta(f)$}\label{theta}

In diesem Kapitel soll einer Abbildung $f:\,M\lto B$ topologischen Räumen unter gewissen Zusatzvoraussetzungen ein Element $\theta(f)\in H^1(B;\Whp{M})$ zugeordnet werden. Dieses Element ist über Fasertransporte der assoziierten Faserung entlang geschlossener Wege im Basisraum definiert und hat folgende Eigenschaften:

\begin{prop}\label{theta:prop_haupt}
Das Element $\theta(f)$ hängt nur von der Homotopieklasse von $f$ in $[M,B]$ ab. Ist $f$ ein differenzierbares Faserbündel von zusammenhängenden und geschlossenen differenzierbaren Mannigfaltigkeiten, so ist $\theta(f)=0$.
\end{prop}

Genauer misst das Element $\theta(f)$ die Whitehead-Torsion der entsprechenden Fasertransporte. Allerdings ist die Whitehead-Torsion zunächst nur für zelluläre Abbildungen zwischen endlichen CW-Komplexen erklärt, während die Fasern von assoziierten Faserungen im Allgemeinen keine CW-Komplexe sind. Daher ist es zunächst nötig, eine Torsion von Abbildungen zwischen Räumen zu erklären, die nur vom Homotopietyp eines endlichen CW-Komplexes sind. Dies gelingt durch die Wahl einer "`einfachen Struktur"' auf diesen Räumen und wird im ersten Teil dieses Kapitels besprochen. Sodann betrachten wir Faserungen $p:\,E\lto B$, deren Faser $F$ den Homotopietyp eines endlichen CW-Komplexes haben, und definieren das Element $\theta(p)\in H^1(B;\Whp{E})$ als Torsion von Fasertransporten nach Wahl einer einfachen Struktur auf $F$. Die Konstruktion der assoziierten Faserung erlaubt es schließlich, diese Definition auf beliebige Abbildungen $f:\,M\lto B$ auszudehnen.

\section{Die Whitehead-Torsion}\label{theta:abschn1}

Ziel dieses Abschnittes ist es, die wichtigsten Eigenschaften der Whitehead-Torsion, die für uns von Belang sein werden, zusammenzufassen. In der Literatur gibt es viele ausführliche Einführungen in das Thema (siehe etwa \cite{milnor:whitehead_torsion}, \cite{cohen} oder \cite{lueck:surgery}). Die topologische Invarianz der Whitehead-Torsion wurde von Chapman in \cite{chapman} bewiesen.

Für eine Gruppe $G$ ist die Whitehead-Gruppe $\Wh(G)$ definiert als der Quotient von $K_1(\ZZ[G])$ durch die Untergruppe, die von den Elementen $[\pm g]$ mit $g\in G$ erzeugt wird. Hierbei bezeichnet $[\pm g]$ die Äquivalenzklasse in $K_1(\ZZ[G])$ der $(1\times 1)$-Matrix mit Eintrag $\pm g\in\ZZ[G]$. Ein Gruppenhomomorphismus $f:\,G\lto H$ induziert vermittels der Funktorialität von $K_1$ einen Gruppenhomomorphismus $\Wh(G)\lto \Wh(H)$, wodurch auch $\Wh$ zu einem Funktor wird.

Für einen Raum $Y$ sei
\[\Whp{Y}= \bigoplus_{C\in\pi_0(Y)} \Wh(\pi_1(C)).\]
Auch die Abbildung $Y\mapsto\Whp{Y}$ ist offenbar funktoriell. Die Whitehead-Torsion ist nun eine Zuordnung, die jeder Homotopie-Äquivalenz $f:\,X\lto Y$ von endlichen CW-Komplexen ein Element $\tau(f)\in\Whp{Y}$ zuordnet. Die wesentlichen Eigenschaften der Whitehead-Torsion sind in folgendem Satz zusammengefasst:

\begin{satz}\label{theta:whitehead_torsion}
Die Whitehead-Torsion $\tau$ hat die folgenden Eigenschaften:
\begin{enumerate}
\item \textup{\textbf{(Homotopie-Invarianz)}} Sind $f,g:\,X\lto Y$ homotope Homotopie-Äquivalenzen zwischen endlichen CW-Komplexen, so ist $\tau(f)=\tau(g)$.
\item \textup{\textbf{(Kettenregel)}} Sind $f:\, X\lto Y$ und $g:\,Y\lto Z$ Homotopie-Äquivalenzen zwischen endlichen CW-Komplexen, so gilt
\[\tau(g\circ f) = \tau(g) + g_*\, \tau(f).\]
\item \textup{\textbf{(Summenformel)}} Seien im folgenden kommutativen Diagramm
\[\xymatrix{
& Y_0 \ar[rrr]^{j_2} \arinclinv'[d]^(0.8){j_1}[dd] & & & Y_2 \ar[dd]^{k_2}\\
X_0 \ar[ur]^{f_0} \ar[rrr]^(0.6){i_2} \arinclinv[dd]^{i_1} & & & X_2 \ar[ru]^(0.4){f_2} \ar[dd]\\
& Y_1 \ar'[rr]^{k_1}[rrr] & & & Y\\
X_1 \ar[ru]^{f_1} \ar[rrr] & & & X \ar[ru]^f
}\]
das vordere und das hintere Quadrat zelluläre Push-outs von endlichen CW-Kom\-plexen. Seien $f_0,f_1,f_2$ Homotopie-Äquivalenzen. Dann ist auch $f$ eine Homotopie-Äquivalenz, und es gilt:
\[\tau(f) = k_{1*}\,\tau(f_1) + k_{2*}\,\tau(f_2) - k_{0*}\,\tau(f_0) \in \Whp{Y}\]
mit $k_0:=k_1\circ j_1$.
\item \textup{\textbf{(Produktformel)}} Seien $f_n:\,X_n\lto Y_n$ für $n\in\{1,2\}$ Homotopie-Äquivalenzen zwischen zusammenhängenden endlichen CW-Komplexen. Dann gilt in $\Whp{(Y_1\times Y_2)}$:
\[\tau(f_1\times f_2) = \chi(Y_1)\cdot i_{2*}\,\tau(f_2) + \chi(Y_2)\cdot i_{1*}\,\tau(f_1),\]
wobei die Abbildungen $i_{1*}$ und $i_{2*}$ von Inklusionen $Y_1\times\{*\}\lto Y_1\times Y_2$ und $\{*\}\times Y_2\lto Y_1\times Y_2$ bezüglich beliebiger Basispunkte induziert werden, und $\chi$ die Eulercharakteristik bezeichnet.
\item \textup{\textbf{(Topologische Invarianz)}} Die Torsion eines Homöomorphismus von endlichen CW-Komplexen verschwindet.
\item \textup{\textbf{(Hindernis-Eigenschaft)}} Eine Homotopie-Äquivalenz von endlichen CW-Komple\-xen hat genau dann verschwindende Torsion, wenn sie homotop zu einer Komposition von elementaren Expansionen und Retrakten ist.
\end{enumerate}
\end{satz}

Dabei nennen wir ein Push-out-Diagramm
\[\xymatrix{
Y_0 \ar[rr]^{\phi} \arinclinv[d] && Y_2 \arinclinv[d] \\
Y_1 \ar[rr]^{\Phi} && Y
}\]
\emph{zelluläres Push-out,} falls $(Y_1,Y_0)$ ein CW-Paar und $Y_2$ ein CW-Komplex ist, sodass die Abbildung $\phi$ zellulär ist. In diesem Fall besitzt auch $Y$ eine CW-Struktur, bezüglich derer die offenen Zellen von $Y_1 - Y_0$ mittels $\Phi$ auf die offenen Zellen von $Y - Y_2$ abgebildet werden und die jeweiligen anklebenden Abbildungen für $Y$ aus den anklebenden Abbildungen für $Y_1$ hervorgehen, indem man sie mittels $\Phi$ weitertransportiert.

Eine \emph{elementare Expansion} $f:\,X\lto Y$ ist eine Inklusion von endlichen CW-Kom\-plexen, wenn $Y$ aus $X$ durch ein Push-out der Form
\[\xymatrix{
S^n_+ \ar[rr]^{\varphi} \arinclinv[d] && X \arinclinv[d]\\
D^{n+1} \ar[rr] && Y
}\]
hervorgeht, in dem $S^n_+$ die obere Hemisphäre der $n$-Sphäre $S_n\subset D^{n+1}$ ist, und die anklebende Abbildung $\varphi$ in das $n$-Gerüst von $X$, den Äquator $S^{n-1}\subset S^n_+$ sogar in das $(n-1)$-Gerüst von $X$ abbildet. Man versieht dann das Paar $(D^{n+1}, S^n_+)$ mit der CW-Struktur, die eine 0-Zelle, den Äquator $S^{n-1}$ als $(n-1)$-Zelle, die beiden Hemisphären $S^n_\pm$ jeweils als $n$-Zellen und $D^{n+1}$ als $(n+1)$-Zelle hat, und erkennt, dass $Y$ aus $X$ durch Ankleben einer $n$ und einer $(n+1)$-Zelle hervorgeht. Man sieht leicht ein, dass es sich hierbei um eine Homotopie-Äquivalenz handelt. --- Ein Homotopie-Inverses einer elementaren Expansion heißt \emph{elementarer Retrakt}. 

Eine Homotopie-Äquivalenz $f:\,X\lto Y$ von endlichen CW-Komplexen heißt \emph{einfach,} wenn es eine Folge $X_0=X, X_1, \ldots, X_n=Y$ von endlichen CW-Komplexen zusammen mit Homotopie-Äquivalenzen $f_i:\,X_i\lto X_{i+1}$ für $i\in\{0,\ldots n-1\}$ gibt, sodass jedes $f_i$ eine elementare Expansion oder ein elementarer Retrakt ist und $f$ homotop zur Komposition $f_{n-1}\circ\ldots\circ f_0$ ist. Die Hindernis-Eigenschaft besagt also, dass $f$ genau dann verschwindende Torsion hat, wenn $f$ einfach ist. Man kann zeigen, dass für eine Abbildung $f:\,X\lto Y$ die Inklusion von $Y$ in den Abbildungszylinder $\Cyl(f)$ eine Komposition von elementaren Expansionen ist; somit verschwindet ihre Whitehead-Torsion.

\section{Einfache Strukturen}

Dieses Unterkapitel wird sich mit einfache Strukturen auf Räumen vom Homotopietyp eines endlichen CW-Komplexes beschäftigen. Die Darstellung geschieht im Sinne von \cite{LNM1408}, Seite 74ff. Jedoch verallgemeinern wir die dortige Situation etwas, indem wir den Begriff einer einfachen Struktur modulo $N$ für eine Untergruppe $N$ in $\Whp{Y}/N$ definieren; eine einfache Struktur ist dann eine einfache Struktur modulo der trivialen Untergruppe. Es wird sich zeigen, dass Satz \ref{theta:whitehead_torsion} eine natürliche Verallgemeinerung für Räume mit einfacher Struktur hat.

Zu einem Raum $X$ vom Homotopietyp eines endlichen CW-Komplexes betrachte alle Paare $(A,f)$ von CW-Komplexen $A$ und Homotopie-Äquivalenzen $f:\,A\lto X$. Zwei Paare $(A,f)$ und $(B,g)$ nennen wir \emph{äquivalent}, falls die Whitehead-Torsion $\tau(g^{-1}\circ f)\in \Whp{B}$ verschwindet. Für eine Untergruppe $N\subset\Whp{X}$ nennen wir $(A,f)$ und $(B,g)$ \emph{äquivalent modulo $N$}, falls gilt: $g_*\, \tau(g^{-1}\circ f)\in N$, wobei $g_*:\,\Whp{B}\lto \Whp{X}$ den von $g$ induzierten Gruppenisomorphismus bezeichnet. Man beachte, dass zwei Paare $(A,f)$ und $(B,g)$ genau dann äquivalent sind, wenn sie modulo der trivialen Untergruppe $\{0\}\subset\Whp{X}$ äquivalent sind. Die Äquivalenz modulo $N$ ist offenbar tatsächlich eine Äquivalenzrelation. 

\begin{defn}\begin{enumerate}\smallsep
\item Eine Äquivalenzklasse (bzw.~eine Äquivalenzklasse modulo $N$) von solchen Paaren $(A,f)$ heißt \emph{einfache Struktur} (bzw.~\emph{einfache Struktur modulo $N$}) auf $Z$. Ist $\xi$ eine einfache Struktur (modulo $N$) auf $X$, so heißt das Paar $(X,\xi)$ \emph{Raum mit einfacher Struktur (modulo $N$)}. 

\item Ist $\varphi:\,(X,\xi)\lto (Y,\eta)$ eine Homotopie-Äquivalenz von Räumen mit einfacher Struktur, so setzen wir
\[\tau(\varphi) = g_*\, \tau(g^{-1}\circ\varphi\circ f)\in \Whp{Y}\]
für Vertreter $f:\,A\lto X$ bzw.~$g:\,B\lto Y$ von $\xi$ bzw.~$\eta$. Die Homotopie-Äquivalenz $\varphi$ heißt \emph{einfach,} falls $\tau(\varphi)=0$. 
\end{enumerate}
\end{defn}

Ist keine Verwechslung möglich, werden wird oft statt $(Z,\xi)$ nur $Z$ schreiben. Die Definition von $\tau(\varphi)$ ist unabhängig von der Wahl der Vertreter $f$ und $g$ nach der Kettenregel. Für eine Untergruppe $N'\subset N$ induziert eine einfache Struktur modulo $N'$ durch die Wahl eines Vertreters der Äquivalenzklasse modulo $N'$ auf kanonische Weise eine einfache Struktur modulo $N$. 

Zu einer gegebenen Homotopie-Äquivalenz $\varphi:\,X\lto Y$ und einer Untergruppe $N\subset\Whp{Y}$ werden wir im Folgenden stets mittels des induzierten Gruppenisomorphismus $\varphi_*:\,\Whp{X}\lto\Whp{Y}$ die Gruppe $N$ als Untergruppe von $\Whp{X}$ auffassen. 

\begin{defn}
Ist $\varphi:\,(X,\xi)\lto (Y,\eta)$ eine Homotopie-Äquivalenz von Räumen mit einfacher Struktur modulo $N$, so heißt $\varphi$ \emph{einfach modulo $N$,} falls gilt: Für eine (und damit alle) Wahlen von einfachen Strukturen $\xi'$ und $\eta'$ auf $X$ bzw.~$Y$, deren induzierte einfache Strukturen modulo $N$ gerade $\xi$ bzw.~$\eta$ sind, ist 
\[\tau\big[\varphi:\,(X,\xi')\lto (Y,\eta')\bigr]\in N.\]
\end{defn}

Ist $\xi$ eine einfache Struktur modulo $N$ von $X$, so bezeichnen wir mit $\varphi_*(\xi)$ die einfache Struktur modulo $N$ von $Y$, bezüglich der $\varphi:\,(X,\xi)\lto (Y,f(\xi))$ einfach modulo $N$ ist. Mit anderen Worten gilt: $(A,\varphi\circ f)$ ist ein Vertreter von $\varphi_*(\xi)$, wenn $(A,f)$ ein Vertreter von $\xi$ ist.

Es sei nun
\begin{equation}\begin{split}\label{theta:push-out}
\xymatrix{
X_0 \ar[rr]^{j_2} \arinclinv[d]_{j_1} & & X_2 \arinclinv[d] \\
X_1 \ar[rr] & & X
}
\end{split}
\end{equation}
ein Push-out von topologischen Räumen, wobei $j_1$ eine Kofaserung sei. Weiter seien die Räume $X_n$ für $n\in\{0,1,2\}$ mit einfachen Strukturen $\xi_n$ ausgestattet. Dann stattet die folgende Konstruktion auch den Raum $X$ auf kanonische Weise mit einer einfachen Struktur $\xi$ aus.

Wähle ein kommutatives Diagramm
\begin{equation}\begin{split}\label{theta:diag1}
\xymatrix{
A_1 \ar[d]^{f_1} & & A_0 \arinclinv[ll]_{i_1} \ar[d]^{f_0} \ar[rr]^{i_2} & & A_2 \ar[d]^{f_2} \\
X_1  & & X_0 \arinclinv[ll]_{j_1} \ar[rr]^{j_2} & & X_2
}
\end{split}
\end{equation}
sodass für $i\in\{0,1,2\}$ das Paar $(A_i,f_i)$ die einfache Struktur $\xi_i$ repräsentiert und die obere Zeile ein zelluläres Push-out-Diagramm von endlichen CW-Komplexen ist. In dieser Situation ist das Push-out $A$ der oberen Zeile auf kanonische Weise ein endlicher CW-Komplex. Das folgende Lemma zeigt, dass die induzierte Abbildung $f:\,A\lto X$ eine einfache Struktur $\xi$ auf $X$ definiert.

\begin{lem}\label{theta:lemma10}
Sind in einem Diagramm der Form \eqref{theta:diag1} von topologischen Räumen alle vertikalen Abbildungen Homotopie-Äquivalenzen und die Abbildungen $i_1$ und $j_1$ Kofaserungen, so ist auch die auf den Push-outs der Zeilen induzierte Abbildung eine Homotopie-Äquivalenz.
\end{lem}

\begin{proof}[Beweis]
Es handelt sich um ein Resultat aus der Theorie der Kofaserungen (siehe etwa \cite{tomdieck}, Seite 192ff.) Hier soll der Beweis nur skizziert werden. Für Kofaserungen existiert ein duales Analogon von Satz \ref{fas:fasertransport} und entsprechend eine duale Version der Proposition \ref{fas:prop1}. Deren Anwendung erlaubt es, die Aussage zurückzuführen auf den Fall, in dem $(A_2,i_2)$ und $(X_2,j_2)$ Inklusionen in Abbildungszylinder sind. In diesem Fall sind jedoch auch $i_2$ und $j_2$ Kofaserungen, und die Behauptung folgt aus der folgenden Proposition \ref{theta:prop3}. 
\end{proof}

\begin{prop}[\cite{tomdieck}, Seite 196]\label{theta:prop3}
Seien im kommutativen Diagramm
\[\xymatrix{
A \arincl[rr]^i \ar[d]^f && X \ar[d]^F\\
B \arincl[rr]^j && Y
}\]
die Abbildungen $i$, $j$ Kofaserungen und die Abbildungen $f$, $F$ Homotopie-Äquivalenzen. Dann ist $(F,f)$ eine Homotopie-Äquivalenz von Paaren. Genauer gilt: Jedes Homotopie-Inverse $g$ von $f$ erweitert sich zu einem Homotopie-Inversen $G$ von $F$.
\end{prop}

Insbesondere gilt also unter den Voraussetzungen von Lemma \ref{theta:lemma10} im Diagramm \eqref{theta:diag1}: Falls auch die Abbildungen $i_2$ und $j_2$ Kofaserungen sind, so existiert ein Homotopie-Inverses von $f$, das von geeigneten Homotopie-Inversen der Abbildungen $f_i$ induziert wird. Das Homotopie-Inverse von $f_0$ ist dabei sogar frei wählbar.

Wir wissen also, dass auch ein Push-out von Räumen mit einfacher Struktur wieder vom Homotopietyp eines endlichen CW-Komplexes ist, falls eine Kofaserungsbedingung erfüllt ist. Nachfolgend zeigen wir in Lemma \ref{theta:lemma4} die Eindeutigkeit der induzierten einfachen Struktur $\xi$; aus der zweifachen Anwendung des danach folgenden Lemmas \ref{theta:lemma3} geht hervor, dass ein Diagramm \eqref{theta:diag1} mit den geforderten Eigenschaften existiert. 

Man beachte, dass auf diese Weise auch die disjunkte Vereinigung $(X,\xi)\amalg (Y,\eta)$ von Räumen mit einfacher Struktur mit einer einfachen Struktur $\xi\amalg\eta$ ausgestattet wird. Diese stimmt offenbar mit der folgenden Konstruktion überein: Zu Vertretern $(A,\alpha)$ bzw.~$(B,\beta)$ von $\xi$ bzw.~$\eta$ definiere $\xi\amalg\eta$ als die Äquivalenzklasse von $(A\amalg B,\alpha\amalg\beta)$. Wir werden im Folgenden stets implizit voraussetzen, dass derartige disjunkte Vereinigungen diese einfache Struktur tragen. 

\begin{defn}
Seien $X_0,X_1,X_2,X$ Räume mit einfacher Struktur. Ein Push-out-Diagramm \eqref{theta:push-out} heißt \emph{einfach} (bzw.~\emph{einfach modulo $N$}), falls $j_1$ eine Kofaserung ist und die durch das Diagramm \eqref{theta:diag1} auf $X$ induzierte Abbildung $f:\,A\lto X$ einfach (bzw.~einfach modulo $N$) ist, mit anderen Worten also die auf $X$ induzierte einfache Struktur (bzw.~einfache Struktur modulo $N$) mit der gegebenen übereinstimmt. 
\end{defn}

Ein endlicher CW-Komplex $X$ hat eine, durch die das Paar $(X,\id_X)$ gegebene, natürliche einfache Struktur, die wir mit $[\id]$ bezeichnen werden. Aus der obigen Konstruktion geht hervor, dass ein zelluläres Push-out von endlichen CW-Komplexen einfach ist, wenn alle Komplexe ihre natürliche einfache Struktur tragen.

\begin{lem}\label{theta:lemma4}
Seien im Diagramm
\begin{equation}\begin{split}\label{theta:diag2}
\xymatrix{
X_1 \ar[d]^{f_1} && X_0 \arinclinv[ll] \ar[d]^{f_0} \ar[rr]^{i_2} && X_2 \ar[d]^{f_2} \\
Z_1 && Z_0 \arinclinv[ll]_{j_1} \ar[rr]^{j_2} && Z_2 \\
Y_1 \ar[u]_{g_1} && Y_0 \arinclinv[ll] \ar[u]_{g_0} \ar[rr]^{k_2} && Y_2 \ar[u]_{g_2}
}\end{split}\end{equation}
die oberste und unterste Zeile zelluläre Push-out-Diagramme von endlichen CW-Kom\-plexen, $j_1$ eine Kofaserung und alle vertikalen Abbildungen Homotopie-Äquivalenzen. Bezeichne mit $X$,$Y$,$Z$ die Push-outs der jeweiligen Zeilen sowie mit $f:\,X\lto Z$ und $g:\,Y\lto Z$ die von den $f_i$ und $g_i$ induzierten Abbildungen. Dann gilt in $\Whp{Z}$:
\[g_*\,\tau(g^{-1}\circ f) = l_{1*}\,g_{1*}\,\tau(g_1^{-1}\circ f_1) + l_{2*}\,g_{2*}\,\tau(g_2^{-1}\circ f_2) - l_{0*}\,g_{0*}\,\tau(g_0^{-1}\circ f_0),\]
wenn $l_i:\,Z_i\lto Z$ für $i\in\{0,1,2\}$ die Strukturabbildungen des Push-outs sind.
\end{lem}

\begin{proof}[Beweis]
Im Spezialfall, dass $j_2$ und $k_2$ Kofaserungen sind, erhält man unter Anwendung von Proposition \ref{theta:prop3}, dass geeignete Homotopie-Inverse der Abbildungen $g_i$ ein Homotopie-Inverses der Abbildung $g$ induzieren. Die Behauptung folgt dann direkt aus der Summenformel für die Whitehead-Torsion. 

Der allgemeine Fall kann nun auf den Spezialfall zurückgeführt werden. Diagramm \eqref{theta:diag2} induziert nämlich ein Diagramm, in dem $(X_2,i_2)$, $(Z_2,j_2)$ und $(Y_2,k_2)$ durch die Inklusionen in die jeweiligen Abbildungszylinder $\Cyl(i_2)$, $\Cyl(j_2)$ und $\Cyl(k_2)$ ersetzt werden. Diese Inklusionen sind Kofaserungen, und die Räume $\Cyl(i_2)$ und $\Cyl(k_2)$ sind endliche CW-Komplexe, sodass wir uns im oben betrachteten Spezialfall befinden. Es bezeichne $p:\,\Cyl(i_2)\lto X_2$ die Projektion des Abbildungszylinders. Diese ist einfach; nach der Summenformel ist dann auch die vom Diagramm
\[\xymatrix{
X_1 \ar@{=}[d] && X_0 \arinclinv[ll] \ar@{=}[d] \arincl[rr] && {\Cyl(i_2)} \ar[d]^p\\
X_1 && X_0 \arinclinv[ll] \ar[rr]^{i_2} && X_2 
}\]
induzierte Homotopie-Äquivalenz einfach; dasselbe gilt für das entsprechende Diagramm mit den Räumen $Y_i$. Nach der Kettenregel folgt damit die Behauptung.
\end{proof}

\begin{lem}\label{theta:lemma3}
Zu einer Abbildung $f:\,(X,\xi)\lto (Y,\eta)$ von Räumen mit einfacher Struktur und einem Vertreter $(A,\alpha)$ von $\xi$ existiert stets ein Vertreter $(B,\beta)$ von $\eta$, sodass $A$ ein Unter-CW-Komplex von $B$ ist und das folgende Diagramm kommutiert:
\[\xymatrix{
A \ar[d]^{\alpha} \arincl[rr] & & B \ar[d]^{\beta} \\
X \ar[rr]^{f} & & Y
}\]
\end{lem}

\begin{proof}[Beweis]
Wähle einen Repräsentanten $(B',\beta')$ von $\eta$, und sei $g:\,A\lto B'$ eine Abbildung mit $\beta'\circ g \simeq f\circ\alpha$. Auf Grund des zellulären Approximationssatzes können wir dabei $g$ als zellulär annehmen. Die Abbildung $g$ faktorisiert über den Abbildungszylinder in einer Komposition $p\circ i:\,A\lto \Cyl(g)\lto B$. Hierbei ist $\Cyl(g)$ ein endlicher CW-Komplex, und $i$ stellt die Inklusion eines Unter-CW-Komplexes dar. Das Diagramm
\[\xymatrix{
A \ar[d]^{\alpha} \arincl[rr]^{i} & & {\Cyl(g)} \ar[d]^{\beta'\circ p} \\
X  \ar[rr]^{f} & & Y
}\]
kommutiert nun bis auf Homotopie; die Kofaserungseigenschaft der Abbildung $i$ erlaubt jedoch eine homotope Abänderung der Abbildung $\beta'\circ p$ in eine Abbildung $\beta$, die das Diagramm echt kommutativ macht. Die Wahl von $B:=\Cyl(g)$ leistet dann das Verlangte, denn die Projektion $p$ ist einfach, und daher repräsentiert $(\Cyl(g),\beta)$ die einfache Struktur $\eta$.
\end{proof}

Schließlich vermerken wir: Sind $(X,\xi)$ und $(Y,\eta)$ Räume mit einfacher Struktur, so erhält $X\times Y$ eine kanonische einfache Struktur $\xi\times\eta$ auf die folgende Weise: Wähle Repräsentanten $f:\,A\lto X$ und $g:\,B\lto Y$ von $\xi$ und $\eta$ und setze $\xi\times\eta := [(A\times B, f\times g)]$. Für wegzusammenhängende Räume $X$ und $Y$ folgt dabei die Wohldefiniertheit direkt aus der Produktformel. Um die Wohldefiniertheit im allgemeinen Fall einzusehen, seien $(A,a)$ und $(A',a')$ bzw.~$(B,b)$ und $(B',b')$ jeweils zwei endliche CW-Modelle von $X$ bzw.~$Y$. Ein endliches CW-Modell $(A,a)$ induziert jedoch endliche CW-Modelle $(A_i,a_i)$ jeder Wegzusammenhangskomponente. Da jeder CW-Komplex die disjunkte Vereinigung seiner Wegzusammenhangskomponenten ist, berechnet sich also die Torsion des Wechsels mit Hilfe der Summenformel wie folgt:
\begin{equation*}
\begin{split}
\tau((a\times b)^{-1}\circ (a'\times b')) & = \tau((a^{-1}\circ a')\times (b^{-1}\circ b'))\\
 & = \tau\biggl[\biggl(\coprod_{i\in \pi_0(X)} a_i^{-1}\circ a'_i\biggr) \times \biggl(\coprod_{j\in\pi_0(Y)} b_j^{-1}\circ b'_j\biggr)\biggr]\\
 & = 0.\end{split}
\end{equation*}

Im Kontext von Räumen mit einfacher Struktur hat nun Satz \ref{theta:whitehead_torsion} die folgende Verallgemeinerung:

\begin{satz} Die Whitehead-Torsion zwischen Räumen mit einfacher Struktur hat die folgenden Eigenschaften:
\begin{enumerate}\setlength{\itemsep}{0pt}
\item \textup{\textbf{(Homotopie-Invarianz)}} Sind $f,g:\,X\lto Y$ homotope Homotopie-Äquivalenzen zwischen Räumen mit einfacher Struktur, so ist $\tau(f)=\tau(g)$.
\item \textup{\textbf{(Kettenregel)}} Sind $f:\, X\lto Y$ und $g:\,Y\lto Z$ Homotopie-Äquivalenzen zwischen Räumen mit einfacher Struktur, so gilt
\[\tau(g\circ f) = \tau(g) + g_*\, \tau(f).\]
\item \textup{\textbf{(Summenformel)}} Seien im folgenden kommutativen Diagramm
\[\xymatrix{
& Y_0 \ar[rrr]^{j_2} \arinclinv'[d]^(0.8){j_1}[dd] & & & Y_2 \ar[dd]^{k_2}\\
X_0 \ar[ur]^{f_0} \ar[rrr]^(0.6){i_2} \arinclinv[dd]^{i_1} & & & X_2 \ar[ru]^(0.4){f_2} \ar[dd]\\
& Y_1 \ar'[rr]^{k_1}[rrr] & & & Y\\
X_1 \ar[ru]^{f_1} \ar[rrr] & & & X \ar[ru]^f
}\]
alle Räume mit einfacher Struktur, alle vertikalen Abbildungen Kofaserungen und die Abbildungen $f_0,f_1,f_2$ Homotopie-Äquivalenzen. Seien das vordere und das hintere Quadrat einfache Push-outs modulo $N$ für eine Untergruppe $N\subset\Whp{Y}$. Dann gilt:
\[\tau(f) \equiv k_{1*}\,\tau(f_1) + k_{2*}\,\tau(f_2) - k_{0*}\,\tau(f_0) \mod N\]
mit $k_0:=k_1\circ j_1$.
\item \textup{\textbf{(Produktformel)}} Seien $f_n:\,X_n\lto Y_n$ für $n\in\{1,2\}$ Homotopie-Äquivalenzen zwischen wegzusammenhängenden Räumen mit einfacher Struktur. Dann gilt bezüg\-lich der kanonischen einfachen Struktur auf den jeweiligen Produkträumen:
\[\tau(f_1\times f_2) = \chi(Y_1)\cdot i_{2*}\,\tau(f_2) + \chi(Y_2)\cdot i_{1*}\,\tau(f_1),\]
wobei die Abbildungen $i_{1*}$ und $i_{2*}$ von Inklusionen $Y_1\times\{*\}\lto Y_1\times Y_2$ und $\{*\}\times Y_2\lto Y_1\times Y_2$ bezüglich beliebiger Basispunkte induziert werden, und $\chi$ die Eulercharakteristik bezeichnet.
\end{enumerate}
\end{satz}

\begin{proof}[Beweis]
Seien $(A,\alpha)$,$(B,\beta)$,$(C,\gamma)$ Vertreter der einfachen Strukturen von $X$,$Y$,$Z$.

(i) Es gilt
\[\tau(f)=\beta_*\,\tau(\beta^{-1}\circ f\circ\alpha)=\beta_*\,\tau(\beta^{-1}\circ g\circ\alpha)=\tau(g).\]

(ii) Dies zeigt die Rechnung
\begin{equation*}\begin{split}
\tau(g\circ f)\; & = \;\gamma_*\,\tau((\gamma^{-1}\circ g \circ\beta)\circ(\beta^{-1}\circ f\circ \alpha)) \\
& =\; \gamma_*\,\tau(\gamma^{-1}\circ g\circ\beta) + g_*\,\beta_*\,\tau(\beta^{-1}\circ f\circ \alpha)\\
& =\; \tau(g) + g_*\,\tau(f).
\end{split}\end{equation*}

(iii) Wähle wie in Diagramm \eqref{theta:diag1} Vertreter $(A_n,\alpha_n)$ bzw.~$(B_n,\beta_n)$ der einfachen Strukturen von $X_n$ bzw.~$Y_n$ mit zugehörigen kommutativen Diagrammen und bezeichne mit $(A,\alpha)$ bzw.~$(B,\beta)$ die jeweiligen Push-outs. Nach Definition sind $\alpha:\,A\lto X$ und $\beta:\,B\lto X$ einfach modulo $N$; es folgt mit Teil (ii):
\[\tau(f)\equiv \beta_*\,\tau(\beta^{-1}\circ f\circ\alpha)\mod N.\]
Aus Lemma \ref{theta:lemma4} folgt schließlich:
\begin{equation*}\begin{split}
& \beta_* \,\tau(\beta^{-1}\circ f\circ \alpha) \\
 = \; &k_{1*}\,\beta_{1*}\, \tau(\beta_1^{-1}\circ f_1\circ \alpha_1) + k_{2*}\,\beta_{2*}\, \tau(\beta_2^{-1}\circ f_2\circ \alpha_2) - k_{0*}\,\beta_{0*}\, \tau(\beta_0^{-1}\circ f_0\circ \alpha_0)\\
 = \; &  k_{1*}\,\tau(f_1) + k_{2*}\,\tau(f_2) - k_{0*}\,\tau(f_0).
\end{split}
\end{equation*}

(iv) Seien $(A_n,\alpha_n)$ und $(B_n,\beta_n)$ Vertreter der einfachen Strukturen von $X_n$ und $Y_n$. Dann sind $A_n$ und $B_n$ endliche zusammenhängende CW-Komplexe. Sei für $n\in\{1,2\}$ die Abbildung $j_{n*}:\,\Whp{B_n}\lto \Whp{(B_1\times B_2)}$ induziert von Inklusionen $\{*\}\times B_2\lto B_1\times B_2$ bzw.~$B_1\times\{*\}\lto B_1\times B_2$. Dann gilt:
\begin{equation*}\begin{split}
\tau(f_1\times f_2)  =\; &(\beta_1\times\beta_2)_*\, \tau((\beta_1\times\beta_2)^{-1}\circ (f_1\times f_2) \circ (\alpha_1\times\alpha_2)) \\
  = \; &(\beta_1\times\beta_2)_*\, \tau((\beta_1^{-1}\circ f_1\circ\alpha_1)\times(\beta_2^{-1}\circ f_2\circ\alpha_2))\\
  = \; &\chi(X_1)\cdot (\beta_1\times\beta_2)_*\, j_{2*}\, \tau(\beta_2^{-1}\circ f_2\circ\alpha_2) \\
&+ \chi(X_2)\cdot (\beta_1\times\beta_2)_*\, j_{1*}\, \tau(\beta_1^{-1}\circ f_1\circ\alpha_1)\\
  = \; &\chi(X_1)\cdot i_{2*}\, \beta_{2*}\, \tau(\beta_2^{-1}\circ f_2\circ\alpha_2) + \chi(X_2)\cdot i_{1*}\, \beta_{1*}\, \tau(\beta_1^{-1}\circ f_1\circ\alpha_1) \\
  = \; &\chi(X_1)\cdot i_{2*}\,\tau(f_2) + \chi(X_2)\cdot i_{1*}\,\tau(f_1).\qedhere
\end{split}\end{equation*}
\end{proof}

Im Folgenden stellen wir noch zwei technische Lemmata bereit, die das Verhalten von einfachen Strukturen modulo $N$ unter mehrfacher Anwendung der Push-out-Konstruk\-tion beschreiben. Diese werden in Kapitel \ref{tau} Anwendung finden.

\begin{lem}["`Transitivität"'] \label{theta:transitivitaet}
Seien $X_i$ für $i\in\{1,\dots,6\}$ Räume mit einfacher Struktur und $N\subset\Whp{X_6}$ eine Untergruppe. Sei im Diagramm
\[\xymatrix{
X_1 \ar[rr]^{k_1} \ar[d]^{j_1} && X_3 \ar[rr]^{k_3} \ar[d]^{j_3} && X_5 \ar[d]^{j_5} \\
X_2 \ar[rr]^{k_2}  && X_4 \ar[rr]^{k_4}  && X_6 
}\]
eine der beiden folgenden Bedingungen erfüllt: \textup{(}i\textup{)} Die Abbildungen $k_1$ und $k_3$ sind Kofaserungen, oder \textup{(}ii\textup{)} die Abbildung $j_1$ eine Kofaserung. Ist das rechte Quadrat ein einfaches Push-out modulo $N$ und das linke Quadrat ein einfaches Push-out modulo $k_{4*}^{-1}(N)$, so ist auch das äußere Rechteck ein einfaches Push-out modulo $N$.
\end{lem}

\begin{proof}[Beweis]
Lemma \ref{theta:lemma3} garantiert die Existenz von Vertretern $(A_n,\alpha_n)$ der einfachen Strukturen von $X_n$ für $n\in\{1,2,3,5\}$, sodass das folgende Diagramm kommutiert und die Abbildungen $i_n$ Inklusionen von Unter-CW-Komplexen sind.
\[\xymatrix{
A_1 \arincl[rr]^{i_1} \arinclinv[d]^{i_2} \ar[rd]^{\alpha_1} && A_3 \arincl[rr]^{i_3} \ar[rd]^{\alpha_3} && A_5 \ar[rd]^{\alpha_5}\\
A_2 \ar[rd]^{\alpha_2} & X_1 \ar[rr]^{k_1} \ar[d]^{j_1} && X_3 \ar[rr]^{k_3} \ar[d]^{j_3} && X_5 \ar[d]^{j_5} \\
& X_2 \ar[rr]^{k_2}  && X_4 \ar[rr]^{k_4}  && X_6
}\]
Seien $A_4$ das Push-out von $(i_1,i_2)$ und $A_6$ das Push-out von $(i_3\circ i_1, i_2)$; bezeichne mit $\alpha_4:\,A_4\lto X_4$ und $\alpha_6:\,A_6\lto X_6$ die entsprechenden induzierten Abbildungen. Dann ist zu zeigen: $\tau(\alpha_6) \in N$.
Da das linke Quadrat einfach modulo $k_{4*}^{-1}(N)$ ist, gilt nach Definition $\tau(\alpha_4)\in k_{4*}^{-1}(N)$. Die Summenformel modulo $N$ im rechten Quadrat angewendet ergibt dann wie gewünscht $\tau(\alpha_6) = k_{4*}\tau(\alpha_4) \in N.$
\end{proof}

\begin{lem}["`Kommutativität"']\label{theta:kommutativitaet}
Seien zu den Abbildungen $f:\,A_0\lto X_0$ und $g:\,B_0\lto X_0$ die folgenden Push-outs von Räumen mit einfacher Struktur gegeben, in denen die vertikalen Inklusionen Kofaserungen sind:
\[\xymatrix{
A_0 \ar[rr]^f \arinclinv[d] && X_0 \arinclinv[d] &&B_0 \ar[rr]^g \arinclinv[d] && X_1 \arinclinv[d]^{j} \\
A_1 \ar[rr]^F && X_1 && B_1 \ar[rr]^G && X_2
}\]
Ist für eine Untergruppe $N\subset\Whp{X_2}$ das rechte Push-out einfach modulo $N$ und das linke Push-out einfach modulo $j_*^{-1}(N)$, so ist auch das folgende Push-out einfach modulo $N$:
\[\xymatrix{
A_0 \amalg B_0 \ar[rr]^{f\amalg g} \arinclinv[d] && X_0 \arinclinv[d] \\
A_1 \amalg B_1 \ar[rr]^{F\amalg G} && X_2
}\]
Es spielt also die Reihenfolge der Abbildungen $f$ und $g$ keine Rolle.
\end{lem}

\begin{proof}[Beweis]
Nach dem Transitivitätslemma \ref{theta:transitivitaet} genügt es zu zeigen, dass in der Folge von Push-out-Diagrammen
\[\xymatrix{
A_0 \amalg B_0 \ar[rr]^{f\amalg g} \arinclinv[d] && X_0 \arinclinv[d] \\
A_1 \amalg B_0 \ar[rr]^{F\amalg g} \arinclinv[d] && X_1 \arinclinv[d]^{j}\\
A_1 \amalg B_1 \ar[rr]^{F\amalg G} && X_2
}\]
das untere Quadrat einfach modulo $N$ und das obere Quadrat einfach modulo $j_*^{-1}(N)$ ist. Für das untere Quadrat folgt dies etwa erneut mit Lemma \ref{theta:transitivitaet}, angewendet auf die Folge von Push-out-Diagrammen:
\[\xymatrix{
B_0 \arincl[rr] \arinclinv[d] && A_1 \amalg B_0 \ar[rr]^{F\amalg g} \arinclinv[d] && X_1 \arinclinv[d]^{j}\\
B_1 \arincl[rr] && A_1 \amalg B_1 \ar[rr]^{F\amalg G} && X_2
}\]
Denn das linke Push-out ist einfach, und damit stimmt die vom rechten Push-out auf $X_2$ induzierte einfache Struktur modulo $N$ mit der vom äußeren Rechteck induzierten einfachen Struktur überein; diese einfache Struktur wiederum ist nach Voraussetzung modulo $N$ die gegebene. --- Der Beweis für das obere Quadrat verläuft entsprechend.
\end{proof}

\section{Einfache Faserungen}

Wir betrachten nun eine Faserung $p:\,E\lto B$ über einem wegzusammenhängenden Raum $B$, deren Faser $F$ den Homotopietyp eines endlichen CW-Komplexes hat. Sei $b\in B$ ein Punkt des Basisraums. Der Fasertransport entlang von Wegen im Basisraum definiert nach Kapitel \ref{fas} eine Abbildung 
\[\omega:\,\pi_1(B,b)\lto [F_b,F_b],\]
wobei wir wie üblich für zwei Räume $X,Y$ mit $[X,Y]$ die Menge der Homotopieklassen von Abbildungen $X\lto Y$ bezeichnen. Wir bezeichnen hier und im Folgenden mit $j_b:\,F_b\lto E$ die Inklusion der Faser und definieren die Abbildung
\[\theta^{b}(p):\,\pi_1(B,b)\lto\Whp{E},\quad x\mapsto j_{b*}\,\tau(\omega(x)).\]
Hierbei sei $\tau$ die Whitehead-Torsion bezüglich einer gewählten einfachen Struktur $\xi$ der Faser $F_b$. Die folgende Rechnung zeigt, dass die Wahl der einfachen Struktur dabei keine Rolle spielt. Seien dazu $\xi$ und $\eta$ zwei einfache Strukturen auf der Faser.
\begin{equation*}\begin{split}
& j_{b*}\, \tau\bigl[(F_b,\eta)\xrightarrow{\omega(x)} (F_b,\eta) \bigr] \\
  = \; & j_{b*}\, \tau\bigl[(F_b,\eta)\xrightarrow{\id} (F_b,\xi) \xrightarrow{\omega(x)} (F_b,\xi) \xrightarrow{\id} (F_b,\eta) \bigr]\\
 = \; & j_{b*}\, \tau\bigl[(F_b,\eta)\xrightarrow{\id} (F_b,\xi)\bigr] + j_{b*}\, \tau\bigl[(F_b,\xi) \xrightarrow{\omega(x)} (F_b,\xi)\bigr] + j_{b*}\,\omega(x)_*\, \tau\bigl[(F_b,\xi)\xrightarrow{\id} (F_b,\eta)\bigr]\\
 = \; & j_{b*}\, \tau\bigl[(F_b,\eta)\xrightarrow{\id} (F_b,\xi)\bigr] + j_{b*}\, \tau\bigl[(F_b,\xi) \xrightarrow{\omega(x)} (F_b,\xi)\bigr] - j_{b*}\, \tau\bigl[(F_b,\eta)\xrightarrow{\id} (F_b,\xi)\bigr]\\
 = \; & j_{b*}\, \tau\bigl[(F_b,\xi) \xrightarrow{\omega(x)} (F_b,\xi)\bigr].
\end{split}\end{equation*}
In die vorletzte Umformung ging dabei ein, dass gemäß Satz \ref{fas:fasertransport} in $[F_b,E]$ die Elemente $[j_b]$ und $[j_b]\circ\omega(x)$ gleich sind, also dieselbe Abbildung $\Whp{F_b}\lto \Whp{E}$ induzieren. 

\begin{lem}\label{theta:lemma5}
\begin{enumerate}\smallsep
\item Die Abbildung $\theta^b(p)$ ist ein Gruppenhomomorphismus.
\item Seien $b'\in B$ ein weiterer Punkt, $\gamma$ ein Weg in $B$ von $b'$ nach $b$ und $c_{\gamma}:\,\pi_1(B,b)\lto \pi_1(B,b')$ der durch Konjugation mit $\gamma$ gegebene Isomorphismus. Dann ist das folgende Diagramm kommutativ:
\[\xymatrix{
\pi_1(B,b) \ar[rr]^{c_{\gamma}} \ar[rd]_{\theta^b(p)} && \pi_1(B,b') \ar[ld]^{\theta^{b'}(p)}\\
& \Whp{E}
}\]
\end{enumerate}
\end{lem}

\begin{proof}[Beweis]
(i) Dies folgt wieder aus der Tatsache, dass $[j_b]=[j_b]\circ\omega(x)\in[F_b,E]$:
\[\theta^b(p)(y*x) = j_{b*}\,\tau(\omega(x)\circ\omega(y)) = j_{b*}\,\tau(\omega(x)) + j_{b*}\,\tau(\omega(y)) = \theta^b(p)(y) + \theta^b(p)(x).\]
(ii) Es sei $\omega_\gamma:\,F_{b'}\lto F_b$ ein Fasertransport entlang $\gamma$. Da die Wahl der einfachen Strukturen auf der Faser keine Rolle spielt, wählen wir einfache Strukturen auf $F_{b'}$ und $F_b$ so, dass $\omega_\gamma$ einfach ist. Dann gilt:
\[\tau((\omega_\gamma)^{-1}\circ\omega(x)\circ\omega_\gamma) = \omega_{\gamma*}^{-1}\, \tau(\omega(x))\]
und somit
\[\theta^{b'}(p)(c_\gamma(x))  = j_{b'*}\,\tau((\omega_\gamma)^{-1}\circ\omega(x)\circ\omega_\gamma) =   j_{b*}\,\tau(\omega(x))= \theta^b(p)(x).\qedhere\]
\end{proof}

Somit induziert für einen wegzusammenhängenden Basisraum $B$ die Abbildung $\theta^b(p)$ einen Gruppenhomomorphismus $\theta(p):\,H_1(B)\lto \Whp{E}$, der zudem unabhängig von der Wahl des Basispunkts ist. Wir werden diese Abbildung auch mit dem universellen Koeffiziententheorem als ein Element $\theta(p)\in H^1(B;\Whp{E})$ interpretieren.

\begin{defn}
Sei $p:\,E\lto B$ eine Faserung über einem zusammenhängenden Basisraum, deren Faser vom Homotopietyp eines endlichen CW-Komplexes ist. Dann heißt $p$ \emph{einfach}, falls $\theta(p):\,H_1(B)\lto \Whp{E}$ verschwindet. Für eine Untergruppe $N\subset\Whp{E}$ heißt die Faserung \emph{einfach modulo $N$,} falls $\theta(p)$ nur Werte in $N$ annimmt.
\end{defn}

Zum Beispiel ist jedes lokal-triviale Faserbündel $p:\,E\lto B$ über einem wegzusammenhängenden Basisraum mit einem endlichen CW-Komplex als Faser eine einfache Faserung. Denn nach Kapitel \ref{fas} kann jeder Fasertransport in einem solchen Bündel als Homöomorphismus gewählt werden; wegen der topologischen Invarianz der Whitehead-Torsion ist die induzierte Abbildung $\pi_1(B,b)\lto \Whp{F_b}$ die Nullabbildung. 

Tatsächlich wird unser Interesse im Folgenden hauptsächlich einfachen Faserungen gelten. Der Grund, auch einfache Faserungen modulo $N$ zu betrachten, liegt darin, dass Einschränkungen von einfachen Faserungen auf zusammenhängende Teilräume des Basisraums im Allgemeinen nicht mehr einfach, sondern nur noch einfach modulo einer entsprechenden Untergruppe sind (siehe Korollar \ref{theta:kor1}). Die für uns relevanten Eigenschaften von einfachen Faserungen übertragen sich jedoch problemlos auf den Fall einfacher Faserungen modulo $N$.

\begin{lem}\label{theta:lemma6}
Sei $p$ wie oben. Für eine Untergruppe $N\subset\Whp{E}$ sind die folgenden Bedingungen äquivalent:
\begin{enumerate}\smallsep
\item $p$ ist einfach modulo $N$.
\item Für je zwei Punkte $b,b'$ in $B$ und einfache Strukturen auf $F_b$ und $F_{b'}$ gilt: Alle Fasertransporte $F_{b'}\lto F_b$ haben die gleiche Torsion modulo $j_{b*}^{-1}(N).$
\end{enumerate}
\end{lem}

\begin{proof}[Beweis]
Für die Implikation von (ii) nach (i) wähle $b=b'$ mit den gleichen einfachen Strukturen auf $F_b=F_{b'}$. Dann hat modulo $j_{b*}^{-1}(N)$ jeder Fasertransport die Torsion der Identität, also gilt 
\[j_{b*}\, \tau(\omega(x))\in N \quad\forall x\in\pi_1(B,b).\]
Dies bedeutet nach Definition, dass $p$ einfach modulo $N$ ist. Sei umgekehrt vorausgesetzt, dass $p$ einfach modulo $N$ ist. Seien zwei Punkte $b,b'$ in $B$ und einfache Strukturen auf $F_b$ und $F_{b'}$ fest gewählt. Für zwei Wege $\gamma_1$ und $\gamma_2$ in $B$ von $b'$ nach $b$ und deren zugehörige Fasertransporte $\omega_{\gamma_1}$ und $\omega_{\gamma_2}$ gilt dann in der Quotientengrppe $\Whp{E}/N$:
\[ [j_{b*}\,\tau(\omega_{\gamma_1})] = [j_{b*}\,\tau(\omega_{\gamma_2}\circ\omega_{\gamma_2}^{-1}\circ \omega_{\gamma_1})]  = [j_{b*}\,\tau(\omega_{\gamma_2})] + [j_{b'*}\,\tau(\omega_{\gamma_1 * \gamma_2^-})] = [j_{b*}\,\tau(\omega_{\gamma_2})]. \]
In dieser Rechnung bezeichnet $\gamma_2^-$ den Weg $\gamma_2$ rückwärts durchlaufen. Bei der letzten Umformung wurde verwendet, dass die Faserung als einfach modulo $N$ vorausgesetzt war. Aus der obigen Rechnung folgt $[\tau(\omega_{\gamma_2})] = [\tau(\omega_{\gamma_1})]\in\Whp{F_b}/j_{b*}^{-1}(N)$.
\end{proof}

Abschließend soll noch die Verträglichkeit von $\theta(p)$ mit Pull-backs untersucht werden.

\begin{prop}\label{theta:prop1}
Sei $f:\,B'\lto B$ eine Abbildung zwischen zwei wegzusammen\-hängenden Räumen. Bezeichne mit $p':\,f^*E\lto B'$ die mittels $f$ zurückgezogene Faserung. Dann ist das Diagramm 
\[\xymatrix{
{\Whp{(f^*E)}} \ar[rr]^{\bar{f}_*} && \Whp{E}\\
H_1(B') \ar[rr]^{f_*} \ar[u]^{\theta(p')} && H_1(B) \ar[u]^{\theta(p)}
}\]
kommutativ, wenn $\bar{f}:\,f^*E\lto E$ die Strukturabbildung des Pull-backs ist.
\end{prop}

\begin{kor}\label{theta:kor1}
In der Situation von Proposition \ref{theta:prop1} gilt: Ist $p$ einfach modulo $N$, so ist $p'$ einfach modulo $\bar{f}_*^{-1}(N)$. Insbesondere: Ist $p$ einfach und $f$ eine Homotopie-Äquivalenz, so ist auch $p'$ einfach.
\end{kor}

Der Beweis der Proposition fußt auf der Verträglichkeit des Pull-backs mit dem Fasertransport:

\begin{proof}[Beweis von Proposition \ref{theta:prop1}]
Zu einem $b\in B'$ wähle eine einfache Struktur auf der Faser $F_{f(b)}$ von $p$ über $f(b)$. Gleichzeitig ist $F_{f(b)}$ Faser von $p'$ über $b$. Wir bezeichnen mit $j_{f(b)}:\,F_{f(b)}\lto E$ und $j'_b:\,F_{f(b)}\lto f^*E$ die jeweiligen Inklusionen. Sei ferner $\omega':\,\pi_1(B',b)\lto [F_{f(b)},F_{f(b)}]$ der Fasertransport von $p'$. Nach Lemma \ref{fas:lemma5} ist das folgende Diagramm kommutativ: 
\[\xymatrix{
& \left[F_{f(b)},F_{f(b)}\right] \\
\pi_1(B',b) \ar[rr]^{f_*} \ar[ru]^{\omega'} && \pi_1(B,f(b)) \ar[lu]_{\omega}
}\]
Folglich gilt für $x\in\pi_1(B,b)$:
\[\bar{f}_*\, \theta^b(p')(x) = \bar{f}_*\, j'_{b*}\, \tau(\omega'(x)) = j_{f(b)*}\, \tau[\omega(f_*(x))] = \theta^{f(b)}(p)[f_*(x)].\qedhere\]
\end{proof}

\section{Das Element $\theta(f)$}

Sei nun $f:\,X\lto B$ eine Abbildung in einen wegzusammenhängendem Bildraum $B$. Betrachte die assoziierte Faserung
\[\xymatrix{
X \ar[rr]^f \ar[rd]_\lambda^\simeq & & B\\
 & X^f \ar[ru]_p
}\]
Wir setzen im Folgenden voraus, dass die Faser $F$ der assoziierte Faserung den Homotopietyp eines endlichen CW-Komplexes hat. Denkt man an Frage, ob die Abbildung $f$ homotop zu einem differenzierbaren Faserbündel ist, ist dies eine sinnvolle Annahme. Denn nach Kapitel \ref{fas} gilt: Ändert man $f$ homotop ab, so sind die assoziierten Faserungen faserhomotopie-äquivalent; insbesondere haben deren Fasern denselben Homotopietyp. Ist $f$ bereits selbst eine Faserung, so stimmt die assoziierte Faserung $p$ bis auf Faserhomotopie-Äquivalenz mit $f$ überein. Die Faser eines differenzierbaren Faserbündels von geschlossenen differenzierbaren Mannigfaltigkeiten ist jedoch selbst eine geschlossene differenzierbare Mannigfaltigkeit und insbesondere vom Homotopietyp eines endlichen CW-Komplexes.

Zur Abbildung $f$ definieren wir nun
\[\theta(f) := \lambda_*^{-1}\,\theta(p) \in H^1(B;\Whp{X}). \]
Im Falle, dass $f=p$ bereits eine Faserung ist, ist die Verträglichkeit dieser Definition von $\theta(p)$ mit der im letzten Unterkapitel gemachten zu zeigen. Tatsächlich gilt folgende Invarianz des Elementes $\theta(p)$ unter Faserhomotopie-Äquivalenzen:

\begin{prop}\label{theta:prop2}
Seien $p:\,E\lto B$ und $p':\,E'\lto B$ Faserungen; sei $f:\,E\lto E'$ eine Faserhomotopie-Äquivalenz. Dann gilt: 
\[\theta(p')=f_*\,\theta(p)\in H^1(B;\Whp{E'}).\] 
Insbesondere gilt: $p'$ ist genau dann einfach modulo $N$, wenn $p$ einfach modulo $N$ ist.
\end{prop}

Proposition \ref{theta:prop2} werden wir am Ende dieses Unterkapitels beweisen. Sie zeigt, dass im Fall einer Faserung $f=p:\,X\lto B$ die Definition von $\theta(f)$ aus diesem Unterkapitel mit der Definition von $\theta(p)$ im letzten Unterkapitel übereinstimmt. Denn in diesem Fall ist $\lambda:X\lto X^f$ eine Faserhomotopie-Äquivalenz nach Proposition \ref{fas:prop2}. --- Wir haben bereits festgestellt, dass $\theta(f)$ verschwindet, falls $f$ ein lokal-triviales Faserbündel mit einem endlichen CW-Komplex als Faser ist. Dies trifft insbesondere auf den Fall eines differenzierbaren Faserbündels von geschlossenen Mannigfaltigkeiten zu und zeigt damit den ersten Teil von Proposition \ref{theta:prop_haupt}. 

\begin{prop}
Das Element $\theta(f)$ hängt nur von der Homotopieklasse von $f$ in $[X,B]$ ab.
\end{prop}

\begin{proof}[Beweis]
Eine Homotopie $H:\,X\times I\lto B$ zwischen $f,g:\,X\lto B$ induziert nach Proposition \ref{fas:prop3} ein bis auf Homotopie kommutatives Diagramm
\[\xymatrix{
& X \ar[ld]_{\lambda} \ar[rd]^{\lambda} \\
X^f \ar[rr]^{\omega_H} && X^g
}\]
mit einer Faserhomotopie-Äquivalenz $\omega_H:\,X^f\lto X^g$ der assoziierten Faserungen. Wieder liefert Proposition \ref{theta:prop2} die Behauptung.
\end{proof}

Damit ist Proposition \ref{theta:prop_haupt} bewiesen. Es bleibt der Beweis von Proposition \ref{theta:prop2} nachzutragen.

\begin{lem}\label{theta:lemma9}
Sei $p:\,E\lto B$ eine Faserung. Seien $f,g:\,X\lto E$ zwei über $B$ faserhomotope Abbildungen. Die Abbildungen $F$ bzw.~$G$ seien Lösungen der HHP
$$\xymatrix{
X\times\{0\} \ar[rr]^f \arinclinv[d] & & E \ar[d]^p & & X\times\{0\} \ar[rr]^g \arinclinv[d] & & E\ar[d]^p\\
X\times I \ar@{.>}[rru]^F \ar[rr]^h & & B & & X\times I \ar@{.>}[rru]^G \ar[rr]^h  && B
}$$
Dann sind $F$ und $G$ über $B$ faserhomotop. Insbesondere sind für alle $t\in I$ die induzierten Abbildungen
\[F_t,G_t: \,X\times\{t\} \lto (h_t)^*E\]
faserhomotop.
\end{lem}

\begin{proof}[Beweis]
Sei $H:\,X\times I\lto E$ eine Faserhomotopie zwischen $f$ und $g$. Setze die Abbildung $h:\,X\times I\lto B$ stationär zu einer Homotopie $h':\,X\times I\times I\lto B$ fort. Setze wie in Kapitel \ref{fas} $\dot{I}:= I\times\partial I\cup \{0\}\times I$, und definiere eine Abbildung $l:\,X\times\dot{I}\lto E$ durch 
\[l\vert_{X\times I\times \{0\}} = F,\quad l\vert_{X\times I\times \{1\}} = G,\quad l\vert_{X\times\{0\}\times I} = H.\]
Aus der HHE der Faserung $p$ erhält man eine Hochhebung $L:\,X\times I\times I\lto E$ von $h'$ nach $E$, die die Abbildung $l$ erweitert. Die Abbildung $L$ ist dann eine Faserhomotopie zwischen $F$ und $G$.

Die Einschränkung von $L$ auf $X\times \{t\}\times I$ definiert zusammen mit der Projektion $X\times \{t\}\times I\lto X$ auf die erste Komponente eine Faserhomotopie zwischen den Abbildungen $F_t$ und $G_t$.
\end{proof}

\begin{lem}\label{theta:lemma8}
Seien $p:\,E\lto B$ und $p':\,E'\lto B$ Faserungen; sei $\varphi:\,E\lto E'$ eine Faserhomotopie-Äquivalenz. Seien $H:\,X\times I\lto B$ eine Homotopie zwischen $f,g:\,X\lto B$ und $\omega_H:\,f^*E\lto g^*E$ bzw.~$\omega'_H:\,f^*E'\lto g^*E'$ Fasertransporte entlang $H$. Dann kommutiert das Diagramm
$$\xymatrix{
f^*E \ar[rr]^{f^*\varphi} \ar[d]_{\omega_H} & & f^*E' \ar[d]^{\omega'_H}\\
g^*E \ar[rr]^{g^*\varphi} & & g^*E'
}$$
bis auf Faserhomotopie.
\end{lem}

\begin{proof}[Beweis]
Sei $\psi:\,E'\lto E$ ein Faserhomotopie-Inverses von $\varphi$. Das mittlere Rechteck im folgenden Diagramm
$$\xymatrix{
f^*E' \times\{0\} \ar[rr]^{f^*\psi}_\simeq \arinclinv[d] & & f^*E\times\{0\} \ar[rr]^{\bar{f}} \arinclinv[d]& & E \ar[rr]^{\varphi}_\simeq \ar[d]^p & & E' \ar[d]^{p'}\\
f^*E'\times I \ar[rr]^{f^*\psi\times\id_I} & & f^*E\times I \ar[r] \ar@{.>}[rru]^{H'}& X\times I \ar[r]^H & B \ar@{=}[rr] & & B
}$$
stellt ein HHP dar; der diagonale Pfeil $H'$ sei eine Lösung desselben. Offenkundig stellt aber auch das ganze Diagramm ein HHP für die Faserung $p'$ dar, das von der Abbildung 
$$H'':=\varphi\circ H'\circ(f^*\psi\times\id_I):\,f^*E'\times I\lto E'$$ 
gelöst wird. Nun sind die Abbildungen $\bar{f}$ und $\varphi\circ\bar{f}\circ f^*\psi$ faserhomotop; somit ist die Abbildung $H''$ ist nach Lemma \ref{theta:lemma9} faserhomotop zur Lösung des folgenden HHP
$$\xymatrix{
f^*E'\times\{0\} \ar[rr]^{\bar{f}} \arincl[d] & & E'\ar[d]\\
f^*E'\times I \ar[r] \ar@{.>}[rru]^{H'''} & X\times I \ar[r]^H & B
}$$
Insbesondere definieren die Auswertungen von $H''$ und $H'''$ an $f^*E'\times\{1\}$ homotope Abbildungen $f^*E'\lto g^*E'$ Im Falle von $H''$ ist dies die Abbildung $g^*\varphi\circ\omega_H\circ f^*\psi$, im Falle von $H'''$ die Abbildung $\omega'_H$.
\end{proof}

\begin{kor}\label{theta:kor2}
Seien $F\lto E\lto B$ und $F'\lto  E'\lto B$ Faserungen; sei $f:\,E\lto E'$ eine Faserhomotopie-Äquivalenz. Es bezeichne $c_f:\,[F_{b_0},F_{b_0}]\lto [F'_{b_0},F'_{b_0}]$ die Konjugation mit $f$. Dann ist das Diagramm \[\xymatrix{
 & \pi_1(B,b_0) \ar[ld]_\omega \ar[rd]^{\omega'}\\
[F_{b_0},F_{b_0}]  \ar[rr]^{c_f} & & [F'_{b_0},F'_{b_0}]
}\]
wenn $\omega$ und $\omega'$ die jeweiligen Fasertransporte in $p$ und $p'$ bezeichnen.
\end{kor}

\begin{proof}[Beweis von Proposition \ref{theta:prop2}]
Wähle eine einfache Struktur $\xi$ auf $F_b$ und statte $F'_b$ mit der einfachen Struktur $f_*(\xi)$ aus. Dann gilt nach Korollar \ref{theta:kor2} bezüglich dieser einfachen Strukturen: $\tau(\omega'(x)) = f_*\,\tau(\omega(x))$ für alle $x\in\pi_1(B,b_0)$. Bezeichnen $j_b:\,F_b\lto E$ und $j'_b:\,F'_b\lto E'$ die Inklusionen der Fasern, so gilt also für alle $x\in H_1(B)$:
\[ \theta(p')(x) = j'_{b*}\,\tau(\omega'(x)) = j'_{b*}\,f_*\,\tau(\omega(x)) = f_*\,j_{b*}\,\tau(\omega(x)) = f_*\,\theta(p)(x),\]
wenn man wieder $\theta(p)$ und $\theta(p')$ als Abbildungen $H_1(B)\lto\Whp{E}$ bzw.~$H_1(B)\lto\Whp{E'}$ interpretiert.
\end{proof}


\chapter{Das Element $\tau(f)$}\label{tau}

Dieses Kapitel hat zum Ziel, einer Abbildung $f:\,M\lto B$ von geschlossenen und zusammenhängenden differenzierbaren Mannigfaltigkeiten mit $\theta(f)=0$ unter gewissen Zusatzvoraussetzungen ein Element $\tau(f)\in\Whp{M}$ zuzuordnen. Wieder spielt dabei die assoziierte Faserung $F\lto M^f\lto B$ der Abbildung $f$ eine wesentliche Rolle. Während das Element $\theta(f)$ die Torsion von Fasertransporten in $M^f$ misst, soll das Element $\tau(f)$ die "`Einfachheit"' von Trivialisierungen repräsentieren, wenn $M^f$ auf geeignete kontrahierbare Unterräume von $B$ eingeschränkt wird. Genauer kann man unter Verwendung einer CW-Struktur des Basisraums Trivialisierungen über den $n$-Zellen von $B$ auf geeignete Weise zusammenfügen, um ein endliches CW-Modell des Totalraums $M^f$ zu erhalten. Die dadurch definierte einfache Struktur $\xi$ auf $M^f$ ist in wichtigen Fällen wohldefiniert, und kann mit der durch $M$ selbst gegebenen einfachen Struktur verglichen werden. In diesem Sinne definieren wird das Element $\tau(f)$ im Wesentlichen als die Whitehead-Torsion der kanonischen Inklusion $\lambda:\, M\lto (M^f,\xi)$. Ganz in Analogie zum Element $\theta(f)$ gilt hier:

\begin{prop}\label{tau:prop_haupt}
Das Element $\tau(f)$ hängt nur von der Homotopieklasse von $f$ in $[M,B]$ ab. Ist $f$ ein differenzierbares Faserbündel, so ist $\tau(f)=0$, falls es definiert ist.
\end{prop}

In den ersten beiden Unterkapiteln geht es darum, die angesprochene einfache Struktur auf $M^f$ zu konstruieren. Diese Konstruktion werden wir allgemein durchführen für einfache Faserungen über einem endlichen zusammenhängenden CW-Komplex, deren Faser den Homotopietyp eines solchen hat. Anschließend beweisen wir einige Eigenschaften der einfachen Struktur auf dem Totalraum; sodann definieren wir das Element $\tau(f)$ und können Proposition \ref{tau:prop_haupt} beweisen.

\section{Eine einfache Struktur auf dem Totalraum}\label{tau:abschn1}

Es sei $p:\,E\lto B$ eine Faserung über einem endlichen zusammenhängenden CW-Komplex $B$, deren Faser den Homotopietyp eines endlichen CW-Komplexes hat. Es sei eine einfache Struktur $\zeta$ auf der Faser $F_b$ über dem Punkt $b\in B$ fest gewählt. Wir werden nun induktiv über das CW-Gerüst des Basisraums eine einfache Struktur auf $E$ konstruieren. Diese einfache Struktur hängt zunächst ab von im Verlauf der Konstruktion getroffenen Wahlen; tatsächlich spielen jedoch für einfache Faserungen, d.~h.~für Faserungen $p$ mit $\theta(p)=0$, diese Wahlen keine Rolle. Für eine ähnliche Konstruktion in einem etwas anderen Kontext vgl.~\cite{LNM1408}, Seite 122ff. Eine explizite induktive Konstruktion von endlichen CW-Modellen des Totalraums in unserem Kontext wird etwa in \cite{lueck:diss}, Kapitel 2, durchgeführt.

Wir stellen ein Lemma voran, das wir im Folgenden häufig verwenden werden:

\begin{lem}[\cite{LNM1408}, Seite 19]\label{tau:lemma1}
Sei $p:\,E\lto B$ eine Faserung; sei
$$\xymatrix{
B_0 \ar[rr]^{j_2} \arinclinv[d]^{j_1} \ar[rrd]^k && B_2 \arinclinv[d]^{k_2} \\
B_1 \ar[rr]^{k_1} && B
}$$
ein Push-out-Diagramm, in dem die vertikalen Inklusionen Kofaserungen sind. Dann ist auch das durch Pull-back entstehende kommutative Diagramm
$$\xymatrix{
k^*E \ar[rr]^{\bar{j}_2} \arinclinv[d]^{\bar{j}_1} && k_2^*E \arinclinv[d]^{\bar{k}_2} \\
k_1^*E \ar[rr]^{\bar{k}_1} && E
}$$
ein Push-out, und die vertikalen Inklusionen sind Kofaserungen.
\end{lem}

Für $n\in\NN$ notieren wir mit $B^{(n)}$ das $n$-Gerüst von $B$ und setzen $E^{(n)}:=p^{-1}(B^{(n)})$. Es bezeichne $j^{(n)}:\, E^{(n)}\lto E$ die Inklusion. Das 0-Gerüst von $B$ ist die disjunkte Vereinigung $B^{(0)} = \coprod_{i\in I_0} \{b_i\}$ von endlich vielen Punkten aus $B$. Wähle für $i\in I_0$ einen Fasertransport $\omega_i:\,F_b\lto F_{b_i}$ und statte $F_{b_i}$ mit der einfachen Struktur $\omega_{i*}(\zeta)$ aus. Diese wurde in Kapitel \ref{theta} definiert als die einfache Struktur auf $F_{b_i}$, bezüglich der $\omega_{\gamma_i}$ einfach ist, wenn $F_b$ die einfache Struktur $\zeta$ trägt. Wir erinnern ebenfalls daran, dass die disjunkte Vereinigung von Räumen mit einfacher Struktur auf kanonische Weise wieder eine einfache Struktur besitzt. Auf diese Weise erhält $E^{(0)}=\coprod_{i \in I_0} F_{b_i}$ eine einfache Struktur $\xi^{(0)}$. 

Sei nun als Induktionsvoraussetzung eine einfache Struktur $\xi^{(n)}$ auf $E^{(n)}$ bereits konstruiert. Der Raum $B^{(n+1)}$ geht aus $B^{(n)}$ durch das Ankleben von endlich vielen Zellen hervor; wähle ein zugehöriges Push-out-Diagramm 
\[\xymatrix{
{\coprod_{i\in I_n}} S^n \ar[rrr]^{\coprod_{i\in I_n} q_i^{(n)}} \arinclinv[d] &&& B^{(n)} \arinclinv[d] \\
{\coprod_{i\in I_n}} D^{n+1} \ar[rrr]^(0.55){\coprod_{i\in I_n} Q_i^{(n+1)}} &&& B^{(n+1)}
}\]
mit einer endlichen Indexmenge $I_n$. Gemäß Lemma \ref{tau:lemma1} ist dann auch das Diagramm
\begin{equation}
\begin{split}\label{tau:push-out}
\xymatrix{
{\coprod_{i\in I_n}}\, q_i^{(n)*} E^{(n)} \ar[rrr]^(0.55){\coprod_{i\in I_n} \overline{q_i^{(n)}}} \arinclinv[d] &&& E^{(n)} \arinclinv[d] \\
{\coprod_{i\in I_n}}\, Q_i^{(n+1)*} E^{(n+1)} \ar[rrr]^(0.6){\coprod_{i\in I_n} \overline{Q_i^{(n+1)}}} &&& E^{(n+1)}
}
\end{split}
\end{equation}
ein Push-out-Diagramm, wobei die vertikalen Inklusionen Kofaserungen sind. Für $i\in I_n$ ist $Q_i^{(n+1)*} E^{(n+1)}$ der Totalraum einer Faserung über dem zusammenziehbaren Raum $D^{n+1}$. Somit induziert nach Proposition \ref{fas:trivialisierung} eine Homotopie der charakteristischen Abbildung $Q_i:\,D^{n+1}\lto B$ zur Punktabbildung $D^{n+1}\lto \{b\} \lto B$ Trivialisierungen 
\begin{equation*}
\begin{split}
T_i:\,F_b\times D^{n+1} &\lto Q_i^{(n+1)*} E^{(n+1)}\\
t_i:\,F_b\times S^n &\lto q_i^{(n)*} E^{(n)}
\end{split}
\end{equation*}
von Faserungen über $D^{n+1}$ bzw.~$S^n$, wobei $t_i$ die Einschränkung von $T_i$ auf $F_b\times S^n$ ist. Insbesondere sind die beiden linken Räume im Diagramm \eqref{tau:push-out} vom Homotopietyp eines endlichen CW-Komplexes. Mit Hilfe der Trivialisierungen statten wir diese Räume mit den einfachen Strukturen $T_{i*}(\zeta\times[\id])$ bzw.~$t_{i*}(\zeta\times[\id])$ aus.

Legen wir also für Räume $Q_i^{(n+1)*} E^{(n+1)}$ und $q_i^{(n)*} E^{(n)}$ diese einfachen Strukturen, für den Raum $E^{(n)}$ die laut Induktionsvoraussetzung gegebene einfache Struktur $\xi^{(n)}$ zu Grunde, erhält auch der Raum $E^{(n+1)}$ durch das Push-out-Diagramm \eqref{tau:push-out} eine einfache Struktur $\xi^{(n+1)}$, und der Induktionsschritt ist vollzogen. Da die CW-Struktur auf $B$ laut Voraussetzung endlich ist, gilt $E^{(N)} = E^{(N+1)} = \ldots = E$ für ein hinreichend großes $N\in\NN$, und offenbar ist dann auch $\xi^{(N)} = \xi^{(N+1)} = \dots$ Somit setzen wir $\xi=\xi^{(N)}$ als einfache Struktur auf $E=E^{(N)}$. 

Diese einfache Struktur hängt zunächst ab von der einfachen Struktur $\zeta$, von der Wahl der charakteristischen Abbildungen $Q_i^{(n)}$ und der Wahl der Trivialisierungen bzw.~Fasertransporte. Man beachte, dass diese Wahlen auch für jeden Unterkomplex $B_1\subset B$ eine einfache Struktur $\xi_1$ auf $E_1:=p^{-1}(B_1)$ induzieren, indem die Indexmengen $I_n$ im Diagramm \eqref{tau:push-out} entsprechend verkleinert werden, die einfachen Strukturen auf den jeweiligen Summanden $Q_i^{(n+1)*}E^{(n+1)}$ und $q_i^{(n)*}E^{(n)}$ jedoch beibehalten werden. Wir untersuchen im Folgenden für den Fall einer einfachen Faserung $p:\,E\lto B$ die Abhängigkeit der einfachen Struktur $\xi$ auf $E$ bzw.~allgemeiner der einfachen Struktur $\xi_1$ auf $E_1$ von den getroffenen Wahlen. 

Seien dazu $b,b'\in B$ zusammen mit einfachen Strukturen $\zeta$ und $\zeta'$ auf den Fasern $F_b$ bzw.~$F_{b'}$ vorgegeben. Zu diesen Wahlen seien $\xi$ und $\eta$ einfache Strukturen auf dem Raum $E$, die induktiv nach der obigen Vorschrift zu den einfachen Strukturen $\zeta$ bzw.~$\zeta'$ konstruiert sind, und seien $\xi_1$ und $\eta_1$ die dazu gehörenden einfachen Strukturen auf $E_1$. 

\begin{lem}\label{tau:lemma2}
Ist für eine Untergruppe $N\subset\Whp{E}$ die Faserung $p:\,E\lto B$ einfach modulo $N$, so ist in $\Whp{E}$ die Identität
\begin{equation}\label{tau:eq2}
j_{1*}\, \tau\bigl[(E_1,\eta_1) \overset{\id}{\lto} (E_1, \xi_1) \bigr] 
\; \equiv\; \chi(B_1)\cdot j_{b*}\,\tau\bigl[(F_{b'},\zeta') \overset{\omega}{\lto} (F_b,\zeta) \bigr]\mod N
\end{equation}
erfüllt, wobei $j_1:\,E_1\lto E$ und $j_b:\,F_b\lto E$ die jeweiligen Inklusionen, $\omega: F_{b'}\lto F_b$ ein (beliebiger) Fasertransport und $\chi$ die Eulercharakteristik bezeichnen.
\end{lem}

Die Unabhängigkeit der rechten Seite in Gleichung \eqref{tau:eq2} von der Wahl des Fasertransports $\omega$ wurde bereits in Lemma \ref{theta:lemma6} gezeigt. Lemma \ref{tau:lemma2}, das im folgenden Unterkapitel bewiesen wird, ist ein zentrales Resultat dieses Teils, denn es zeigt den für uns im Folgenden wichtigen Satz:

\begin{satz}\label{tau:satz1}
Ist $p:\,E\lto B$ eine einfache Faserung modulo $N$, so hängt die einfache Struktur $\xi_1$ auf $E_1$ modulo $j_*^{-1}(N)$ nur von der Wahl von $b$ und der einfachen Struktur modulo $j_{b*}^{-1}(N)$ auf $F_b$ ab.
\end{satz}

Diese einfache Struktur $\xi_1$ modulo $j_*^{-1}(N)$ auf $E_1$ werden wir als \emph{kanonische einfache Struktur modulo $j_*^{-1}(N)$ zu $(b,\zeta)$} bezeichnen und schreiben $\xi_1=\xi_{B_1}(b,\zeta)$. Im Fall $B_1=B$ schreiben wir $\xi_B(b,\zeta)=\xi(b,\zeta)$. Ist $p$ einfach, so ist die einfache Struktur $\xi$ auf $E$ nur abhängig von $b$ und der einfachen Struktur $\zeta$ (modulo $\ker j_{b*}$) auf $F_b$.  

\begin{satz}\label{tau:satz2}
Ist eine der folgenden beiden Bedingungen erfüllt:
\begin{enumerate}\smallsep
\item $p$ ist einfach und $\chi(B)=0$,
\item die von der Inklusion induzierte Abbildung $j_*:\,\Whp{F}\lto\Whp{E}$ ist die Nullabbildung,
\end{enumerate}
so hat $E$ eine kanonische einfache Struktur $\xi=\xi(b,\zeta)$, die unabhängig von der Wahl von $b$ und $\zeta$ ist.
\end{satz}

\section{Wohldefiniertheit der einfachen Struktur}

Dieses Unterkapitel ist dem Beweis von Lemma \ref{tau:lemma2} gewidmet. Es seien daher im Folgenden die Voraussetzungen dieses Lemmas stets erfüllt. Wir werden die Gleichung \eqref{tau:eq2} induktiv über die Dimension von $B_1$ als CW-Komplex beweisen. 

Ein 0-dimensionaler endlicher CW-Komplex ist eine endliche disjunkte Vereinigung von Punkten. Gilt somit $\dim B_1=0$, so ist der Raum $E_1=\coprod_{i\in I_{-1}} F_{b_i}$ die disjunkte Vereinigung von Fasern, und die einfachen Strukturen $\xi_1$ bzw.~$\eta_1$ gehen aus einfachen Strukturen auf den jeweiligen Fasern hervor. Diese wiederum werden durch die Wahl von Fasertransporten $\omega_i$ bzw.~$\omega'_i$ aus den einfachen Strukturen $\zeta$ bzw.~$\zeta'$ erhalten. Damit gilt:
\begin{align*}
& j_{1*}\, \tau\Bigl[(\coprod_{i\in I_{-1}} F_{b_i},\eta_1) \overset{\id}{\lto} (\coprod_{i\in I_{-1}} F_{b_i}, \xi_1) \Bigr] \\
=\; &\sum_{i\in I_0} j_{b_{i}*}\,\tau\bigl[(F_{b_i},\omega'_{i*}(\zeta')) \overset{\id}{\lto} (F_{b_i},\omega_{i*}(\zeta)) \bigr]\\
=\; &\sum_{i\in I_0} j_{b*}\,\tau\bigl[(F_{b'},\zeta') \xrightarrow{\omega_i^{-1}\circ\omega'_i} (F_b,\zeta) \bigr].
\end{align*}
Für die Rechnung modulo $N$ spielt jedoch gemäß Lemma \ref{theta:lemma6} die Wahl der jeweiligen Fasertransporte keine Rolle. Daraus folgt Gleichung \eqref{tau:eq2} für $\dim B_1=0$, denn die Eulercharakteristik $\chi(B_1)$ ist in diesem Fall gleich der Mächtigkeit von $I_{-1}$.

Für den Induktionsschluss nehmen wir nun an, dass Gleichung \eqref{tau:eq2} für alle Unter-CW-Komplexe von $B$ der Dimension kleiner oder gleich $n$ bereits bewiesen ist. Sei $B_1$ nun ein $(n+1)$-dimensionaler Unter-CW-Komplex von $B$; wir wollen das $n$-Gerüst von $B_1$ mit $B_2$ bezeichnen. Die (induktive) Konstruktion von $\xi_1$ und $\eta_1$ erfolgte über die Konstruktion von einfachen Strukturen auf $E_2:=p^{-1}(B_2)$, die wir mit $\xi_2$ bzw.~$\eta_2$ bezeichnen. Für diese einfachen Strukturen ist nach Induktionsvoraussetzung Gleichung \eqref{tau:eq2} entsprechend gültig. Die jeweiligen einfachen Strukturen auf $E_1$ erhält man dann aus denen von $E_2$ durch Push-out-Diagramme vom Typ \eqref{tau:push-out}.

Ein zentraler Punkt des Induktionsschlusses wird sein, die Unabhängigkeit der einfachen Struktur von der Wahl der charakteristischen Abbildungen $Q_i$ zu zeigen. Seien dazu
\[\xymatrix{
{\copi} S^n \ar[rrr]^{\copi q_i} \arinclinv[d] &&& B_2 \arinclinv[d] && {\copii} S^n \ar[rrr]^{\copii q'_i} \arinclinv[d] &&& B_2 \arinclinv[d] \\
{\copi} D^{n+1} \ar[rrr]^{\copi Q_i} &&& B_1 && {\copii} D^{n+1} \ar[rrr]^{\copii Q'_i} &&& B_1
}\]
die in der Konstruktion von $\xi_1$ bzw.~$\eta_1$ verwendeten Darstellungen von $B_1$ als Push-outs. Zunächst erhalten wir Homöomorphismen
\[\xymatrix{
{\copi} \bigl(D^{n+1} - S^n\bigr) \ar[rr]^(0.6){\copi Q_i} && B_1 - B_2 && {\copii} \big(D^{n+1} - S^n\bigr) \ar[ll]_(0.6){\copii Q'_i},
}\]
die die Indexmengen $I_n$ und $I'n$ mit den Wegzusammenhangskomponenten von $B_1 - B_2$ identifizieren. Wir können also $I_n=I'_n$ mit $Q_i(D^{n+1})=Q'_i(D^{n+1})$ für alle $i\in I_n$ annehmen. Weiter können wir durch geeignete Homöomorphismen $(D^{n+1}, S^n)\lto (D^{n+1}, S^n)$ erreichen, dass für alle $i\in I_n$ gilt: $Q_i(0)=Q'_i(0)$.

Um die unterschiedlichen charakteristischen Abbildungen vergleichen zu können, werden wir im Folgenden das $n$-Gerüst $B_2$ durch die charakteristischen Abbildungen $Q_i$ in den Komplex $B_1$ "`hineinragen"' lassen und an dieses den Rest des $(n+1)$-Gerüstes mit Hilfe entsprechend modifizierter charakteristischer Abbildungen ankleben. Dies funktioniert im Einzelnen wie folgt:

Sei für $t\in [0,1]$ die Abbildung $\mu_t : D^{n+1}\lto D^{n+1}$ gegeben durch die Multiplikation mit $t$, und setze
\newlength{\mylengthone}
\newlength{\mylengthtwo}
\settowidth{\mylengthone}{$\mu_t(D^{n+1}),\quad S^n(t)$}
\settowidth{\mylengthtwo}{$Q_i\circ\mu_t,\quad q_i(t)$}
\begin{align*}
D^{n+1}(t) &= \mu_t(D^{n+1}),\hspace{10em}\hspace{-\mylengthone} S^n(t) = \mu_t(S^n),\\
Q_i(t) &= Q_i \circ \mu_t,\hspace{10em}\hspace{-\mylengthtwo}  q_i(t) = Q_i(t)\vert_{S^n},\\
B_2(t) &= B_2 \cup \bigcup_{i\in I_n} Q_i\big(\overline{D^{n+1} - D^{n+1}(t)}\big), \\
E_2(t) &= p^{-1}\big( B_2(t) \big).
\end{align*}
Mit diesen Bezeichnungen erhalten wir für $t>0$ eine neue Darstellung von $B_1$ als Push-out, dessen vertikale Inklusionen Kofaserungen sind, und das für $t=1$ ist das linke der oberen beiden Diagramme übergeht:
\[\xymatrix{
{\copi} S^n \arinclinv[d] \ar[rrr]^{\copi q_i(t)} &&& B_2(t) \arinclinv[d] \\
{\copi} D^{n+1} \ar[rrr]^{\copi Q_i(t)} &&& B_1
}\]
Der Vorteil dieses neuen Push-outs ist, dass die charakteristische Abbildungen $Q_i(t)$ für $t<1$ topologische Einbettungen sind. Mit Lemma \ref{tau:lemma1} erhalten wir ein entsprechendes Push-out-Diagramm auf Ebene der Totalräume, dessen vertikale Inklusionen Kofaserungen sind: 
\begin{equation} \label{tau:push-out2}
\begin{split}
\xymatrix{
{\copi} q_i(t)^* E_2(t) \arinclinv[d] \ar[rrr]^(0.6){\copi \overline{q_i(t)}} &&& E_2(t) \arinclinv[d] \\
{\copi} Q_i(t)^* E_1 \ar[rrr]^(0.6){\copi \overline{Q_i(t)}} &&& E_1
}
\end{split}
\end{equation}
Man beachte, dass dabei für $0<t<1$ die horizontalen Abbildungen wieder topologische Einbettungen sind. 

Wir versehen nun die Räume aus Diagramm \eqref{tau:push-out2} mit einer einfachen Struktur. Zunächst ist offenbar die Inklusion $B_2\lto B_2(t)$ und damit auch die Inklusion $E_2\lto E_2(t)$ eine Homotopie-Äquivalenz; wir statten über die Inklusion den Raum $E_2(t)$ mit der einfachen Struktur $\xi_2$ von $E_2$ aus. Weiter ist die Abbildung $\mu_t$ auf offensichtliche Weise zur Identität auf $D^{n+1}$ homotop, und diese Homotopie induziert für alle $i\in I$ eine Homotopie $H_i:\,D^{n+1}\times I\lto B_1$ zwischen den Abbildungen $Q_i(t)$ und $Q_i$. Diese Homotopie schränkt sich zu einer Homotopie $h_i:\,S^n\times I\lto B_2(t)$ auf $S^n$ ein. Betrachte die zugehörigen Faserhomotopie-Äquivalenzen
\begin{align*}
\omega_{H_i} :\,Q_i(t)^*E_1 &\lto Q_i^*E_1,\\
\omega_{h_i} :\,q_i(t)^*E_2(t) &\lto q_i^*E_2.
\end{align*}
Wir können gemäß Lemma \ref{fas:lemma4} annehmen, dass $\omega_{h_i}$ für alle $i\in I$ die Einschränkung von $\omega_{H_i}$ ist. In der Konstruktion der einfachen Struktur $\xi_1$ haben wir einfache Strukturen auf den Bildräumen von $\omega_{H_i}$ und $\omega_{h_i}$ konstruiert, die wir nun mit Hilfe der jeweiligen Abbildungen auf die Urbildräume transportieren.

\begin{lem}\label{tau:lemma3}
Mit diesen einfachen Strukturen ist Diagramm \eqref{tau:push-out2} ein Push-out von Räumen mit einfacher Struktur.
\end{lem}

\begin{proof}[Beweis]
Wir nennen $\xi'_1$ die einfache Struktur, die Diagramm \eqref{tau:push-out2} auf $E_1$ induziert. Es ist dann zu zeigen:
\[\tau\bigl[ (E_1,\xi'_1) \overset{\id}{\lto} (E_1,\xi_1) \bigr]=0.\]
Betrachte dazu die folgende natürliche Transformation von Push-out-Diagrammen:
\[\xymatrix{
{\copi} Q_i(t)^* E_1 \arinclinv[d] && {\copi} q_i(t)^* E_2(t) \arinclinv[ll] \arincl[rrr]^(0.6){\copi \overline{q_i(t)}} \arinclinv[d]^{j_0} &&& E_2(t)\ar@{=}[d]\\
{\copi} Q_i^* E_1 && {\copi} h_i^* E_2(t) \arinclinv[ll] \ar[rrr]^(0.6){\copi \bar{h}_i} &&& E_2(t) \\
{\copi} Q_i^* E_1 \ar@{=}[u] &&  {\copi} q_i^* E_2 \arinclinv[ll] \arincl[u]_{j_1} \ar[rrr]^(0.6){\copi \bar{q}_i} &&& E_2 \arincl[u]
}\]
Man beachte, dass man den Raum $Q_i(t)^* E_1= \mu_t^* Q_i^* E_1$ als Teilraum von $Q_i^* E_1$ auffassen kann. Ebenso ist auch $h_i^*E_2(t)$ ein Teilraum von $Q_i^*E_1$, denn $h_i$ ist eine Einschränkung von $Q_i$. Alle vertikalen Abbildungen des Diagramms sind Homotopie-Äquivalenzen, und die induzierten Abbildungen der jeweiligen Push-outs sind jeweils die Identität auf $E_1$.

Die Komposition $j_1^{-1}\circ j_0$ ist nach Kapitel \ref{fas} homotop zu $\omega_{h_i}$ und hat daher verschwindende Torsion. Die Inklusion in der rechten Spalte des oberen Diagramms ist nach Konstruktion einfach. Es ist also noch zu zeigen, dass auch die Inklusion in der linken Spalte einfach ist; aus der Summenformel folgt nämlich dann die Behauptung.

Tatsächlich ist die Inklusionsabbildung homotop zur Abbildung $\omega_{H_i}$ von oben. Dies folgt direkt aus Satz \ref{fas:fasertransport}, Teil (i), denn mit der Identifikation $Q_i(t)^*E_1 = \mu_t^* Q_i^* E_1$ ist die Abbildung $\bar{\mu}_t:\,Q_i(t)^*E_1\lto Q_i^*E_1$ die Inklusionsabbildung, und es gilt:
\[\bar{\mu}_t \simeq \bar{\mu}_1\circ\omega_{H_i} = \omega_{H_i}.\qedhere\]
\end{proof}

Dieselben Konstruktionen führen wir genauso für die gestrichten Abbildungen $Q'_i$ durch: Wir erhalten Räume $B_2(t)'$ und $E_2(t)'$ sowie Abbildungen $Q'_i(t)$ und $q'_i(t)$, und Lemma \ref{tau:lemma3} gilt entsprechend.

Die nächste Etappe besteht nun darin, die Torsion der Identität $(E_1,\eta_1)\lto (E_1,\xi_1)$ durch Anwendung der Summenformel in Diagrammen vom Typ \eqref{tau:push-out2} aus der Torsion von $\id:\,(E_2,\eta_2)\lto (E_2,\xi_2)$ zu berechnen, die laut Induktionsvoraussetzung bereits bekannt ist. Dieses Ziel wird mit Gleichung \eqref{tau:eq3} unten erreicht und wird den Induktionsschritt im Beweis von Gleichung \eqref{tau:eq2} erlauben.

\begin{lem}
Es gibt stetige, streng monoton wachsende Funktionen $\varepsilon,\delta:\,[0,1]\lto [0,1]$ mit
\[B_2(t)'\subset B_2(\varepsilon(t))\subset B_2(\delta(t))'\qquad (t\in[0,1]).\]
\end{lem}

\begin{figure*}[ht]
  \centering
  \includegraphics[width=\textwidth]{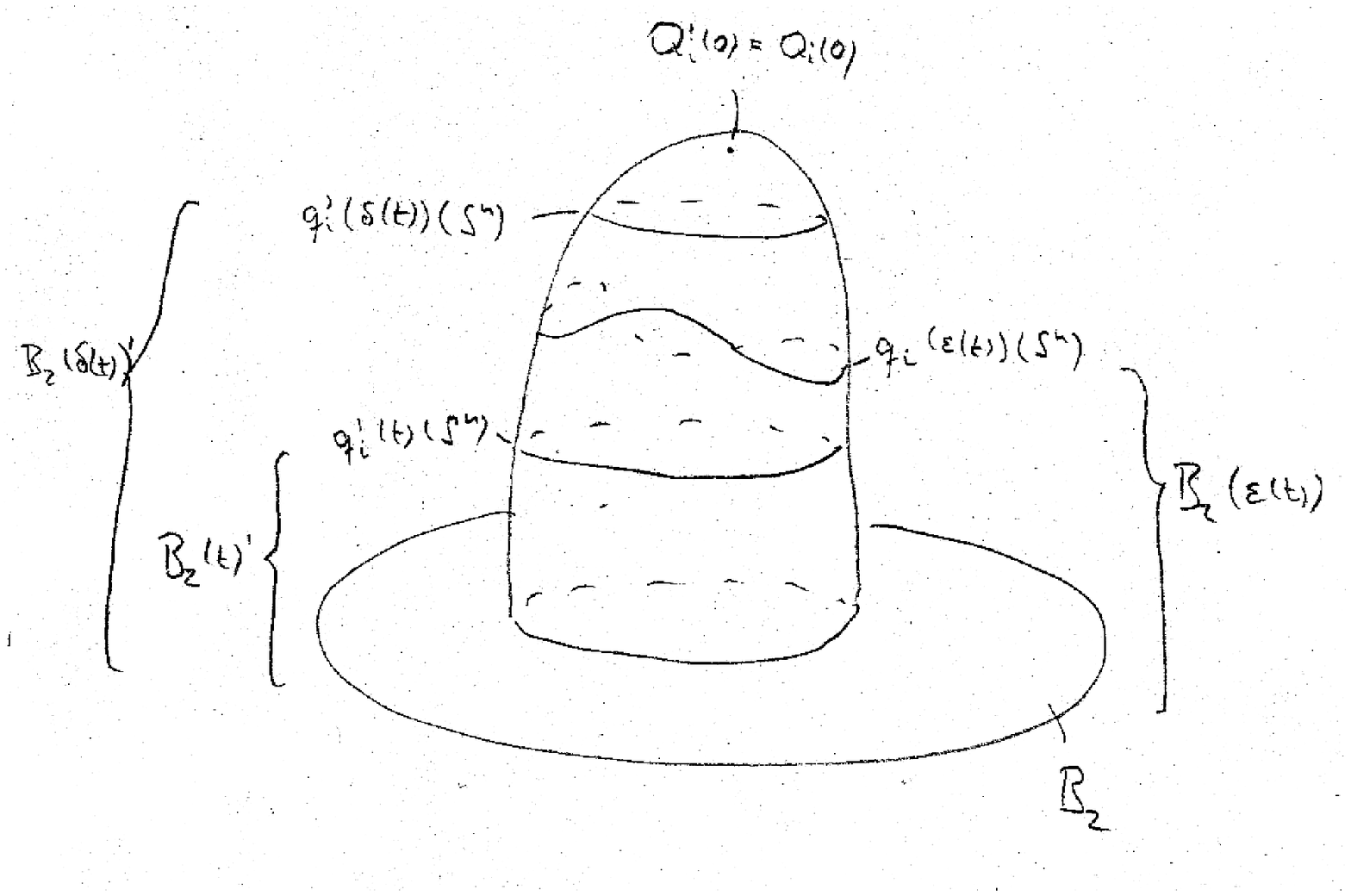}
  \caption{}
  \label{tau:fig1}
\end{figure*}

\begin{proof}[Beweis]
Die Abbildung 
\[\bar{\varepsilon}:\,[0,1)\lto [0,1),\quad t\mapsto \min \bigl\{\norm{x};\; x\in\bigcup_{i\in I_n} Q_i^{-1} (Q_i'(\overline{D^{n+1}-D^{n+1}(t)})) \bigr\}\]
ist offenbar stetig und monoton wachsend. Tatsächlich ist sie sogar streng monoton wachsend: Für $t=0$ ist nämlich $\bar{\varepsilon}(t)=0$; andererseits gilt für $t>0$ auch $\bar{\varepsilon}(t)>0$, und das Minimum von $\norm{x}$ wird stets auf dem Rand der betrachteten Menge angenommen. Somit erweitert sich die Abbildung $\bar{\varepsilon}$ zu einer stetigen und streng monoton wachsenden Funktion $\varepsilon:\,[0,1]\lto[0,1]$. Wir zeigen, dass diese Funktion $\varepsilon$ die geforderte Bedingung $B_2(t)'\subset B_2(\varepsilon(t))$ erfüllt. Nach Definition ist dazu die Inklusion
\[B_2\cup \bigcup_{i\in I_n}Q_i'(\overline{D^{n+1}-D^{n+1}(t)}) \subset B_2\cup \bigcup_{i\in I_n} Q_i(\overline{D^{n+1}-D^{n+1}(\varepsilon(t))})\]
zu zeigen. Sei dazu $p$ ein Element der linken Seite. Für $p\in B_2$ ist nichts zu zeigen. Sei also $p\notin B_2$, d.~h.~$p=Q_i'(x')$ für ein $i\in I_n$ mit einem $x'\in\overline{D^{n+1}-D^{n+1}(t)}$, $\norm{x'}<1$. Dann ist zugleich $p=Q_i(x)$ mit einem $x\in D^{n+1}$, $\norm{x}<1$, und es gilt
\[x\in Q_i^{-1} (Q_i'(\overline{D^{n+1}-D^{n+1}(t)})).\]
Nach Definition von $\varepsilon$ ist also $\varepsilon(t)\leq\norm{x}<1$ und $p\in Q_i(\overline{D^{n+1}-D^{n+1}(\varepsilon(t))})$. Dies war zu zeigen. --- Die Funktion $\delta$ definiert man nun analog.
\end{proof}

Fig.~\ref{tau:fig1} veranschaulicht die Situation für ein festes $t\in (0,1)$ und das Ankleben einer $(n+1)$-Zelle an $B_2$. --- Mit derartigen Funktionen $\varepsilon$ und $\delta$ definieren wir eine Homotopie $k:\,B_1\times I\lto B_1$ durch:
\begin{enumerate}\setlength{\itemsep}{0pt}
\item $k_t\vert_{B_2(t)'} = \id_{B_2(t)'}$ für alle $t\in [0,1]$, insbesondere $k_0 = \id_{B_1}$,
\item $k_t\circ Q'_i(s) = Q'_i(\frac{ts}{\norm{s}})$ für alle $t\in [0,1]$, $i\in I_n$ und $\delta(t)\leq\norm{s}<t$,
\item $k_t\circ Q'_i(s) = Q'_i(\frac{ts}{\delta(t)})$ für alle $t\in [0,1]$, $i\in I_n$ und $0<\norm{s}\leq\delta(t)$,
\item $k_t(Q'_i(0)) = Q_i'(0)$ für alle $i\in I_n$.
\end{enumerate}

\begin{lem}
Die Funktion $k:\, B_1\times I\lto B_1$ ist wohldefiniert und stetig.
\end{lem}

\begin{proof}[Beweis]
Wir betrachten zunächst den Bereich $B_2\times I$. Dort ist $k$ wohldefiniert, und zur Stetigkeit genügt es wegen Bedingung (i) zu zeigen, dass für $(\norm{s},t)\to (1,1)$ auch $k_t(Q'_i(s))$ gegen $Q'_i(s)$ konvergiert. Dies ist jedoch im Gebiet von (ii) klar, und das Gebiet von (iii) ist nur relevant, falls auch $\delta(t)\to 1$. Demnach ist $k$ auf $B_2\times I$ stetig. 

Für jede Zusammenhangskomponente von $(B_1-B_2)\times I$ ist nach Anwendung des Homöo\-morphismus $Q'_i:\,D^{n+1}-S^n\lto B_1-B_2$ an die Wohldefiniertheit und Stetigkeit der Funktion
\[(D^{n+1}-S^n)\times I\lto D^{n+1}-S^n,\quad (s,t)\mapsto
\begin{cases}
s &, \quad\norm{s}\geq t,\\
\frac{ts}{\norm{s}} &, \quad\delta(t)\leq\norm{s}\leq t,\\
\frac{ts}{\delta(t)}& , \quad 0<\norm{s}\leq\delta(t),\\
0& , \quad s=0
\end{cases}\]
zu zeigen. Dies ist unmittelbar nachzuprüfen.
\end{proof}

Wähle nun ein $t\in (0,1)$ und setze $f:=k_t:\,B_1\lto B_1$. Die Abbildung $f$ ist vermöge der Homotopie
\[k':\,B_1\times I\lto B_1,\quad (x,t')\mapsto k(x,t\cdot t')\]
homotop zur Identität. Man beachte außerdem, dass die Abbildung $f$ den Teilraum $B_2(\varepsilon(t))$ auf $B_2(t)'$, die Teilräume $q_i(\varepsilon(t))(S^n)$ für alle $i\in I_n$ auf $q'_i(t)(S^n)$ und die Teilräume $Q_i(\varepsilon(t))(D^{n+1})$ auf $Q'_i(t)(D^{n+1})$ abbildet. $f$ faktorisiert dann als die folgende natürliche Transformation von Push-out-Diagrammen:
\begin{equation}\label{tau:diag1}
\begin{split}
\xymatrix{
{\copi} D^{n+1} \ar[d]^{\copi V_i} &&{\copi} S^n \arinclinv[ll] \ar[d]^{\copi v_i} \arincl[rrr]^{\copi q_i(\varepsilon(t))} &&& B_2(\varepsilon(t))\ar[d]^r \\
{\copi} D^{n+1}  &&{\copi} S^n \arinclinv[ll]\arincl[rrr]^{\copi q'_i(t)} &&& B_2(t)'
}
\end{split}
\end{equation}

Wir bezeichnen mit $(\sigma_{D^{n+1}},\sigma_{S^n}):\, (D^{n+1},S^n)\lto (D^{n+1},S^n)$ die Spiegelung an der Äquatorebene $x_{n+1}=0$. Es gilt:

\begin{lem}\label{tau:lemma4}
Für alle $i\in I_n$ existiert eine Homotopie $(L,l):\,(D^{n+1},S^n)\times I\lto (D^{n+1},S^n)$ mit Start $(L_0,l_0)=(V_i,v_i)$ und mit Ziel entweder $(L_1,l_1)=(\id_{D^{n+1}},\id_{S^n})$ oder $(L_1,l_1)=(\sigma_{D^{n+1}},\sigma_{S^n}).$
\end{lem}

\begin{proof}[Beweis]
Der Beweis fußt auf der Feststellung, dass $v_i$ eine Homotopie-Äquivalenz ist. Im kommutativen Diagramm
\[\xymatrix{ S^n \arincl[r]^(0.3){\mu_{\varepsilon(t)}} \ar[d]^{v_i} & D^{n+1} - (S^n\cup \{0\}) \ar[r]^(0.4){Q_i} & Q_i(D^{n+1}) - \big(Q_i(S^n) \cup \{Q_i(0)\} \big) \ar[d]^f \\
S^n \arincl[r]^(0.3){\mu_t}  & D^{n+1} - (S^n\cup \{0\}) \ar[r]^(0.4){Q'_i} & Q_i(D^{n+1}) - \big(Q_i(S^n) \cup \{Q_i(0)\} \big)
}\]
sind nämlich die Abbildungen $Q_i$ und $Q'_i$ Homöomorphismen, die Abbildungen $\mu_{\varepsilon(t)}$ und $\mu_t$ Homotopie-Äquivalenzen, und die Abbildung $f$ (gemeint ist die entsprechende Einschränkung der oberen Abbildung $f$) vermöge der entsprechenden Einschränkung von $k$ homotop zur Identität, also ebenfalls eine Homotopie-Äquivalenz. Es folgt somit, dass auch $v_i$ eine Homotopie-Äquivalenz ist und somit als Abbildungsgrad den Wert 1 oder $-1$ hat. 

Wir setzen $\alpha=\id_{D^{n+1}}$, falls der Abbildungsgrad von $v_i$ gleich 1 ist, ansonsten $\alpha=\sigma_{D^{n+1}}$. Dann ist in beiden Fällen $v_i$ homotop zu $\alpha\vert_{S^n}$. Mit Hilfe der Kofaserungseigenschaft des Paares $(D^{n+1}, S^n)$ weiten wir diese Homotopie zu einer Homotopie $L':\, D^{n+1}\times I\lto D^{n+1}$ zwischen $V_i$ und einer Abbildung $L'_1:\, D^{n+1}\lto D^{n+1}$ aus, deren Einschränkung auf $S^n$ gleich $\alpha\vert_{S^n}$ ist. Dann ist
$$L'':\, D^{n+1}\times I\lto D^{n+1},\quad (x,t)\mapsto \big(t\cdot \alpha(x) + (1-t)\cdot L'_1(x)\big).$$
eine Homotopie relativ $S^n$ zwischen $L'_1$ und $\alpha$, und die Zusammensetzung von $L'$ und $L''$ ist die gesuchte Homotopie $L$. Aus der Konstruktion folgt, dass die Einschränkung von $L$ eine Homotopie $l:\, S^n\times I \lto S^n$ zwischen $v_i$ und $\alpha\vert_{S^n}$ ist.
\end{proof}

Wir werden im Folgenden ohne Beschränkung der Allgemeinheit davon ausgehen, dass $(V_i,v_i)$ homotop zur Identität ist; andernfalls wende den Homöomorphismus $(\sigma_{D^{n+1}},\sigma_{S^n})$ auf das Urbild von $(Q_i,q_i)$ an. --- Die Homotopie $k'$ induziert einen Fasertransport $\omega_{k'}:\,E_1\lto f^*E_1$ und eine Abbildung
$$g:=\bar{f}\circ\omega_{k'}:\,E_1\lto E_1$$
über $f$, die homotop zur Identität auf $E_1$ ist. Zudem ist $k'$ eine Homotopie relativ $B_2$, sodass die Einschränkung von $g$ zu einer Abbildung $E_2\lto E_2$ ebenfalls homotop zur Identität ist. Auch die Abbildung $g$ faktorisiert nun in der folgenden natürlichen Transformation von Push-out-Diagrammen, die über dem Diagramm \eqref{tau:diag1} liegt, d.~h.~die entsprechenden Abbildungen der Basisräume überdeckt.
\begin{equation}\label{tau:diag2}
\begin{split}
\xymatrix{
{\copi} Q_i(\varepsilon(t))^*E_1\ar[d]^{\copi \Phi_i}  &{\copi} q_i(\varepsilon(t))^*E_2(\varepsilon(t)) \ar[d]^{\copi \phi_i} \arinclinv[l]  \arincl[rr]^(0.65){\copi \overline{q_i(\varepsilon(t))}} && E_2(\varepsilon(t)) \ar[d]^{\psi}\\
{\copi} Q'_i(t)^*E'_1  &{\copi} q'_i(t)^*E_2(t)'  \arinclinv[l]\arincl[rr]^(0.65){\copi \overline{q_i(t)}} && E_2(t)'
}
\end{split}
\end{equation}

In diesem Diagramm sind alles vertikalen Abbildungen Homotopie-Äquivalenzen, da sie von Homotopie-Äquivalenzen der jeweiligen Basisräume induziert werden. Es folgt aus der Homotopieinvarianz und der Summenformel, wenn $j:\,E_2(t)\lto E_1$ die Inklusion bezeichnet:
\begin{equation}\begin{split}\label{tau:eq3}
&j_{1*}\,\tau\bigl[(E_1,\eta_1)\overset{\id}{\lto} (E_1,\xi_1)\bigr]\\
= \; &j_{1*}\,\tau\bigl[(E_1,\eta_1)\overset{g}{\lto} (E_1,\xi_1)\bigr]\\
= \; & \sum_{i\in I_n} j_{1*}\, \overline{Q'_i(t)}_*\, \tau(\Phi_i) - \sum_{i\in I_n} j_{1*}\, j_*\, \overline{q'_i(t)}_*\, \tau(\phi_i)  + j_{1*}\, j_*\, \tau(\psi)\\
=\; &  \sum_{i\in I_n} j_{1*}\, \overline{Q'_i(t)}_*\, \tau(\Phi_i)\, - \sum_{i\in I_n} j_{1*}\, j_*\, \overline{q'_i(t)}_*\, \tau(\phi_i) + j_{2*}\, \tau\bigl[(E_2,\eta_2)\overset{\id}{\lto} (E_2,\xi_2)\bigr].
\end{split}\end{equation}

Hierbei bezeichnet $j_2:\,E_2\lto E$ die Inklusion. In die letzte Umformung ging ein, dass die Einschränkung von $\psi$ auf $E_2$ homotop zur Identität ist, und die beiden Räume in der rechten Spalte von Diagramm \eqref{tau:diag2} ihre einfache Struktur aus der Inklusion von $E_2$ erhalten haben. 

Nun haben wir unser obiges Ziel erreicht, den Wechsel der einfachen Struktur auf $E_1$ mit Hilfe des Wechsels der einfachen Struktur auf $E_2$ zu berechnen. Um den Induktionsschritt durchzuführen, müssen wir allerdings noch die Torsionen der Abbildungen $\Phi_i$ und $\phi_i$ bestimmen. Dazu ändern wir zunächst diese Abbildungen entlang der Homotopien $L_i$ des Basisraums homotop in Faserhomotopie-Äquivalenzen $\bar{\Phi}_i$ und $\bar{\phi}_i$ ab, deren Torsionen wir anschließend leicht berechnen können. Dies geschieht auf die folgende Weise:

Die Konstruktion von $f$ hat bereits eine Homotopie $k'$ zwischen $f$ und $\id_{B^{n+1}}$ geliefert. Andererseits zeigt Lemma \ref{tau:lemma4} nach unseren gemachten Annahmen die Existenz von Homotopien $(L_i,l_i)$ für $i\in I_n$ zwischen $(V_i,v_i)$ und $(\id_{D^{n+1}},\id_{S^n})$. Wir bilden nun die zu diesen Homotopien gehörenden Faserhomotopie-Äquivalenzen und erhalten für $i\in I_n$ die folgende Abbildung:
\[\bar{\Phi}_i:\, Q_i(\varepsilon(t))^*E_1 \xrightarrow{Q_i(\varepsilon(t))^*\omega_{k'}}  Q_i(\varepsilon(t))^* f^* E_1= V_i^* Q'_i(t)^* E_1 \xrightarrow{\omega_{L_i}}  Q'_i(t)^* E_1.\]
Wir nennen $\bar{\phi}_i:\, q_i(\varepsilon(t))^*E_2(\varepsilon(t)) \lto q'_i(t)^* E_2(t)'$ die Einschränkung von $\bar{\Phi}_i$. 

\begin{lem}\label{tau:lemma5}
\begin{enumerate}\setlength{\itemsep}{0pt}
\item Es gibt eine Homotopie zwischen $(\Phi_i,\phi_i)$ und $(\bar{\Phi}_i,\bar{\phi}_i)$.
\item Mit den Bezeichnungen aus \eqref{tau:eq3} gilt:
\begin{align}\label{tau:eq4}
\begin{split}
& j_{1*} \,\overline{Q'_i(t)}_*\, \tau(\Phi_i) - j_{1*}\, j_*\, \overline{q'_i(t)}_*\, \tau(\phi_i) \\
= \; &(-1)^{n+1} \cdot j_{b*}\,\tau\bigl[(F_{b'},\zeta') \overset{\omega_i}{\lto} (F_b,\zeta) \bigr]
\end{split}
\end{align}
für einen Fasertransport $\omega_i$ von $b'$ nach $b$.
\end{enumerate}
\end{lem}

Teil (ii) dieses Lemmas ist der letzte nötige Schritt, um den Induktionsschluss durchzuführen. Denn sei $\omega:\,F_{b'}\lto F_b$ ein beliebiger Fasertransport. Dann gilt mit den Gleichungen \eqref{tau:eq3} und \eqref{tau:eq4} unter Verwendung von Lemma \ref{theta:lemma6}:
\begin{equation*}\begin{split}
&j_{1*}\,\tau\bigl[(E_1,\eta_1)\overset{\id}{\lto} (E_1,\xi_1)\bigr]\\
\equiv\; & (-1)^{n+1}\cdot \vert I_n\vert \cdot j_{b*}\,\tau\bigl[(F_{b'},\zeta') \overset{\omega}{\lto} (F_b,\zeta) \bigr] + \chi(B_2)\cdot j_{b*}\,\tau\bigl[(F_{b'},\zeta') \overset{\omega}{\lto} (F_b,\zeta) \bigr]\\
\equiv\; & \chi(B_1)\cdot j_{b*}\,\tau\bigl[(F_{b'},\zeta') \overset{\omega}{\lto} (F_b,\zeta) \bigr]\mod N,
\end{split}\end{equation*}
denn es ist $\chi(B_1)-\chi(B_2)=(-1)^{n+1}\cdot \vert I_n \vert$. Damit ist der Beweis von Lemma \ref{tau:lemma2} abgeschlossen.

\begin{proof}[Beweis von Lemma \ref{tau:lemma5}]
(i) Zunächst stellt man fest, dass die Komposition 
$$Q_i(\varepsilon(t))^*E_1 \xrightarrow{Q_i(\varepsilon(t))^*\omega_{k'}}  Q_i(\varepsilon(t))^* f^* E_1 \xrightarrow{\bar{f}\circ\overline{Q_i(\varepsilon(t))}} E_1$$
gleich der Abbildung $\overline{Q'_i(t)}\circ \Phi_i$ ist. Dies folgt aus der Definition von $\Phi_i$ bzw.~$\psi$. Folglich ist die Komposition
$$Q_i(\varepsilon(t))^*E_1 \xrightarrow{Q_i(\varepsilon(t))^*\omega_{k'}}  Q_i(\varepsilon(t))^* f^* E_1\\ 
  = V_i^* Q'_i(t)^* E_1 \xrightarrow{\bar{V}_i}  Q'_i(t)^* E_1$$
gleich der Abbildung $\Phi_i$, denn deren jeweilige Kompositionen mit $\overline{Q'_i(t)}$ stimmen überein, und $\overline{Q'_i(t)}$ ist injektiv. Demnach erhält man die Abbildung $\bar{\Phi}_i$ aus $\Phi_i$, indem man lediglich die letzte Abbildung $\bar{V}_i$ durch $\omega_{L_i}$ ersetzt.

Mit diesen Observationen folgt leicht die Aussage (i), denn die Abbildung $\bar{V}_i$ ist homotop zu $\omega_{L_i}$ nach Satz \ref{fas:fasertransport}, was eine Homotopie zwischen $\bar{\Phi}_i$ und $\Phi_i$ induziert. Diese Konstruktion ist verträglich mit den jeweiligen Einschränkungen auf $q_i(\varepsilon(t))^*E_2(\varepsilon(t))$.

(ii) Wegen (i) kann man in Gleichung \eqref{tau:eq4} die Abbildungen $\Phi_i$ und $\phi_i$ durch die Varianten mit Querstrich ersetzen. Vollzieht man die Definition der einfachen Strukturen auf $Q_i(\varepsilon(t))^*E_1$ und $Q'_i(t)^* E_1$ nach, erkennt man, dass sie über gewisse Trivialisierungen $T_i:\,F_b\times D^{n+1}\lto Q_i(\varepsilon(t))^*E_1$ und $T_i':\,F_{b'}\times D^{n+1}\lto Q'_i(t)^*E_1$ definiert sind. Nach Proposition \ref{fas:trivialisierung} existiert nun ein bis auf Faserhomotopie kommutatives Diagramm
$$\xymatrix{
F_b\times D^{n+1} \ar[rr]^{T_i} \ar[d]^{\omega_i\times \id_{D^{n+1}}} && Q_i(\varepsilon(t))^*E_1 \ar[d]^{\bar{\Phi}_i}\\
F_{b'}\times D^{n+1} \ar[rr]^{T_i'} && Q'_i(t)^*E_1
}$$
mit einem Fasertransport $\omega_i:\,F_b\lto F_{b'}$, denn auch die Komposition $\bar{\Phi}_i\circ T_i$ ist eine Trivialisierung. Es folgt:
\begin{align*}
& j_{1*}\, \overline{Q'_i(t)}_*\, \tau(\Phi_i) - j_{1*}\, j_*\, \overline{q'_i(t)}_*\, \tau(\phi_i) \\
=\; & j_{1*}\, \overline{Q'_i(t)}_*\, (T_i')_*\, \tau(\omega_i\times\id_{D^{n+1}}) - j_{1*}\, j_*\, \overline{q'_i(t)}_*\, (t_i')_*\, \tau(\omega_i\times\id_{S^n}) \\
=\; & \chi(D^{n+1}) \cdot j_{b*}\, \tau(\omega_i) - \chi(S^{n})\, \cdot j_{b*}\, \tau(\omega_i) \\
=\; &(-1)^{n+1} \cdot j_{b*}\, \tau(\omega_i).
\end{align*}
Damit ist Gleichung \eqref{tau:eq4} bewiesen.
\end{proof}

\section{Einige Eigenschaften dieser Konstruktion}

Nun werden wir einige Resultate beweisen, die im Zusammenhang mit der Konstruktion von einfachen Strukturen aus Unterkapitel \ref{tau:abschn1} stehen. Tatsächlich verhalten sich die oben konstruierten einfachen Strukturen "`gutartig"' in Bezug auf Push-outs, Faserhomotopie-Äquivalenzen, Trivialisierungen und Pull-backs. Diese Verträglichkeitsresultate sind die Grundlage für entsprechende Eigenschaften des noch zu definierenden Elements $\tau(f)$. Die teilweise etwas technische Formulierung sollte nicht darüber hinwegtäuschen, dass es sich stets um vollkommen kanonische Konstruktionen handelt.

Sei für den Moment noch $p:\,E\lto B$ eine nicht notwendig einfache Faserung über einem endlichen zusammenhängenden CW-Komplex; sei die Faser $F$ vom Homotopietyp eines endlichen CW-Komplexes und eine einfache Struktur $\zeta$ auf der Faser $F_b$ über einem Punkt $b\in B$ fixiert. Unsere Konstruktion stattet dann auch den Raum $E$ mit einer einfachen Struktur $\xi$ aus, für die nun allerdings nicht mehr das Eindeutigkeitsresultat von Lemma \ref{tau:lemma2} gültig ist. Ist weiter $B_1\subset B_2\subset B$ eine Kette von Unter-CW-Komplexen, so erhalten durch die im Verlauf der Konstruktion von $\xi$ getroffenen Wahlen auch die Räume $E_i:=p^{-1}(B_i)$ für $i=1,2$ wie in Unterkapitel \ref{tau:abschn1} beschrieben einfache Strukturen $\xi_i$.

Nun gilt, dass man sich in dieser Situation die einfache Struktur $\xi_2$ auf $E_2$ auch aus der einfachen Struktur $\xi_1$ auf $E_1$ konstruiert denken kann. Tatsächlich ist $(B_2,B_1)$ auch ein relativer CW-Komplex, wobei die Anklebeabbildungen für die relativen Zellen genau die anklebenden Abbildungen für die entsprechenden Zellen von $B$ sind. Entsprechend kann man unter Anwendung von Lemma \ref{tau:lemma1} auch $E_2$ aus $E_1$ durch endlich viele Push-outs von Typ des Diagramms \eqref{tau:push-out} erhalten. Somit erhält der Raum $E_2$ aus dem Raum $E_1$ und der einfachen Struktur $\xi_1$ induktiv über die Push-out-Diagramme vom Typ \eqref{tau:push-out} eine einfache Struktur $\xi'_2$, wenn man die Räume in der linken Spalte von \eqref{tau:push-out} mit denselben einfachen Strukturen versieht wie bereits bei der Konstruktion von $\xi$.

Da sich die jeweils auftretenden Push-out-Diagramme der beiden Konstruktionen nur in ihrer Reihenfolge, nicht jedoch in den einfachen Strukturen der anzuklebenden Räume unterscheiden, notieren wir als Konsequenz von Lemma \ref{theta:kommutativitaet}:

\begin{bem}\label{tau:lemma6}
Diese beiden Strukturen $\xi_2$ und $\xi'_2$ auf $E_2$ stimmen überein.
\end{bem}

Wir nehmen nun für den Rest dieses Unterkapitels an, dass die Faserung $p:\,E\lto B$ einfach modulo $N$ für eine Untergruppe $N\subset\Whp{E}$ ist. Wir untersuchen nun das Verhalten der kanonischen einfachen Strukturen unter Push-outs. Sei dazu $f:\,A\lto B$ eine Abbildung von einem zusammenhängenden CW-Komplex $A$. Seien $A_2$ und $B_2$ Unter-CW-Komplexe von $A$ bzw.~$B$, und das Diagramm
\begin{equation}\begin{split}\label{tau:push-out-von-cw-komplexen}
\xymatrix{
A_1 \ar[rr]^{f_1} \arinclinv[d] && B_1 \arinclinv[d]\\
A_2 \ar[rr]^{f_2} && B_2
}
\end{split}\end{equation}
sei ein Push-out mit den folgenden Eigenschaften:
\begin{itemize}
\item $(A_2,A_1)$ und $(B_2,B_1)$ sind endliche CW-Paare und
\item Ist $Q:\,D^n\lto A_2$ eine charakteristische Abbildung einer $n$-Zelle des relativen Komplexes $(A_2,A_1)$, so ist $f_2\circ Q$ eine charakteristische Abbildung einer $n$-Zelle des relativen Komplexes $(B_2,B_1)$. Dies bedeutet insbesondere, dass die offenen relativen $n$-Zellen von $(A_2,A_1)$ unter $f_2$ in Bijektion zu den offenen relativen $n$-Zellen von $(B_2,B_1)$ stehen.
\end{itemize}
Diese beiden Bedingungen garantieren die Verträglichkeit der CW-Strukturen auf den einzelnen Räumen unter der Push-out-Struktur. Zum Beispiel ist jedes zelluläre Push-out ein Push-out von CW-Komplexen. Ein Push-out, das diese Bedingungen erfüllt, wollen wir kurz ein \emph{Push-out von CW-Komplexen} nennen.

Setze $E_i:=p^{-1}(B_i)$ für $i=1,2$ und seien $j_i:\,E_i\lto E$ die jeweiligen Inklusionen. Nach Wahl von $a\in A_1$ sowie einer einfachen Struktur $\zeta$ auf $F_{f(a)}$ stattet unsere obige Konstruktion die Räume $E_i$ mit einfachen Strukturen $\xi_i=\xi_{B_i}(f(a),\zeta)$ aus, die eindeutig modulo $j_{i*}^{-1}(N)$ sind.

Man beachte, dass die Faser $F_{f(a)}$ zugleich die Faser der zurückgezogenen Faserung $f^*E$ über dem Punkt $a$ ist; seien $\eta$ eine dazu konstruierte einfache Struktur auf $f^*E$ und $\eta_i$ für $i=1,2$ die dazu gehörenden einfachen Strukturen auf $f_i^*E_i = (f^*E)\vert_{A_i}$. (Diese einfachen Strukturen sind wohldefiniert modulo geeigeneter Untergruppen, denn die Faserung $f^*E$ ist einfach modulo $\bar{f}_*^{-1}(N)$; jedoch werden wir diese Tatsache hier nicht benötigen.) 

\begin{lem}[Verträglichkeit mit Push-outs]\label{tau:lemma10}
Ist Diagramm \eqref{tau:push-out-von-cw-komplexen} ein Push-out von CW-Komplexen, so ist bezüglich der einfachen Strukturen $\eta_i$ und $\xi_i$ ist das Push-out
\[\xymatrix{
f_1^*E_1 \ar[rr]^{\bar{f}_1} \arinclinv[d] && E_1 \arinclinv[d]\\
f_2^*E_2 \ar[rr]^{\bar{f}_2} && E_2
}\]
einfach modulo $j_{2*}^{-1}(N)$.
\end{lem}

\begin{proof}[Beweis]
Nach Bemerkung \ref{tau:lemma6} kann man sich auch $\eta_2$ vorstellen als eine zum relativen CW-Gerüst von $(A_2,A_1)$ konstruierte, zu $(b,\zeta)$ und $\eta_1$ gehörende einfache Struktur auf $f_2^*E_2$. Gemäß dem Transitivitätslemma \ref{theta:transitivitaet} genügt es nun, die Behauptung für den Fall zu zeigen, in dem $A_2$ aus $A_1$ nur durch das Ankleben von $n$-Zellen hervorgeht. Seien $(Q_i,q_i):\,(D^n,S^{n-1})\lto (A_2,A_1)$ für $i$ in einer endlichen Indexmenge $I_n$ die entsprechenden charakteristischen Abbidungen, die bei der Konstruktion von $\eta_2$ aus $\eta_1$ verwendet wurden. Betrachte das Diagramm
\begin{equation}\begin{split}\label{tau:diag3}\xymatrix{
{\copi} q_i^*f_1^*E_1 \ar[rr]^(0.55){\copi \bar{q}_i} \arinclinv[d] && f_1^*E_1 \ar[rr]^{\bar{f}_1} \arinclinv[d] && E_1 \arinclinv[d]\\
{\copi} Q_i^*f_2^*E_2 \ar[rr]^(0.55){\copi \bar{Q}_i} && f_2^*E_2 \ar[rr]^{\bar{f}_2} && E_2
}\end{split}\end{equation}
In der Konstruktion von $\eta$ wurden die Räume in der linken Spalte mit Hilfe von Trivialisierungen mit einfachen Strukturen ausgestattet. Mit diesen einfachen Strukturen vollzieht das linke Quadrat dann die Konstruktion von $\eta_2$ nach; demnach ist das linke Push-out einfach. 

Sei $\xi'_2$ die einfache Struktur, mit der das rechte Push-out den Raum $E_2$ ausstattet. Wir müssen zeigen, dass $\xi_2$ und $\xi'_2$ modulo $j_{2*}^{-1}(N)$  übereinstimmen. Dazu genügt es nach Bemerkung \ref{tau:lemma6} zu zeigen, dass $\xi'_2$ eine einfache Struktur zu $(f(a),\zeta)$ ist, die der relative CW-Komplex $(B_2,B_1)$ zusammen mit $\xi_1$ auf $E_2$ induziert. 

Um dies einzusehen, bemerken wir durch eine erneute Anwendung von Lemma \ref{theta:transitivitaet}, dass $\xi'_2$ die einfache Struktur ist, die das gesamte Push-out \eqref{tau:diag3} auf $E_2$ induziert. Nach Voraussetzung sind jedoch die Abbildungen $(f_2\circ Q_i,f_1\circ q_i)$ für $i\in I_n$ zugleich charakteristische Abbildungen für $(B_2,B_1)$ als relativer CW-Komplex. Das Diagramm \eqref{tau:diag3} vollzieht somit tatsächlich die Konstruktion einer einfachen Struktur zu $(f(a),\zeta)$ auf $E_2$ aus der einfachen Struktur $\xi_1$ auf $E_1$ nach. Dies war zu zeigen. 
\end{proof}

Als nächstes untersuchen wir die Verträglichkeit der kanonischen einfachen Struktur mit Faserhomotopie-Äquivalenzen. Sei dazu nun $p:\,E\lto B$ eine einfache Faserung modulo $N$. Wir nehmen die Existenz einer Faserhomotopie-Äquivalenz $f:\,E'\lto E$ mit einer Faserung $p':\,E'\lto B$ an. Ist $B_1\subset B$ ein Unterkomplex, so erhält der Raum $E\vert_{B_1}$ zu $b\in B_1$ und einer einfachen Struktur $\zeta$ auf der Faser $F_b$ von $p$ über $b$ eine kanonische einfache Struktur $\xi_{B_1}(b,\zeta)$ modulo $j_{1*}^{-1}(N)$, wenn $j_1:\,E\vert_{B_1}\lto E$ wieder die Inklusion bezeichnet.

Gemäß Proposition \ref{theta:prop2} ist nun aber auch die Faserung $p'$ einfach modulo $f_*^{-1}(N)$, und zu $(b,f_*^{-1}(\zeta))$ erhält somit auch der Raum $E'\vert_{B_1}$ eine kanonische einfache Struktur $\xi_{B_1}(b,f_*^{-1}(\zeta))$ modulo $j_{1*}^{-1}(N)$.

\begin{lem}[Verträglichkeit mit Faserhomotopie-Äquivalenzen]\label{tau:lemma8}
Es gilt
\[\xi_{B_1}(b,f_*^{-1}(\zeta)) = f_*^{-1}\,\xi_{B_1}(b,\zeta)\]
als Identität von einfachen Strukturen modulo $j_{1*}^{-1}(N)$.
\end{lem}

\begin{proof}[Beweis]
Wir notieren $F_b$ bzw.~$F'_b$ für $b\in B$ stets die jeweiligen Fasern von $p$ bzw.~$p'$, sowie $E_1:=E\vert_{B_1},\xi_1=\xi_{B_1}(b,\zeta)$ und analog für die gestrichten Varianten $E'_1$ und $\xi'_1$. Wir zeigen nun induktiv über die Dimension $n$ von $B_1$, dass die Abbildung $f:\,(E'_1,\xi'_1)\lto (E_1,\xi_1)$ einfach modulo $j_{1*}^{-1}(N)$ ist. 

Für $n=0$ ist
\[f:\,E'_1=\coprod_{i\in I_{-1}} F'_{b_i} \lto \coprod_{i\in I_{-1}} F_{b_i} = E_1.\]
Die Räume $F_{b_i}$ und $F'_{b_i}$ erhalten ihre jeweiligen einfachen Strukturen durch die Wahl von Fasertransporten $F_{b_i}\lto F_b$ bzw.~$F'_{b_i}\lto F'_b$. Auf Grund der Eindeutigkeit der einfachen Struktur $\xi_1$ modulo $j_{1*}^{-1}(N)$ können wir davon ausgehen, dass hierbei jeweils Fasertransporte entlang derselben Wege $\gamma_i$ von $b_i$ nach $b$ in $B$ verwendet wurden. Wir bezeichnen mit $\omega_{\gamma_i}$ und $\omega'_{\gamma_i}$ Fasertransporte entlang $\gamma_i$ in $E$ bzw.~$E'$. Nun ist das Diagramm 
\[\xymatrix{
F'_b \ar[rr]^f \ar[d]^{\omega'_{\gamma_i}} && F_b \ar[d]^{\omega_{\gamma_i}}\\
F'_{b_i} \ar[rr]^f && F_{b_i}
}\]
nach Lemma \ref{theta:lemma8} bis auf Homotopie kommutativ. Nach der Kettenregel ist also die Behauptung für $n=0$ richtig.

Sei nun vorausgesetzt, dass die Behauptung für alle Unterkomplexe von $B$ der Dimension $n$ erfüllt ist, und sei $B_1$ ein Unterkomplex von $B$ der Dimension $n+1$. Sei $B_2$ das $n$-Gerüst von $B_1$ und setze $E_2:=E\vert_{B_2}$ sowie $E'_2:=E'\vert_{B_2}$. Nach Wahl von charakteristischen Abbildungen $Q_i$ erhalten $E_1$ und $E'_1$ ihre einfachen Strukturen aus denen von $E_2$ bzw.~$E'_2$ durch ein Push-out wie in Diagramm \eqref{tau:push-out}. Die Abbildung $f:\, E'_1\lto E_1$ ist zugleich die von der universellen Eigenschaft des Push-outs induzierte Abbildung in folgendem kommutativen Diagramm:
\begin{equation}\begin{split}\label{tau:diag4}
\xymatrix{
{\copi}Q_i^*E'_1 \ar[d]^{\copi Q_i^*f} && {\copi}q_i^*E'_2 \arinclinv[ll] \ar[rr]^(0.6){\copi \bar{q}_i} \ar[d]^{\copi q_i^*f} && E'_2 \ar[d]^f\\
{\copi}Q_i^*E_1 && {\copi}q_i^*E_2 \arinclinv[ll] \ar[rr]^(0.6){\copi \bar{q}_i} && E_2
}\end{split}\end{equation}

Nach der Induktionsvoraussetzung und der Summenformel für die Whitehead-Torsion genügt es also zu zeigen, dass die linke und die mittlere vertikale Abbildung in Diagramm \eqref{tau:diag4} einfache Homotopie-Äquivalenzen sind. Wir zeigen diese Aussage für Abbildungen $Q_i^*f$; den Beweis für die Abbildungen $q_i^*f$ erhält man auf dieselbe Weise, indem man alle auftretenden Homotopie-Äquivalenzen entsprechend einschränkt. 

Die Räume in der linken Spalte von \eqref{tau:diag4} erhalten ihre einfachen Strukturen über Trivialisierungen. Wieder können wir davon ausgehen, dass diese Trivialisierungen entlang derselben Nullhomotopie der Basisabbildungen $Q_i$ erhalten wurde. Gemäß Lemma \ref{theta:lemma8} kommutiert dann für alle $i\in I_n$ das folgende Diagramm bis auf Homotopie.
$$\xymatrix{
F'_b\times D^{n+1} \ar[rrr]^{f\times\id_{D^{n+1}}} \ar[d]^{T'_i} &&& F_b\times D^{n+1} \ar[d]^{T_i} \\
Q_i^*E'_1 \ar[rrr]^{Q_i^*f} &&& Q_i^*E_1
}$$
Die obere horizontale Abbildung ist nach der Produktregel einfach. Wiederum liefert nun die Kettenregel die Behauptung. 
\end{proof}

Die triviale Faserung $p:\,B\times F\lto B$, ist als Faserbündel einfach. Damit induziert jede Wahl von $b\in B$ und einer einfachen Struktur $\zeta$ auf $F$ eine einfache Struktur $\xi(b,\zeta)$ auf $B\times F$. Diese ist als Folge von Lemma \ref{tau:lemma2} offenbar unabhängig von der Wahl des Basispunktes. Hier soll nun diese Struktur $\xi(b,\zeta)$ mit der einfache Struktur $[\id]\times\zeta$ auf $B\times F$ verglichen werden. Tatsächlich zeigen wir wenigstens unter gewissen Voraussetzungen an den Basisraum (die in den für uns relevanten Anwendungen erfüllt sind) dass die beiden einfachen Strukturen dieselben sind. Dazu bezeichnen wir einen CW-Komplex $X$ als \emph{gut,} wenn man alle Anklebeabbildungen $q_i^{(n)}:\,S^n\lto X^{(n)}$ zellulär wählen kann (d.~h.~so wählen kann, dass sie den Basispunkt der $S^n$ ins 0-Gerüst $X^{(0)}$ von $X$ abbilden). Der technische Vorteil von guten CW-Komplexen besteht darin, dass die definierenden Push-outs
\[\xymatrix{
{\copi} S^n \ar[rrr]^{\copi q_i^{(n)}} \arinclinv[d] &&& B^{(n)} \arinclinv[d] \\
{\copi} D^{n+1} \ar[rrr]^{\copi Q_i^{(n+1)}} &&& B^{(n+1)}
}\]
zelluläre Push-outs sind, sodass die Summenformel für die Whitehead-Torsion anwendbar ist. 

\begin{lem}\label{tau:lemma9}
Für die triviale Faserung $p:\,B\times F\lto B$ über einem guten endlichen zusammenhängenden CW-Komplex und jeden Unterkomplex $B_1\subset B$ gilt
\[\xi_{B_1}(b,\zeta) = [\id]\times\zeta\]
als Gleichheit von einfachen Strukturen modulo $\ker j_{1*}$ auf $B_1\times F$, wenn $j_1:\,B_1\times F\lto B\times F$ die Inklusion bezeichnet.
\end{lem}

\begin{kor}[Verträglichkeit mit Trivialisierungen]\label{tau:kor1}
Seien $B$ ein kontrahierbarer endlicher guter CW-Komplex und $p:\,E\lto B$ eine Faserung. Seien $B_1\subset B$ ein Unterkomplex, $b\in B_1$ und $\zeta$ eine einfache Struktur auf der Faser $F_b$. Dann ist jede Trivialisierung 
\[T:\,(B_1,[\id])\times (F_b,\zeta) \lto (E\vert_{B_1},\xi_{B_1}(b,\zeta)) \]
einfach modulo $\ker j_{1*}$ mit der Inklusion $j_1:\,E\vert_{B_1}\lto E$.
\end{kor}

Das Korollar folgt mit der Verträglichkeit mit Faserhomotopie-Äquivalenzen, denn $T$ ist insbesondere eine Faserhomotopie-Äquivalenz, und es gilt $\zeta=T_*(\zeta)$, da jeder Fasertransport $F_b\lto F_b$ homotop zur Identität und damit einfach ist.

\begin{proof}[Beweis von Lemma \ref{tau:lemma9}]
Wieder geht der Beweis induktiv über die Dimension von $B_1$. Wir können davon ausgehen, dass $F$ ein endlicher CW-Komplex mit seiner natürlichen einfachen Struktur $[\id]$ ist, sonst wähle einen Repräsentanten $(A,f)$ von $\zeta$ und wende Lemma \ref{tau:lemma8} auf die Faserhomotopie-Äquivalenz $\id_B\times f:\,B\times A\lto B\times F$ an. 

Die Behauptung für $\dim B_1=0$ ist dann trivial. Nehmen wir als Induktionsvoraussetzung an, dass die Behauptung für alle Unter-CW-Komplexe der Dimension $n$ bewiesen ist. Ist $B_1\subset B$ dann ein Unterkomplex der Dimension $n+1$, bezeichne mit $B_2$ das $n$-Gerüst von $B_1$, mit $j_2:\,B_2\times F\lto B\times F$ die Inklusion und wähle zelluläre charakteristische Abbildungen $(Q_i,q_i):\,(D^{n+1},S^n)\lto (B_1,B_2)$. Nach Induktionsvoraussetzung stimmt die einfache Struktur $\xi_{B_1}(b,\zeta)$ wenigstens modulo $\ker j_{2*}$ mit der natürlichen Stuktur auf $B_2\times F$ überein. Andererseits erhält $B_1\times F$ seine einfache Struktur aus der einfachen Struktur $\xi_{B_2}(b,\zeta)$ von $B_2\times F$ durch das Push-out
\[\xymatrix{
{\copi} S^n\times F \ar[rrr]^(0.55){\copi q_i\times\id_F} \arinclinv[d] &&& B_2\times F \arinclinv[d]\\
{\copi} D^{n+1}\times F \ar[rrr]^(0.58){\copi Q_i\times\id_F} &&& B_1\times F
}\]
Dies ist nun ein zelluläres Push-out endlicher CW-Komplexe, wobei die Räume in der linken Spalte ihre natürliche einfache Struktur tragen. Damit induziert das Push-out auch die natürliche einfache Struktur modulo $\ker j_{1*}$ auf $B_1\times F$.
\end{proof}

Abschließend wollen wir uns mit der Abhängigkeit der einfachen Struktur $\xi(b,\zeta)$ auf $E$ von der Wahl der CW-Struktur auf $B$ beschäftigen. Sei hierzu wieder $p$ eine einfache Faserung modulo $N$, und betrachte das Pull-back von $p$ mit einer einfache Homotopie-Äquivalenz $f:\,B'\lto B$ zu einer Faserung $p':\,f^*E\lto B'$ sowie die von $f$ induzierte Homotopie-Äquivalenz $\bar{f}:\,f^*E\lto E$. Sei $\zeta$ eine einfache Struktur von $p^{-1}(f(b'))=(p')^{-1}(b')$ Wir untersuchen nun den Zusammenhang zwischen den kanonischen einfachen Strukturen $\xi(f(b'),\zeta)$ auf $E$ und $\xi'(b,\zeta)$ auf $f^*E$. Dazu beachte man, dass das Pull-back $p':\,f^*E\lto B'$ von $p$ gemäß Korollar \ref{theta:kor1} ebenfalls einfach modulo $N$ ist.

\begin{lem}[Verträglichkeit mit Pull-backs]\label{tau:lemma7}
Ist $f$ einfach, so gilt 
\[\bar{f}_*\,\xi'(b,\zeta) = \xi(f(b),\zeta)\]
als Gleichheit einfacher Strukturen modulo $N$ auf $E$.
\end{lem}

\begin{proof}[Beweis]
Es ist zu zeigen, dass die Abbildung $\bar{f}:\,f^*E\lto E$ bezüglich der jeweiligen einfachen Strukturen einfach modulo $N$ ist. Dies vollzieht sich in drei Schritten. Um die Diskussion bezüglich der Basispunkte zu vereinfachen, bemerken wir, dass die Wahl einer einfachen Struktur $\zeta_b$ auf $F_b$ vermöge der Wahl von Fasertransporten für alle $b'\in B$ eine (im Allgemeinen von der Wahl der Wege abhängige) einfache Struktur $\zeta_{b'}$ auf $F_{b'}$ induziert. Lemma \ref{tau:lemma2} zeigt, dass die einfachen Strukturen $\xi(b,\zeta_b)$ und $\xi(b',\zeta_{b'})$ modulo $N$ übereinstimmen. In diesem Sinne werden wir im Folgenden die Basispunkte nicht mehr aufführen, sondern stets implizit annehmen, dass die jeweils auftretenden Basispunkte $b'$ die durch Fasertransport erhaltene einfache Struktur $\zeta_{b'}$ tragen.

(i) Wir zeigen: \textit{Lemma \ref{tau:lemma7} ist richtig, falls $f:\,B'\lto B$ eine elementare Expansion ist. }Sei
\begin{equation}\begin{split}\label{tau:diag5}\xymatrix{
S^n_+ \ar[rr]^q \arinclinv[d] && B' \arinclinv[d]^{f}\\
D^{n+1} \ar[rr]^Q && B
}\end{split}\end{equation}
eine elementare Expansion (vgl.~Unterkapitel \ref{theta:abschn1}). Wir erinnern daran, dass hier $D^{n+1}$ die CW-Struktur trägt, die aus einer 0-Zelle, einer $(n-1)$-Zelle (dem Äquator von $S^n$), zwei $n$-Zellen (nämlich $S^n_\pm$) und eine $(n+1)$-Zelle besteht. Bezüglich dieser CW-Struktur ist $D^{n+1}$ ein guter CW-Komplex im Sinne von Lemma \ref{tau:lemma9}, und Diagramm \eqref{tau:diag5} ist ein Push-out von CW-Komplexen. Aus der Verträglichkeit mit Push-outs geht hervor, dass das Diagramm
\[\xymatrix{
q^*(E\vert_{B'}) \ar[rr]^{\bar{q}} \arinclinv[d]^{\tilde{f}} && E\vert_{B'} \arinclinv[d]^{\bar{f}}\\
Q^*E \ar[rr]^{\bar{Q}} && E
}\]
ein einfaches Push-out modulo $N$ ist. Daher genügt es, die Einfachheit der Abbildung $\tilde{f}$ zu verifizieren, wie eine Anwendung der Summenformel in der folgenden natürlichen Transformation von Push-out-Diagrammen zeigt.
\[\xymatrix{
q^*(E\vert_{B'}) \arinclinv[d]^{\tilde{f}} \ar@{=}[rr] && q^*(E\vert_{B'}) \ar@{=}[d] \ar[rr]^{\gamma} && E\vert_{B'} \ar@{=}[d]\\
Q^*E && q^*(E\vert_{B'}) \arinclinv[ll] \ar[rr]^{\gamma} && E\vert_{B'}
}\]
Um die Einfachheit von $\tilde{f}$ einzusehen, wähle Trivialisierungen $(T,t)$ von $(Q^*E,q^*(E\vert_{B'}))$ und betrachte das folgende Diagramm:
\[\xymatrix{
q^*(E\vert_{B'}) \ar[d]^{\tilde{f}} \ar[rr]^t && F_b\times S^n_+ \arinclinv[d]\\
Q^*E \ar[rr]^T && F_b\times D^{n+1}
}\]
Gemäß der Verträglichkeit mit Trivialisierungen ist dann die Abbildung $T$ einfach, und die Abbildung $t$ ist wenigstens einfach modulo $\ker \tilde{f}_*$. Die rechte vertikale Abbildung ist einfach nach der Produktformel. Also ist nach der Kettenregel auch $\tilde{f}$ einfach, und (i) ist bewiesen.

(ii) Wir zeigen: \textit{Lemma \ref{tau:lemma7} ist richtig, falls $f:\,B'\lto B$ eine elementarer Retrakt ist. } Seien dazu $g:\,B\lto B'$ eine elementare Expansion, deren Homotopie-Inverses $f$ ist, und $H:\,B\times I\lto B$ eine Homotopie zwischen $f\circ g$ und $\id_B$. Dann ist die Abbildung $\bar{f}:\,f^*E\lto E$ homotop zur Komposition
\[\omega_H\circ\bar{g}^{-1}:\,f^*E\lto g^*f^*E \lto E.\]
Es wurde bereits bewiesen, dass $\bar{g}$ (und damit auch dessen Homotopie-Inverses) einfach modulo $\bar{f}_*^{-1}(N)$ ist. Die Abbidung $\omega_H$ ist als Faserhomotopie-Äquivalenz gemäß Lemma \ref{tau:lemma8} einfach modulo $N$. Daraus folgt die Behauptung (ii).

(iii) Wir zeigen nun den allgemeinen Fall. Eine einfache Homotopie-Äquivalenz $f$ zwischen endlichen CW-Komplexen ist homotop zu einer Komposition $g=g_1\circ\ldots\circ g_n$ von elementaren Expansionen und Retrakten. Die homotope Abänderung von $f:\,B'\lto B$ in die Abbildung $g$ induziert indes eine Faserhomotopie-Äquivalenz $\varphi:\,f^*E\lto g^*E$, und es gilt $\bar{g}\circ\varphi\simeq \bar{f}$. Mit Lemma \ref{tau:lemma8} erkennt man, dass $\bar{f}$ genau dann einfach modulo $N$ ist, wenn auch $\bar{g}$ einfach modulo $N$ ist. Unter der Identifikation $g^*E=g_n^*\dots g_1^*E$ gilt aber $\bar{g}=\bar{g}_1\circ\ldots\circ \bar{g}_n$, und es wurde bereits gezeigt, dass jede der Abbildungen $\bar{g}_i$ einfach modulo $\bar{g}_{i-1*}^{-1}\ldots \bar{g}_{1*}^{-1}(N)$ ist. Also ist $\bar{g}$ einfach modulo $N$.
\end{proof}

\begin{prop}\label{tau:prop1}
Die einfache Struktur $\xi(b,\zeta)$ modulo $N$ auf $E$ ist unabhängig von der Wahl der endlichen CW-Struktur auf $B$.
\end{prop}

\begin{proof}[Beweis]
Seien $B'$ und $B$ dieselben topologischen Räume, jedoch mit anderer endlicher CW-Struktur. Gemäß der topologischen Invarianz der Whitehead-Torsion ist die Identität $B'\lto B$ einfach. Die Anwendung von Lemma \ref{tau:lemma7} auf das Pull-back
\[\xymatrix{
E' \ar[d]^p \ar[rr]^{\id} && E \ar[d]^p \\
B' \ar[rr]^{\id} && B 
} \]
zeigt nun die Behauptung.
\end{proof}

\section{Das Element $\tau(f)$}

Nun können wir das Element $\tau(f)$ definieren. Sei dazu $f:\,M\lto B$ eine (stetige) Abbildung zwischen differenzierbaren geschlossenen und zusammenhängenden Mannigfaltigkeiten, die die Voraussetzungen an die Definition des Elements $\theta(f)$ erfüllt. Es gelte $\theta(f)=0$. Um keine Wahlen bezüglich einfacher Strukturen treffen zu müssen, nehmen wir außerdem an, dass eine der beiden Bedingungen erfüllt ist:
\begin{itemize}\smallsep
\item Die Eulercharakteristik von $B$ ist 0, oder
\item Ist $F\lto M^f\lto B$ die assoziierte Faserung, so ist die von der Inklusion induzierte Abbildung $\Whp{F}\lto \Whp{M^f}$ die Nullabbildung. Dies ist etwa dann erfüllt, wenn $\Whp{F}=0$ ist, oder wenn die von $f$ induzierte Abbildung $f_*:\,\pi_1(M)\lto\pi_1(B)$ injektiv ist.
\end{itemize}
Satz \ref{tau:satz2} zeigt, dass dann $M^f$ eine kanonische einfache Struktur hat, die unabhängig von der Wahl einer einfachen Struktur auf der Faser $F$ ist. Die differenzierbare Mannigfaltigkeit $M$ besitzt eine Triangulierung, die $M$ mit der Struktur eines endlichen CW-Komplexes und damit auch mit einer natürlichen einfachen Struktur ausstattet (vgl.~\cite{munkres}, Seite 101). Für die kanonische Inklusion $\lambda:\,M\lto M^f$ definieren wir nun
\[\tau(f):=\lambda_*^{-1}\,\tau(\lambda)\in\Whp{M}.\]

\begin{bem}
Das Element $\tau(f)$ hängt nur von der Homotopieklasse von $f$ in $[M,B]$ ab.
\end{bem}

\begin{proof}[Beweis]
Sei $g:\,M\lto B$ homotop zu $f$. Nach Proposition \ref{fas:prop3} existiert ein bis auf Homotopie kommutatives Diagramm
\[\xymatrix{
& M \ar[ld]_{\lambda_1} \ar[rd]^{\lambda_2} \\
M^f \ar[rr]^{\varphi} && M^g
}\]
mit einer Faserhomotopie-Äquivalenz $\varphi$, wobei $\lambda_1$ und $\lambda_2$ die jeweiligen kanonischen Inklusionen bezeichnen. In Folge von Lemma \ref{tau:lemma8} ist die Abbildung $\varphi$ einfach. Unter Ausnutzung der Kettenregel für die Whitehead-Torsion folgt die Behauptung.
\end{proof}

Um den Beweis von Proposition \ref{tau:prop_haupt} abzuschließen, ist noch zu zeigen, dass für ein differenzierbares Faserbündel $f$ das Element $\tau(f)$ verschwindet. Wir stellen eine weitere Bemerkung voran:

Ist $f:\,M\lto B$ eine Faserung, so ist ist diese nach Proposition \ref{fas:prop2} über die kanonische Inklusion $\lambda:\,M\lto M^f$ faserhomotopie-äquivalent zur assoziierten Faserung $p:\,M^f\lto B$. Sind also die Voraussetzungen zur Definition von $\tau(f)$ erfüllt, so erhält neben $M^f$ auch $M$ durch die Konstruktion aus Unterkapitel \ref{tau:abschn1} eine wohldefinierte einfache Struktur $\xi$. Gemäß der Verträglichkeit mit Faserhomotopie-Äquivalenzen ist dabei die kanonische Inklusion $\lambda:\,(M,\xi)\lto M^f$ einfach. Es folgt:

\begin{bem}
Für eine Faserung $f:\,M\lto B$ gilt:
\[\tau(f) = \tau\bigl[(M,[\id])\overset{\id}{\lto} (M,\xi)\bigr].\]
\end{bem}

Ist $f:\,M\lto B$ ein differenzierbares Faserbündel, müssen wir also die beiden auf $M$ gegebenen einfachen Strukturen $[\id]$ und $\xi$ vergleichen. Im Falle, dass das Faserbündel trivial ist, also $M\cong B\times N$ über $B$, so besagt Lemma \ref{tau:lemma9}, dass die beiden einfachen Strukturen übereinstimmen. Man beachte hierzu, dass eine Triangulierung einer Mannigfaltigkeit dieser die Struktur eines guten CW-Komplexes im Sinne von Lemma \ref{tau:lemma9} gibt. Insbesondere ist die Behauptung also für zusammenziehbare Basisräume richtig.

Der allgemeine Fall folgt unter Verwendung einer Henkelzerlegung der $m$-dimensiona\-len differenzierbaren Mannigfaltigkeit $B$ (siehe etwa~\cite{matsumoto}, Seite 79ff.) Dies ist eine Kette Räumen
\[D^m=B_0\subset B_1\subset \ldots \subset B_n=B,\]
sodass für alle $j\in\{0,\ldots,n-1\}$ der Raum $B_j$ eine Mannigfaltigkeit mit Rand ist und eine differenzierbare Einbettung $\phi_j:\,S^k\times D^l\lto \partial B_j$ zusammen mit einem (topologischen) Push-out 
\begin{equation}\begin{split}\label{tau:push-out3}\xymatrix{
S^k\times D^l \arinclinv[d] \arincl[rr]^{\phi_j} && B_j \arinclinv[d]\\
D^{k+1}\times D^l \arincl[rr] && B_{j+1}
}\end{split}\end{equation}
existiert. Die Zahlen $k$ und $l$ können dabei von $j$ abhängen, es gilt jedoch $k+l=n-1$. Außerdem können wir davon ausgehen, dass alle Mannigfaltigkeiten $B_j$ zusammenhängend sind (vgl.~\cite{matsumoto}, Seite 128f.). 

Nun sind die jeweiligen Urbilder $E_j:=p^{-1}(B_j)$ ebenfalls Mannigfaltigkeiten, und die Einschränkung von $p$ auf $E_j$ ist ebenfalls ein Faserbündel $p:\,E_j\lto B_j$. Der Raum $E_n$ hat also ebenfalls zwei einfache Strukturen, nämlich die natürliche, die von einer Triangulierung herrührt, sowie nach Wahl eines Punktes $b\in B_0$ die einfache Struktur $\xi_j=\xi(b,[\id])$, die der Raum aus der Konstruktion in Unterkapitel \ref{tau:abschn1} erhält. (Die Wahl des Basispunktes spielt nach Lemma \ref{tau:lemma2} dabei in Wirklichkeit keine Rolle, da alle Fasertransporte in einem Faserbündel als Homöomorphismen einfach sind.) Wir zeigen nun induktiv über $j$, dass diese beiden einfachen Strukturen auf $E_j$ übereinstimmen. Der Induktionsanfang wurde bereits bewiesen, denn $B_0=D^m$ ist kontrahierbar.

Nehmen wir als Induktionsvoraussetzung an, dass die Behauptung für ein $j\leq n-1$ gültig ist. Wir behaupten zunächst, dass für geeignete Triangulierungen der Räume im Push-out-Diagramm \eqref{tau:push-out3} dieses Diagramm ein Push-out von CW-Komplexen im Sinne von Lemma \ref{tau:lemma10} ist. Dazu zeigen wir: Die Räume in Diagramm \eqref{tau:push-out3} besitzen Triangulierungen, sodass die Spalten CW-Paare sind und die Abbildung $\phi_j$ zellulär ist.

Wir beginnen mit einer glatten Triangulierung von $D^{k+1}\times D^l$, die $S^k\times D^l$ und $S^k\times S^{l-1}$ als Unterkomplexe enthält. Nach \cite{munkres}, Seite 101, kann für eine glatte Mannigfaltigkeit mit Rand jede glatte Triangulierung des Randes zu einer glatten Triangulierung der gesamten Mannigfaltigkeit ausgedehnt werden. Erweitere also die gegebene Triangulierung von $\phi_j(S^k\times S^{l-1})$ zu einer glatten Triangulierung von $\partial B_j - \opn{int}\phi_j(S^k\times D^l)$. Durch Zusammensetzen erhält man so eine glatte Triangulierung von $\partial B_j$, bezüglich der die Abbildung $\phi_j$ zellulär ist. Diese Triangulierung kann nun zu einer Triangulierung von $B_j$ ausgedehnt werden, und man erhält wiederum durch Zusammensetzen eine Triangulierung von $B_{j+1}$, die $B_j$ als Unterkomplex enthält.

Betrachte nun das folgende Push-out:
\[\xymatrix{
E\vert_{S^k\times D^l} \arincl[rr] \arinclinv[d] && E_j \arinclinv[d] \\
E\vert_{D^{k+1}\times D^l} \arincl[rr] && E_{j+1}
}\]
Zunächst ist es wieder möglich, alle Räume so zu triangulieren, dass das Push-out zellulär ist, indem dasselbe Verfahren wie für Diagramm \eqref{tau:push-out3} angewendet wird: Ausgehend von einer glatten Triangulierung von $E\vert_{D^{k+1}\times D^l} \cong F\times D^{k+1}\times D^l$, die $E\vert_{S^k\times D^l}$ als Unterkomplex enthält, erhält man eine glatte Triangulierung von $\partial E_j = E_j\vert_{\partial B_j}$, die man auf $E_j$ ausdehnt; das Push-out selbst definiert dann eine Triangulierung auf $E_{j+1}$. Damit ist die von den natürlichen einfachen Strukturen auf $E_j$, $E\vert_{S^k\times D^l}$ und $E\vert_{D^{k+1}\times D^l}$ durch das Push-out auf $E_{j+1}$ induzierte einfache Struktur gleich der natürlichen.

Unter Anwendung der Verträglichkeit mit Push-outs (Lemma \ref{tau:lemma10}) erkennt man weiter, dass dieses Push-out einfach ist, wenn alle Räume ihre zu $(b,[\id])$ gehörenden kanonischen einfachen Strukturen tragen. Nach Induktionsvoraussetzung ist die zu $(b,[\id])$ gehörende kanonische einfache Struktur auf $E_j$ gleich der natürlichen. Die gleiche Behauptung ist richtig für die Räume $E\vert_{D^{k+1}\times D^l}$ und $E\vert_{S^k\times D^l}$, denn es handelt sich um triviale Faserbündel, für die wir diese Behauptung bereits verifiziert haben. Nach der Summenformel stimmen dann auch für $E_{j+1}$ die natürliche und die zu $(b,[\id])$ gehörende einfache Struktur überein, und der Induktionsschluss ist vollzogen. Wir haben also gezeigt:

\begin{prop}
Ist $f:\,E\lto B$ ein Faserbündel, und sind die Voraussetzungen an die Definition von $\tau(f)$ erfüllt, so gilt $\tau(f)=0$.
\end{prop}


\chapter{Basisraum $S^1$}\label{far}

Für den Fall einer Abbildung $f:\,M\lto S^1$ einer geschlossenen und zusammenhängenden differenzierbaren Mannigfaltigkeit der Dimension $n\geq 6$ in die 1-Sphäre hat F.~T.~Farrell in \cite{farrell} bereits die Frage untersucht, ob $f$ homotop zu einem Faserbündel ist. Als Resultat erhielt er --- unter der Annahme, dass die von $f$ induzierte Abbildung $\pi_1(M)\lto\pi_1(S^1)$ surjektiv ist --- drei Bedingungen, die zusammen notwendig und hinreichend dafür sind, diese Frage positiv zu beantworten. Die erste Bedingung ist eine Endlichkeitsbedingung an einen Raum $X$, der, wie wir zeigen werden, homotopie-äquivalent zur Faser der zu $f$ assoziierten Faserung ist. Die beiden weiteren Bedingungen sind jeweils das Verschwinden von $K$-theoretisch definierten Hindernissen. 

In diesem Kapitel sollen nun die von uns in den vorangegangenen Kapiteln erhaltenen Hindernisse $\theta(f)$ und $\tau(f)$ für den Basisraum $S^1$ spezialisiert und mit den von Farrell definierten Hindernissen verglichen werden. Tatsächlich wird sich zeigen, dass unser im vorigen Kapitel betrachtetes Element die Informationen der beiden Hindernisse von Farrell zusammenfasst. Insbesondere sind also unsere Elemente $\theta(f)$ und $\tau(f)$ die einzigen Hindernisse, die einer homotopen Abänderung von $f$ in ein differenzierbares Faserbündel entgegenstehen.

\section{Einige Notationen}\label{s1:abschn1}

Sei $p_f:\,M^f\lto S^1$ die assoziierte Faserung zur Abbildung $f$; wir identifizieren dabei die Faser $F$ mit der Faser über dem Basispunkt der punktierten $S^1$. Es sei vorausgesetzt, dass eine (endliche) einfache Struktur auf der Faser $F$ existiert, und sei ein Vertreter $z:\,Z\lto F$ dieser einfachen Struktur fest gewählt. Der in Kapitel \ref{theta} konstruierte Homomorphismus $\theta(f):\,\pi_1(S^1)\lto \Whp{M}$ ist bereits vollständig durch das Bild eines Erzeugers bestimmt. Sei $T$ der zur gegebenen Orientierung von $S^1$ gehörende Erzeuger von $\pi_1(S^1)$; wir setzen
\[\tau_1 := \theta(f)(T) \in\Whp{M}.\]

Unter der Voraussetzung, dass die Torsion $\tau_1$ verschwindet, können wir den Raum $M^f$ mit einer kanonischen einfachen Struktur ausstatten; diese ist unabhängig von der Wahl der einfachen Struktur auf $F$, weil die Eulercharakteristik der 1-Sphäre verschwindet. Wir werden im Folgenden stets die CW-Struktur von $S^1$ zu Grunde legen, die aus genau einer 0-Zelle und genau einer 1-Zelle besteht. Sei $\sigma:\,I\lto S^1$ die dazu gehörende anklebende Abbildung der 1-Zelle an die 0-Zelle, die die Enden des Einheitsintervalls verklebt. Vollzieht man die Konstruktion der kanonischen einfachen Struktur auf dem Totalraum der Faserung $p_f$ nach, erkennt man, dass ein Vertreter derselben von einer natürlichen Transformation von Push-out-Diagrammen des folgenden Typs induziert wird:

\begin{equation}\begin{split}\label{s1:diag3}\xymatrix{
Z\times I \ar[d]_{\simeq}^{\varphi} && Z\amalg Z \ar[d]^{z\amalg (T\circ z)}_{\simeq} \arinclinv[ll] \arincl[rr] && {\Cyl(\id\amalg T')} \ar[d]_{\simeq}^{\psi}\\
\sigma^*M^f && F\amalg F \arinclinv[ll] \ar[rr]^{\id\amalg\id} && F
}\end{split}\end{equation}

Hierbei sind $T:\,F\lto F$ ein Fasertransport entlang des Erzeugers $T\in\pi_1(S^1)$ und $T':\,Z\lto Z$ eine zelluläre Approximation der Abbildung $z^{-1}\circ T\circ z:\,Z\lto Z$. Die Abbildung $\varphi$ repräsentiert eine Trivialisierung der Faserung $\sigma^*M^f\lto I$; die Abbildung $\psi$ ist Vertreter der einfachen Struktur auf $F$. Der Urbildraum von $\psi$ ist hierbei der Abbildungszylinder der Abbildung $\id\amalg T':\,Z\amalg Z\lto Z$. 

Das Push-out der oberen Zeile in \eqref{s1:diag3} ist offenbar homöomorph zum Abbildungstorus
\[\Torus(Z,T') = Z\times I\; / \;\bigl\langle (z,0)\sim (T'(z),1) \;\bigl\vert\; z\in Z \bigr\rangle,\]
wobei man den Raum $Z\times I$ aus \eqref{s1:diag3} mit dem Teilraum $Z\times [\frac14,\frac34]\subset\Torus(Z,T')$ identifizieren kann. Wir notieren $x:\,\Torus(Z,T')\lto M^f$ die von Diagramm \eqref{s1:diag3} induzierte Abbildung. Weiter definiert die Abbildung $\sigma\circ\Proj:\,Z\times I\lto I\lto S^1$ zusammen mit der Punktabbildung $\Cyl(\id\amalg T')\lto \{*\}\subset S^1$ eine Abbildung $y:\,\Torus(Z,T')\lto S^1$, die das Diagramm
\[\xymatrix{
\Torus(Z,T') \ar[rr]^(0.6)x_(0.6){\simeq} \ar[rd]_y && M^f \ar[ld]^{p_f}\\
& S^1
}\]
kommutativ macht. Schließlich setzen wir unter Verwendung der kanonischen Inklusion $\lambda:\,M\lto M^f$:
\[\tau_2 = \lambda_*^{-1}\, x_*\,\tau(x^{-1}\circ\lambda)\in\Whp{M}.\]
Aus den Ergebnissen der vorherigen Kapitel folgt:

\begin{prop}
Die Torsionen $\tau_1$ und $\tau_2$ sind wohldefiniert und hängen nur von der Homotopieklasse von $f$ in $[M,S^1]$ ab.
\end{prop}

Wollen wir untersuchen, ob $f$ homotop zu einem differenzierbaren Faserbündel ist, so können wir ohne Beschränkung der Allgemeinheit davon ausgehen, dass die von $f$ induzierte Abbildung $f_*:\,\pi_1(M)\lto\pi_1(S^1)$ surjektiv ist. Denn ist das Bild von $f_*$ in $\pi_1(S^1)=\ZZ$ eine Untergruppe $n\ZZ$, so existiert ein Lift $g$ der Abbildung $f$ zur $n$-fachen Überlagerung:
\[\xymatrix{
&&S^1\ar[d]&z \ar@{|->}[d]\\
M \ar[rr]^f \ar@{.>}[rru]^g && S^1 & z^n
}\]
Sodann ist die von $g$ induzierte Abbildung der Fundamentalgruppen eine Surjektion, und $g$ ist genau dann homotop zu einem differenzierbaren Faserbündel, wenn auch $f$ es ist. Ist das Bild von $f_*$ in $\pi_1(S^1)$ trivial, so ist $f$ mit Sicherheit nicht homotop zu einem Faserbündel: In diesem Fall existiert ein Lift $g$ in die universelle Überlagerung $\RR$ von $S^1$; eine Abbildung einer kompakten Mannigfaltigkeit nach $\RR$ kann aber kein Faserbündel sein.

\section{Wechsel des Basisraums und Farrells erste Bedingung}

Es wird im Folgenden oft vorteilhaft sein, die gegebenen Abbildungen nach $S^1$ zu Abbildungen nach $\RR$ zurückzuziehen. Dieser Wechsel des Basisraums soll hier näher betrachtet werden.

Sei $f:\,M \lto S^1$ eine Abbildung von punktierten Räumen; es sei $f_*:\,\pi_1(M)\lto \pi_1(S^1)$ surjektiv. Es bezeichne $X$ das Pull-back von $f$ mit der Standard-Abbildung $\pi:\,(\RR,0)\lto S^1$:
\begin{equation}\begin{split}\label{s1:defn_X}\xymatrix{
X \ar[d]^{\bar{f}} \ar[rr]^p && M\ar[d]^f\\
{\RR} \ar[rr]^{\pi}  && S^1
}\end{split}\end{equation}

\begin{bem}\label{s1:basiswechsel}
Die Abbildung $p:\,X\lto M$ ist eine Überlagerung, und die Sequenz
\[\xymatrix{
1 \ar[r] & \pi_1(X) \ar[r]^{p_*} & \pi_1(M) \ar[r]^{f_*} & \pi_1(S^1) \ar[r] & 1
}\]
ist exakt. 
\end{bem}

\begin{proof}[Beweis]
Ein Pull-back einer Überlagerung ist wieder eine Überlagerung. Nach Anwenden des Funktors $\pi_1$ auf das Diagramm (\ref{s1:defn_X}) erkennt man wegen $\pi_1(\RR) = 0$, dass $f_*\circ p_* = 0$. Sei umgekehrt ein geschlossener Weg $\gamma$ in $M$ gegeben mit $f_*([\gamma])=0\in\pi_1(S^1)$. Dann existiert eine Hochhebung von $\gamma$ zu einem geschlossenen Weg $\gamma'$ in $\RR$, und $\gamma$ zusammen mit $\gamma'$ definieren auf Grund der universellen Eigenschaft des Pull-backs einen geschlossenen Weg $\gamma''$ in $X$. Nach Konstruktion gilt $p_*([\gamma''])=[\gamma]\in\pi_1(M)$.
\end{proof}

Aus der Exaktheit der obigen Sequenz folgt nach Wahl eines Schnittes von $f_*$ die Existenz eines Isomorphismus
$$\pi_1(M) \cong G \rtimes_\alpha \ZZ$$
der Fundamentalgruppe von $M$ mit einem semi-direkten Produkt aus $G := \opn{ker}f_*$ und der unendlichen zyklischen Gruppe. Insbesondere operiert $\ZZ$ auf $G$ durch Konjugation. Sei $T$ der orientierte Erzeuger der Fundamentalgruppe $\pi_1(S^1)=\ZZ$. Die Operation von $T$ auf $G$ ist ein Gruppenhomomorphismus 
\[\alpha:\,G\lto G,\quad g\mapsto TgT^{-1}\] 

Betrachten wir als Anwendung den Fall, dass $M=\Torus(Z,T)$ der Abbildungszylinder einer zellulären Selbstabbildung eines endlichen zusammenhängenden CW-Komplexes ist, wobei $f$ nun kanonische Abbildung 
\[\Torus(Z,T) \lto S^1,\quad [(z,t)]\mapsto \pi(t) \mathrm{\quad mit\;}(z,t)\in Z\times I\]
ist. $M$ ist ein punktierter Raum nach Wahl eines Basispunkts in $Z=Z\times\{0\}$. Das Pull-back $X$ entsteht dann aus einer durch $\ZZ$ indizierten disjunkten Vereinigung von Abbildungszylindern von $T$, indem man den Boden des $n$-ten Zylinders mit dem Deckel der $(n+1)$-ten Zylinders identifiziert. Wie oben notieren wir $\bar{f}:\,X\lto\RR$ die Strukturabbildung des Pull-backs.

\begin{lem}\label{far:lemma1}
Identifiziere $Z$ mit dem Teilraum $Z\times\{0\}\subset X$. Ist $T$ eine Homotopie-Äquivalenz, so induziert auch die Inklusion $Z\lto X$ eine Homotopie-Äquivalenz.
\end{lem}

\begin{proof}[Beweis]
Nach einem wohlbekannten Satz (vgl.~etwa \cite{tomdieck}, Seite 258) genügt es, die folgenden beiden Bedingungen zu verifizieren: (a) Die von der Inklusion induzierte Abbildung $\pi_1(Z)\lto\pi_1(X)$ ist ein Isomorphismus, und (b) fasst man die universelle Überlagerung $\tilde{Z}$ von $Z$ als Teilraum der universellen Überlagerung $\tilde{X}$ von $X$ auf, so gilt $H_*(\tilde{X},\tilde{Z})=0$. 

Für $n\in\ZZ$ ist $\bar{f}^{-1}([n,n+1])$ gleich dem $n$-ten Abbildungszylinder in $X$. Da $T$ eine Homotopie-Äquivalenz ist, sind die Inklusionen des Deckels und des Bodens in den Abbildungszylinder von $T$ Deformationsretrakte. Induktiv sieht man, dass die Inklusionen $Z\lto \bar{f}^{-1}([-n,n])$ für alle natürlichen Zahlen $n$ Deformationsretrakte sind, insbesondere induzieren diese jeweils Isomorphismen der Fundamentalgruppen. 

Um Behauptung (a) einzusehen, sei nun $\gamma:\,S^1\lto X$ ein Vertreter eines Elements von $\pi_1(X)$. Wegen der Kompaktheit von $S^1$ ist der Bildbereich von $\bar{f}\circ\gamma\subset\RR$ beschränkt, sodass $\bar{f}\circ\gamma$ ein Element von $\pi_1(\bar{f}^{-1}([-n,n]))$ für ein hinreichend großes $n$ definiert. Nach dem eingangs Gesagten ist dann $\gamma$ homotop relativ Basispunkt zu einer Schleife in $Z$. Die so definierte Abbildung $\pi_1(X)\lto\pi_1(Z)$ ist ein Inverses der von der Inklusion induzierten Abbildung.

Behauptung (b) gilt, weil $\tilde{X}$ der direkte Limes der Urbilder von $\bar{f}^{-1}([-n,n])$ in $\tilde{X}$ für $n\to\infty$ ist.
\end{proof}

Um auf die in Unterkapitel \ref{s1:abschn1} betrachtete Situation zurückzukommen, existiert also eine kurze exakte Sequenz
\begin{equation}\label{s1:sequenz1}\xymatrix{
1 \ar[r] & \pi_1(Z) \ar[r]^(0.35){i_*} & \pi_1(\Torus(Z,T')) \ar[r]^(0.65){y_*} & \pi_1(S^1) \ar[r] & 1,
}\end{equation}
deren Abbildung $i_*$ von der Inklusion induziert wird. Dazu beachte man, dass die Abbildung $y$ homotop zur soeben betrachteten Abbildung $f:\,\Torus(Z,T')\lto S^1$ ist.

Wir haben auch die Auswirkungen des Basiswechsels auf die assoziierte Faserung von $f$ zu untersuchen. Es wird sich zeigen, dass die Bildung der assoziierten Faserung mit der Bildung des Pull-backs unter $\pi:\,\RR\lto S^1$ verträglich ist. 

Betrachte dazu wieder das Pull-back von $f$ mit $\pi$ wie im Diagramm \eqref{s1:defn_X}. Seien $p_f:\,M^f \lto S^1$ bzw.~$p_{\bar{f}}:\,X^{\bar{f}}\lto \RR$ die assoziierten Faserungen zu $f$ bzw.~$\bar{f}$. Es bezeiche $\pi^I:\,\RR^I\lto (S^1)^I$ die Abbildung, die $\gamma:\,I\lto\RR$ auf $\pi\circ\gamma$ abbildet. Weiter seien $\alpha:\,M^f\lto (S^1)^I$ und $\beta:\,M^f\lto M$ die kanonischen Abbildungen der assoziierten Faserung $M^f$, die gestrichten Varianten die der Faserung $X^{\bar{f}}$. Das Diagramm
\[\xymatrix{
& M^f \ar[rrr]^{\alpha} \ar'[d]^(0.8){\beta}[dd] & & & (S^1)^I \ar[dd]^{p_0}\\
X^{\bar{f}} \ar@{.>}[ur]^{u} \ar[rrr]^(0.6){\bar{\alpha}} \ar[dd]^{\bar{\beta}} & & & {\RR^I} \ar[ru]^(0.4){\pi^I} \ar[dd]^(0.3){p_0}\\
& M \ar'[rr]^f[rrr] & & & S^1\\
X \ar[ru]^p \ar[rrr]^{\bar{f}} & & & {\RR} \ar[ru]^{\pi}
}\]
induziert dann eine Abbildung $u:\,X^{\bar{f}}\lto M^f$, die offenbar das Diagramm
\begin{equation}\begin{split}\label{s1:pull-back1}\xymatrix{
X^{\bar{f}} \ar[rr]^{u} \ar[d]^{p_{\bar{f}}} & & M^f \ar[d]^{p_f}\\
{\RR} \ar[rr]^{\pi} & & S^1
}\end{split}\end{equation}
kommutativ macht, d.~h.~eine Abbildung von Faserungen ist. Direkt aus der Konstruktion folgt, dass das Diagramm
\begin{equation}\begin{split}\label{s1:pull_back2}\xymatrix{
X \ar[rr]^p \arinclinv[d]^{\lambda} & & M \arinclinv[d]^{\lambda} \\
X^{\bar{f}} \ar[rr]^{u} & & M^f
}\end{split}\end{equation}
ebenfalls kommutiert.

\begin{lem}\label{s1:lemma2}
Das Diagramm (\ref{s1:pull-back1}) ist ein Pull-back.
\end{lem}

\begin{proof}[Beweis]
Wir zeigen die universelle Eigenschaft des Pull-backs. Seien $\varphi_1:\,Z\lto \RR$ und $\varphi_2:\,Z\lto M^f$ zwei Abbildungen mit $p_f\circ\varphi_2 = \pi\circ\varphi_1$. Es ist zu zeigen, dass es eine eindeutige Abbildung $\Phi:\,Z\lto X^{\bar{f}}$ gibt mit $u\circ\Phi = \varphi_2$ und $p_{\bar{f}}\circ\Phi = \varphi_1$. 

Aus den Voraussetzungen folgt, dass das Diagramm
\begin{equation}\begin{split}\label{s1:diag4}\xymatrix{
Z \ar[rr]^{\alpha\circ\varphi_2} \ar[d]^{\varphi_1} & & (S^1)^I \ar[d]^{p_1}\\
{\RR} \ar[rr]^{\pi} & & S^1
}\end{split}\end{equation}
kommutativ ist. Unter Ausnutzung des Exponentialgesetzes für kompakt erzeugte Räume ist Diagramm (\ref{s1:diag4}) äquivalent zum folgenden HHP für die Überlagerung $\pi$:
\begin{equation}\begin{split}\label{s1:diag2}\xymatrix{
Z\times\{1\} \ar[rr]^{\varphi_1} \arinclinv[d] & & {\RR} \ar[d]^{\pi}\\
Z\times I \ar@{.>}[rru]^{\Phi_1} \ar[rr]_{\alpha\circ\varphi_2} & & S^1
}\end{split}\end{equation}
Wir interpretieren die Lösung dieses HHP als Abbildung $\Phi_1:\,Z\lto \RR^I$. Weiter sei $\Phi_2:\,Z\lto X$ induziert von den Abbildungen $\beta\circ\varphi_2:\,Z\lto M$ und $p_0\circ\Phi_1:\,Z\lto \RR$. $\Phi_1$ und $\Phi_2$ induzieren dann zusammen eine Abbildung $\Phi:\,Z\lto X^{\bar{f}}$, von der man verifiziert, dass sie die geforderten Eigenschaften hat.

Es bleibt noch die Eindeutigkeit der Abbildung $\Phi$ zu zeigen. Sei dazu $\tilde{\Phi}$ eine zweite Abbildung mit $u\circ\tilde{\Phi} = \varphi_2$ und $p_{\bar{f}}\circ\tilde{\Phi} = \varphi_1$. Dann ist $\bar{\alpha}\circ\tilde{\Phi}:\,Z\lto \RR^I$ eine Lösung des HHP (\ref{s1:diag2}), dessen Lösung aber eindeutig ist, sodass gilt: $\bar{\alpha}\circ\tilde{\Phi}=\bar{\alpha}\circ\Phi$. Dann ist jedoch 
\[ \bar{f}\circ\bar{\beta}\circ\tilde{\Phi} = p_0\circ\bar{\alpha}\circ\tilde{\Phi} = p_0\circ\bar{\alpha}\circ\Phi = \bar{f}\circ\bar{\beta}\circ\Phi\]
und
\[p\circ\bar{\beta}\circ\tilde{\Phi} = \beta\circ u\circ\tilde{\Phi}= \beta\circ u\circ\Phi = p\circ\bar{\beta}\circ\Phi,\]
sodass mit der universellen Eigenschaft des Pull-backs auch $\bar{\beta}\circ\tilde{\Phi}=\bar{\beta}\circ\Phi$ gilt. Insgesamt folgt $\tilde{\Phi}=\Phi$.
\end{proof}

\begin{bem}
Offenbar besitzt Lemma \ref{s1:lemma2} eine Verallgemeinerung auf den Fall einer beliebigen Abbildung $f:\,X\lto Y$ und einer Überlagerung $\pi:\,\tilde{Y}\lto Y$.
\end{bem}

Wir haben bereits gesehen, dass man die Strukturabbildung $p:\,X\lto M$ auch als die zu $G$ gehörende Überlagerung von $M$ betrachten kann, wobei $G$ der Kern der Abbildung $f_*:\,\pi_1(M)\lto\pi_1(S^1)$ sei. Lemma \ref{s1:lemma2} besagt, dass man $X^{\bar{f}}$ auch als Pull-back von $M^f$ mit $\pi$ auffassen kann. Damit ist im Sinne von Bemerkung \ref{s1:basiswechsel} die kanonische Abbildung 
\[u:\,X^{\bar{f}} \lto M^f\] 
die zum Kern der Abbildung $p_{f*}:\,\pi_1(M^f)\lto\pi_1(S^1)$ gehörende Überlagerung. In dieser Sichtweise bedeutet die Kommutativität des Diagramms (\ref{s1:pull_back2}):

\begin{bem}\label{s1:bem1}
Die kanonische Inklusion $\lambda:\,X\lto X^{\bar{f}}$ überdeckt die kanonische Inklusion $\lambda:\,M\lto M^f$.
\end{bem}

Als Konsequenz erhält man die folgende natürliche Transformation von exakten Sequenzen, deren unterste Zeile die Sequenz \eqref{s1:sequenz1} ist und in der die Abbildung $j_*$ von der Inklusion $F\lto X^{\bar{f}}$ induziert wird:
\[\xymatrix{
1 \ar[r] & \pi_1(X) \ar[r]^{p_*} \ar[d]^{\lambda_*} & \pi_1(M) \ar[r]^{f_*} \ar[d]^{\lambda_*} & \pi_1(S^1) \ar[r] \ar@{=}[d] & 1\\
1 \ar[r] & \pi_1(X^{\bar{f}}) \ar[r]^{u_*} & \pi_1(M^f) \ar[r]^{p_{f*}} & \pi_1(S^1) \ar[r] & 1\\
1 \ar[r] & \pi_1(Z) \ar[u]_{j_*\circ z_*} \ar[r]^(0.35){i_*} &\pi_1(\Torus(Z,T')) \ar[u]_{x_*}\ar[r]^(0.65){y_*} &\pi_1(S^1) \ar@{=}[u] \ar[r] &1
}\]
Dabei induziert ein Schnitt der Abbildung $f_*$ einen Schnitt der Abbildungen $p_{f*}$ und $y_*$. Wir erinnern daran, dass die Wahl von derartigen Schnitten eine Struktur von $\pi_1(M)$, $\pi_1(M^f)$ bzw.~$\pi_1(\Torus(Z,T')$ als semi-direktes Produkt induziert.

\begin{kor}
Die Abbildungen $\lambda_*:\,\pi_1(M)\lto\pi_1(M^f)$ und $x_*:\,\pi_1(\Torus(Z,T'))\lto \pi_1(M^f)$ erhalten die semi-direkte Produktstruktur.
\end{kor}

Sei nun $M$ eine geschlossene und zusammenhängende differenzierbare Mannigfaltigkeit. Die erste in \cite{farrell} genannte Bedingung dafür, dass $f$ homotop zu einer Faserung ist, lautet nun, dass $X$ den Homotopietyp eines endlichen CW-Komplexes hat. Tatsächlich ist dies äquivalent zu unserer, an die Definition von $\tau_1$ geknüpfte Bedingung, dass die Faser der assoziierten Faserung den Homotopietyp eines endlichen CW-Komplexes hat, denn es gilt:

\begin{prop}\label{s1:prop1}
Die Faser $F$ ist homotopie-äquivalent zu $X$. 
\end{prop}

\begin{proof}[Beweis]
Die Faserung $p_{\bar{f}}:\,X^{\bar{f}}\lto \RR$ ist ein Pull-back von $p_f:\,M^f\lto S^1$, sodass deren Fasern identizifiert werden können. Nun ist $X^{\bar{f}}$ homotopie-äquivalent zur Faser $F$, da der Basisraum der Faserung zusammenziehbar ist. Schließlich ist die Inklusion $\lambda:\,X\lto X^{\bar{f}}$ ebenfalls eine Homotopie-Äquivalenz.
\end{proof}

\section{Splittings und Farrells zweite Bedingung}

An dieser Stelle seien einige Begriffe und Ergebnisse aus Farrells Arbeit \cite{farrell} kurz zusammengefasst. Die Seitenangaben beziehen sich dabei immer auf den genannten Artikel. Ist $b\in S^1$ ein regulärer Wert einer differenzierbaren Abbildung $f:\,M\lto S^1$, so ist $N:=f^{-1}(b)$ eine $(n-1)$-dimensionale Untermannigfaltigkeit von $M$ mit trivialem Normalenbündel $\nu$. Ein solches Paar $(N,\nu)$ heiße im Folgenden \emph{Splitting}. Die Pontrjagin-Thom-Konstruktion (siehe etwa \cite{ranicki}, Seite 128) etabliert 1:1-Korrespondenz zwischen gerahmten Bordismusklassen von Splittings und Homotopieklassen von Abbildungen $M\lto S^1$. Das erste Ziel ist es, durch geschickte Wahl elementarer gerahmter Bordismen die Eigenschaften eines gegebenen Splitting zu verbessern. 

\begin{prop}[Seite 325]
Induziert $f$ eine Surjektion der Fundamentalgruppen, und ist $X$ dominiert von einem endlichen CW-Komplex, so gilt: Ein Splitting $N$ kann so verbessert werden, dass $N$ zusammenhängend ist und die Inklusion $N\lto X$ einen Isomorphismus der Fundamentalgruppen induziert.
\end{prop}

Wir betrachten nun zu diesem verbesserten Splitting das Pullback von $M$ zu einem Raum $X$ über $\RR$ und bezeichnen mit $T:\,X\lto X$ die zum orientierten Erzeuger $T$ von $\pi_1(S^1)$ gehörende Decktransformation. Es sei eine Hochhebung von $N$ zu einem Teilraum von $X$ fest gewählt; diese Hochhebung werden wir ebenfalls mit $N$ bezeichnen. Der Raum $X-N$ zerfällt in zwei zusammenhängende Komponenten $A$ und $B$, wobei $A$ diejenige Komponente sei mit $A\subset T(A)$. Dann induzieren auch die Inklusionen  $A\lto X$ und $B\lto X$ Isomorphismen der Fundamentalgruppen; mittels dieser Isomorphismen identifizieren wir die Fundamentalgruppen miteinander und nennen sie $G=\pi_1(X)$. 

Zur universellen Überlagerung $\tilde{X}\lto X$ bezeichne $\tilde{N}$, $\tilde{A}$ bzw.~$\tilde{B}$ das Urbild von $N$, $A$ bzw.~$B$ in $\tilde{X}$. Es sei $H_*(X,A;\ZZ[G]) = H_*(\tilde{X},\tilde{A})$ die zelluläre Homologie der universellen Überlagerungen als $\ZZ[G]$-Modul. Für eine natürliche Zahl $s$ heiße ein Splitting $N$ \emph{$s$-bi-zusammenhängend}, falls gilt: 
\begin{itemize}\setlength{\itemsep}{0pt}
\item $N$ ist zusammenhängend,
\item Die Inklusion $N\lto X$ induziert einen Isomorphismus der Fundamentalgruppen, 
\item $H_j(X,A;\ZZ[G]) = 0$ für $j\leq s$,
\item $H_j(X,B;\ZZ[G]) = 0$ für $j\leq n - s - 1$.
\end{itemize} 

\begin{prop}[Seite 329]
Ein Splitting $N$ kann so verbessert werden, dass es 2-bi-zusammenhängend ist.
\end{prop}

Zu einem gegebenen Splitting $N$ bezeichne mit $M_N$ die Mannigfaltigkeit mit Rand $\overline{T(A)-A}$. Mann kann $M_N$ auch als Bordismus zwischen $N$ und $T(N)$ auffassen. Farrell zeigt (Seite 329), dass ein Splitting $N$ genau dann $s$-bi-zusammenhängend ist, wenn der Bordismus $M_N$ eine Henkelzerlegung hat, die nur aus $s$- und $(s+1)$-dimensionalen Henkeln besteht. Insbesondere folgt, dass für ein $s$-bi-zusammenhängendes Splitting der Modul $H_*(X,A;\ZZ[G])$ im Grad $s+1$ konzentriert ist. Es gilt:

\begin{lem}[Seiten 326, 329]\label{far:lemma3}
Ist $N$ ein $s$-bi-zusammenhängendes Splitting, so ist $H_*(X,A;\ZZ[G])$ ein endlich erzeugter projektiver $\ZZ[G]$-Modul.
\end{lem}

Die Moduln $H_*(X,T^j(A);\ZZ[G])$ bilden für $j\in\NN$ unter den von den Inklusionen induzierten Abbildungen ein direktes System mit
\[\dirlim_{j\to\infty} H_*(X,T^j(A);\ZZ[G]) = H_*(X,X;\ZZ[G]) = 0.\]
Aus der endlichen Erzeugtheit der Moduln folgt damit, dass für ein hinreichend großes $j$ die von der Inklusion induzierte Abbildung $H_*(X,A;\ZZ[G])\lto H_*(X,T^j(A);\ZZ[G])$ die Nullabbildung ist. Dies hat als Konsequenz:

\begin{bem}\label{far:kurze_exakte_sequenzen}
Die kanonische Abbildung $H_*(X;\ZZ[G])\lto H_*(X,A;\ZZ[G])$ ist die Nullabbildung, und es existieren kurze exakte Sequenzen
$$\xymatrix{
0 \ar[r] & H_{*+1}(X,A;\ZZ[G]) \ar[r] & H_*(A;\ZZ[G]) \ar[r] & H_*(X;\ZZ[G]) \ar[r] & 0.
}$$
\end{bem}

Denn die Abbildung $H_*(X;\ZZ[G])\lto H_*(X,A;\ZZ[G])$ faktorisiert über die soeben betrachtete Abbildung $H_*(X,T^{-j}(A);\ZZ[G])\lto H_*(X,A;\ZZ[G])$, die für hinreichend großes $j$ die Nullabbildung ist. Die lange exakte Homologiesequenz des Paares $(\tilde{X},\tilde{A})$ zerfällt sodann in die genannten kurzen exakten Sequenzen. --- Für eine weitere unmittelbare Konsequenz betrachte den von einer Hochhebung von $T^{-1}: X\lto X$ auf die universelle Überlagerung induzierten $\ZZ$-linearen Endomorphismus $t_*^{-1}$ von $H_*(X,A;\ZZ[G])$. Dieser lässt sich als Komposition 
$$\xymatrix{
H_*(X,A;\ZZ[G])\ar[r] & H_*(X,T(A);\ZZ[G]) \ar[r]^(0.55){t_*^{-1}} & H_*(X,A;\ZZ[G])
}$$ 
schreiben, wobei die erste Abbildung von der Inklusion induziert ist, und es folgt:

\begin{bem}
$t_*^{-1}:\, H_*(X,A;\ZZ[G])\lto H_*(X,A;\ZZ[G])$ ist nilpotent.
\end{bem} 

Mit Hilfe eines $s$-bi-zusammenhängendes Splitting definiert Farrell nun in der folgenden Weise ein Hindernis $c(f)$ dafür, dass $f$ homotop zu einer Faserung ist. Wir beschreiben zunächst die Gruppe $C(\ZZ[G],\alpha)$, in dem dieses Hindernis definiert ist. Für einen Ring $R$ mit 1 und einen Automorphismus $\alpha$ von $R$ betrachte die Kategorie $\mathcal{C}(R,\alpha)$ mit folgenden Objekten und Morphismen: Objekte sind Paare $(P,f)$ von endlich erzeugten projektiven $R$-Rechtsmoduln $P$ und nilpotenten, $\alpha$-semilinearen Endomorphismen $f$ von $R$, d.~h.~es gilt $f(xr) = f(x)\alpha(r)$ für alle $r\in R, x\in P$. Morphismen $(P,f)\lto (P',f')$ in dieser Kategorie sind $R$-Modulhomomorphismen $\varphi:\,P\lto P'$ mit $f'\circ\varphi = \varphi\circ f$. Es handelt sich bei $\mathcal{C}(R,\alpha)$ um eine abelsche Kategorie. Der Raum $C(R,\alpha)$ besteht nun aus Isomorphieklassen von Objekten in $\mathcal{C}(R,\alpha)$, geteilt durch die Äquivalenzrelation, die von folgenden beiden Bedingungen erzeugt wird:
\begin{enumerate}\setlength{\itemsep}{0pt}
\item Für einen freien Modul $F$ ist $(P,f)\sim (P\oplus F,f\oplus 0)$
\item Für eine exakte Sequenz $0\lto (P',f')\lto (P,f)\lto (P'',f'')\lto 0$ in $\mathcal{C}(R,\alpha)$ gilt $(P,f)\sim (P'\oplus P'', f'\oplus f'')$.
\end{enumerate}
Bezüglich der Verknüpfung, die durch die direkte Summe von Moduln und Morphismen gegeben ist, ist dies in der Tat eine Gruppe (Seite 318).

In unserem Fall ist $R=\ZZ[G]$ der Gruppenring der Fundamentalgruppe; der Morphismus $\alpha$ ist gegeben durch die Aktion von $T\in \pi_1(S^1)$ auf $G=\pi_1(X)$ durch Konjugation. Das Hindernis $c(f)$ ist nun definiert als
\[(-1)^{s+1}\,\bigl[(H_{s+1}(X,A;\ZZ[G]),t_*^{-1})\bigr]\in C(\ZZ[G],\alpha).\] 
Dieses Element ist wohldefiniert und hängt lediglich ab von der Homotopieklasse von $f$ in $[M,S^1]$. Ferner gilt:

\begin{satz}[Seite 338] Die folgenden Bedingungen sind äquivalent:
\begin{enumerate}\smallsep
\item $c(f)=0$.
\item Es gibt ein Splitting $N'$, sodass der Bordismus $M_{N'}$ ein $h$-Kobordismus ist.
\end{enumerate}
\end{satz}

Dieser Satz ist eines der Hauptresultate aus der zitierten Arbeit von Farrell. Tatsächlich zeigt Farrell das etwas allgemeinere Resultat, dass man jeden Vertreter $(P,\phi)$ von $c(f)\in C(\ZZ[G],\alpha)$ "`als Splitting realisieren"' kann, d.~h.~dass es zu jedem solchen Vertreter ein $s$-bi-zusammenhängendes Splitting $N'$ von $f$ gibt, sodass für dieses Splitting ein Isomorphismus $(H_{s+1}(X,A';\ZZ[G]),t_*^{-1})\cong(P,\phi)$ existiert. Ist dann $c(f)=0$, so ist das Paar $(0,0)$ aus Nullraum und Nullabbildung ein Vertreter von $c(f)$, folglich existiert ein $s$-bi-zusammenhängendes Splitting $N'$ mit $H_*(X,A';\ZZ[G])=0$. Dies wiederum ist äquivalent zur Aussage, dass der entsprechende Bordismus $M_{N'}$ ein $h$-Kobordismus ist.

Zur Abbildung $f:\,M\lto S^1$ definieren wir die Abbildung $-f:\,M\lto S^1$ als die von $\pi\circ (-\bar{f}):\,X\lto\RR\lto S^1$ induzierte Abbildung. Ein Splitting für $f$ ist offenbar auch ein Splitting für $-f$, in dem die Rollen von $A$ und $B$ vertauscht sind, insbesondere ist ein $s$-bi-zusammenhängendes Splitting für $f$ ist auch ein $(n-s-1)$-bi-zusammenhängendes Splitting für $-f$. Dementsprechend gilt
\[c(-f)= (-1)^{n-s}\,\bigl[(H_{n-s}(X,B;\ZZ[G]),t_*)\bigr]\in C(\ZZ[G],\alpha^{-1}).\]
Aus dem oben zitierten Satz schließen wir, dass $c(f)$ genau dann verschwindet, wenn auch $c(-f)$ verschwindet. Der algebraische Übergang von $c(f)$ nach $c(-f)$ wird von Farrell untersucht, spielt für unsere Zwecke jedoch keine Rolle.

Der Zusammenhang zwischen $c(f)$, bzw.~genauer $c(-f)$, zur von uns definierten Torsion $\tau_2$ wird nun durch den von Farrell auf Seite 322 für eine Gruppe $G$ und einen Automorphismus $\alpha$ von $G$ definierten Gruppenhomomorphismus 
\[p:\,\Wh(G\rtimes_\alpha\ZZ)\lto C(\ZZ[G],\alpha^{-1})\] 
hergestellt. 

\begin{satz}\label{s1:satz2}
Für $\tau_1=0$ gilt:
\[c(-f) = p(\tau_2).\]
\end{satz}

Insbesondere impliziert das Verschwinden von $\tau_2$ auch das Verschwinden des Hindernisses $c(f)$. Der Beweis von Satz \ref{s1:satz2} gelingt mit Hilfe von Methoden, die Farrell selbst bereitstellt, und wird den Rest dieses Unterkapitels einnehmen.

Wir stellten bereits fest, dass das Diagramm
\begin{equation}\begin{split}\label{far:abb_ueber_S1}\xymatrix{
M \ar[rr]^{\lambda} \ar[rrd]_f & & M^f  \ar[d]^{p_f} & & \Torus(Z,T') \ar[ll]_x \ar[lld]^y \\
& & S^1
}\end{split}\end{equation}
kommutiert und dass man $\pi_1(M)$ und $\pi_1(\Torus(Z,T'))$ mittels der von $x^{-1}\circ\lambda$ induzierten Abbildung miteinander und mit $G\rtimes_\alpha \ZZ$ identifizieren kann, wobei $G=\pi_1(X)$. In dieser Notation ist $\tau_2$ definiert als die Whitehead-Torsion der Abbildung $x^{-1}\circ\lambda$.

Da sowohl $c(f)$ als auch $\tau_2$ nur von der Homotopieklasse von $f$ abhängen, können wir ohne Beschränkung der Allgemeinheit die folgende Situation annehmen: Es gibt ein $b\in S^1$, sowie eine abgeschlossene zusammenhängende Umgebung $U$ von $b$, die nicht den Basispunkt der $S^1$ enthält, sodass $f$ auf der Teilmenge $f^{-1}(U)\subset M$ differenzierbar und submersiv ist und $N:=f^{-1}(b)$ ein $s$-bi-zusammenhängendes Splitting von $M$ für $f$ ist. 

Dann ist $f^{-1}(U)$ über $U$ diffeomorph zu einem Produkt $U\times N$. Zu einer gegeben glatten Triangulierung von $N$ erhält man eine Produkt-Triangulierung von $U\times N$, sodass $\partial U\times N$ und $\{b\}\times N$ Unterkomplexe sind. Erweitere nun diese Triangulierung von $\partial U\times N$ zu einer glatten Triangulierung von $f^{-1}(\overline{S^1-U}\times N)$. Dies ist möglich nach \cite{munkres}, Seite 101. Auf diese Weise erhält man eine Triangulierung von $M$, die $N$ als Unterkomplex enthält.

Wir übernehmen zum Splitting $N$ die Bezeichnungen $X$, $A$, $B$ von oben. Man beachte, dass die gewählte Triangulierung von $M$ eine Triangulierung von $X$ induziert.

\begin{lem}\label{far:lemma2}
Für eine geeignete endliche CW-Struktur auf $\Torus(Z,T')$ kann $x^{-1}\circ\lambda$ homotop so durch eine zelluläre Abbildung $g:\,M\lto \Torus(Z,T')$ approximiert werden, dass gilt:
\begin{enumerate}\smallsep
\item Das Diagramm
\[\xymatrix{
f^{-1}(U) \ar[rr]^g \ar[rd]_f & & y^{-1}(U) \ar[ld]^y \\
& U
}\]
ist kommutativ.
\item $g$ bildet den Teilraum $f^{-1}(\overline{S^1-U})\subset M$ auf $y^{-1}(\overline{S^1-U})$ ab.
\end{enumerate}
\end{lem}

Dieses Lemma wird am Ende dieses Unterkapitels bewiesen. Während ein beliebiges Homotopie-Inverses von $x$ keine Abbildung mehr über $S^1$ sein wird, garantiert Lemma \ref{far:lemma2} immerhin, dass in einer Umgebung von $N$ das Homotopie-Inverse von $x$ faserweise gewählt werden kann. 

Betrachte das Pull-back $X'$ von $\Torus(Z,T')$ mit $\pi:\,\RR\lto S^1$; identifiziere ferner $Z$ mit einem Lift von $Z$ nach $X'$. In Lemma \ref{far:lemma1} haben wir festgestellt, dass die Inklusion $Z\lto X'$ auf Ebene der Fundamentalgruppen ein Isomorphismus ist. $Z$ teilt den Raum $X'$ in zwei zusammenhängende Komponenten $A'$ und $B'$; es sei $B'$ diejenige Komponente mit $T(B')\subset B'$. Es existiert dann ein Lift der Abbildung $g$ aus dem obigen Lemma zu einer Abbildung $g_1:\,X\lto X'$ mit
\[g_1^{-1}(Z)=N,\quad g_1^{-1}(A')=A,\quad g_1^{-1}(B')=B.\]
Betrachte nun das kommutative Diagramm mit exakten Zeilen, dessen obere Zeile die kurze exakte Sequenz aus Bemerkung \ref{far:kurze_exakte_sequenzen} ist:
\[\xymatrix{
0 \ar[r] & H_{j+1}(X,B;\ZZ[G]) \ar[r] \ar[d]^{(g_1,g_1\vert_B)_*} & H_j(B;\ZZ[G]) \ar[r] \ar[d]^{(g_1\vert_B)_*} & H_j(X;\ZZ[G]) \ar[r] \ar[d]^{g_{1*}}_{\cong} & 0 \\
& 0 \ar[r] & H_j(B';\ZZ[G]) \ar[r] & H_j(X';\ZZ[G]) \ar[r] & 0
}\]
Dabei ist $H_{j+1}(X',B';\ZZ[G])=\dirlim_{n\to\infty} H_{j+1}(T^{-n}(B'),B';\ZZ[G])=0$. Da $N$ ein $s$-bi-zusammenhängendes Splitting für $f$ und damit ein $(n-s-1)$-bi-zusammenhängendes Splitting für $-f$ ist, ist $H_{j+1}(X,B;\ZZ[G])$ nur für $j=n-s-1$ von 0 verschieden, und aus dem Diagramm folgt:
\begin{equation}\label{far:iso}
H_{n-s}(X,B;\ZZ[G]) \;\cong\; \ker \bigl((g_1\vert_B)_*:\,H_*(B;\ZZ[G])\lto H_*(B';\ZZ[G])\bigr).
\end{equation}
Weiter gilt:
\begin{enumerate}\setlength{\itemsep}{0pt}
\item $(g_1\vert_B)_*:\,H_j(B;\ZZ[G])\lto H_j(B';\ZZ[G])$ ist surjektiv für alle $j$,
\item $(g_1\vert_B)_*$ ist injektiv außer für $j=n-s-1$,
\item $\ker (g_1\vert_B)_*$ ist ein projektiver $\ZZ[G]$-Modul (nach Lemma \ref{far:lemma3}),
\item $g$ ist zellulär (nach Lemma \ref{far:lemma2}),
\item $N$, $A$ und $B$ sind Unter-CW-Komplexe von $X$ (nach Wahl der Triangulierung von $M$).
\end{enumerate}
Wir wenden nun das folgende Lemma von Farrell (zusammen mit einer Bemerkung auf Seite 333) an:

\begin{lem}[Seite 332f.]
Unter den Voraussetzungen \textup{(}i\textup{)} bis \textup{(}v\textup{)} von oben ist $\ker (g_1\vert_B)_*$ ein endlich erzeugter projektiver $\ZZ[G]$-Modul. Die Aktion von $T\in\pi_1(S^1)$ auf $G$ induziert einen $\alpha^{-1}$-semi-linearen Endomorphismus $t_*$ auf $\ker (g_1\vert_B)_*$, und es gilt:
\[p(\tau(g)) = (-1)^{n-s}\, \bigl[(\ker (g_1\vert_B)_*, t_*)\bigr] \in C(\ZZ[G],\alpha^{-1}).\]
\end{lem}

Es bleibt zu bemerken, dass der Isomorphismus \eqref{far:iso} als Teil der langen exakten Homologiesequenz mit der Aktion von $\pi_1(S^1)$ verträglich ist, sodass gilt:
\[p(\tau_2) = (-1)^{n-s}\, \bigl[(\ker (g_1\vert_B)_*, t_*)\bigr] = (-1)^{n-s}\,\bigl[(H_{n-s}(X,B;\ZZ[G]),t_*)\bigr]=c(-f).\]
Damit ist Satz \ref{s1:satz2} bewiesen. Es bleibt noch ein Beweis nachzutragen:

\begin{proof}[Beweis von Lemma \ref{far:lemma2}]
Fixiere eine endliche CW-Struktur auf $\Torus(Z,T')$, sodass der Teilraum $y^{-1}(U)=Z\times U$ ein Unterkomplex ist, der eine Produkt-CW-Struktur trägt. Die Abbildung $x:\,\Torus(Z,T')\lto M^f$ wird induziert vom folgenden Diagramm:
\[\xymatrix{
y^{-1}(U) \ar[d]^{x_1}_{\simeq} && y^{-1}(\partial U) \ar[d]^{x_0}_{\simeq} \arinclinv[ll] \arincl[rr] && y^{-1}(\overline{S^1-U}) \ar[d]^{x_2}_{\simeq} \\
p_f^{-1}(U) && p_f^{-1}(\partial U) \arinclinv[ll] \arincl[rr] && p_f^{-1}(\overline{S^1-U})
}\]
Hierbei bezeichnen die Abbildungen $x_i$ die jeweiligen Einschränkungen von $x$. 

Gemäß Proposition \ref{theta:prop3} wird also ein Homotopie-Inverses von $x$ induziert von geeigneten Homotopie-Inversen der Abbildungen $x_i$. Sei zunächst $(x_1^{-1},x_0^{-1})$ faserhomotopie-invers zu $(x_1,x_0)$. Dann ist wegen der Zusammenziehbarkeit von $U$ die Abbildung
\[x_1^{-1}\circ\lambda:\, N\times U\lto Z\times U\]
homotop über $U$ zu einer Abbildung des Typs $f'\times\id_U$ mit einer zellulären Abbildung $f':\,N\lto Z$ (vgl.~etwa den Beweis von Proposition \ref{fas:trivialisierung}.) Wir bezeichnen mit $(H_1,H_0):\,N\times(U,\partial U)\times I\lto Z\times(U,\partial U)$ die entsprechende Homotopie über $U$. 

Sei weiter $x_2^{-1}$ ein Homtopie-Inverses von $x_2$, sodass $(x_2^{-1},x_0^{-1})$ homotopie-invers als Abbildung von Paaren zu $(x_2,x_0)$ ist, und sei $x^{-1}$ das von den $x_i^{-1}$ induzierte Homotopie-Inverse von $x$. Auf Grund der Kofaserungseigenschaft des Paares $(f^{-1}(\overline{S^1-U}),f^{-1}(\partial U))$ kann die Homotopie $H_0$ zu einer Homtopie $H_2:\,f^{-1}(\overline{S^1-U})\times I \lto y^{-1}(\overline{S^1-U})$ erweitert werden, sodass die Auswertung von $H_2$ an 0 gleich $x_2^{-1}\circ\lambda$ ist. Nach dem zellulären Approximationssatz können wir außerdem ohne Beschränkung der Allgemeinheit davon ausgehen, dass die Auswertung von $H_2$ an 1 zellulär ist. 

Die Homotopien $H_i$ definieren dann zusammen eine Homotopie $M\times I\lto \Torus(Z,T')$, deren Auswertung an 0 gleich $x^{-1}\circ\lambda$ ist und deren Auswertung an 1 alle geforderten Eigenschaften erfüllt. 
\end{proof}

\section{Farrells dritte Bedingung}

Wir gehen ab jetzt davon aus, dass die Abbildung $f:\,M\lto S^1$ die Bedingung $c(f)=0$ erfüllt. Gemäß dem vorigen Unterkapitel existiert also ein Splitting $N$, sodass $M_N$ ein $h$-Kobordismus ist. Sei ein derartiges Splitting für den Rest des Kapitels fest gewählt. Zu diesem Splitting übernehmen wir wieder die Bezeichnungen $A$ und $B$ als die beiden Zusammenhangskomponenten von $X-N$. Die kanonische Inklusion $\lambda:\,M\lto M^f$ induziert als Abbildung über $S^1$ eine Abbildung $\sigma^*\lambda:\,M_N=\sigma^*M\lto \sigma^*M^f$ zwischen den mit $\sigma:\,I\lto S^1$ zurückgezogenen Räumen.

\begin{lem}
\begin{enumerate}\smallsep
\item Die Inklusion $N\lto X$ ist eine Homotopie-Äquivalenz.
\item Die Abbildungen
\[\lambda\vert_N:\,N\lto F,\quad \sigma^*\lambda:\,M_N\lto\sigma^*M^f\]
sind Homotopie-Äquivalenzen. Insbesondere definiert also $\lambda\vert_N$ eine einfache Struktur auf $F$.
\end{enumerate}
\end{lem}

\begin{proof}[Beweis]
(i) Zunächst gilt, da $M_N=\overline{T(A)-A}$ ein $h$-Kobordismus ist, dasselbe für die $h$-Kobordismen $X_n:=\overline{T^n(A)-T^{-n}(A)}$ für alle $n\in\NN$. Infolgedessen ist auch die Inklusion $N\subset X_n$ für alle $n\in\NN$ eine Homotopie-Äquivalenz. Der Beweis der Behauptung (i) entspricht nun dem Beweis von Lemma \ref{far:lemma1}.

(ii) Sei $\bar{f}:\,X\lto\RR$ wieder die von $f$ induzierte Abbildung der Pull-backs. Die kanonische Inklusion $\lambda:\,X\lto X^{\bar{f}}$ überdeckt nach Bemerkung \ref{s1:bem1} die kanonische Inklusion $M\lto M^f$, sodass unter den Identifikation $M_N=X\vert_I$ und $\sigma^*M^f=\pi^*M^f\vert_I$ die Abbildung $\sigma^*\lambda$ der Einschränkung der kanonischen Inklusion von $X\lto X^{\bar{f}}$ entspricht. Im Diagramm
\[\xymatrix{
N \arincl[rr]^{\simeq} \ar[d]^{\lambda\vert_N} && M_N \arincl[rr]^{\simeq} \ar[d]^{\sigma^*\lambda} && X\ar[d]^{\lambda}_{\simeq} \\
F \arincl[rr]^{\simeq} && \sigma^*M^f \arincl[rr]^{\simeq} && X^{\bar{f}}
}\]
sind nun die oberen horizontalen Abbildungen Homotopie-Äquivalenzen nach Teil (i). Für die untere horizontalen Abbildungen gilt dasselbe, da $X^{\bar{f}}\lto\RR$ eine Faserung über einem zusammenziehbaren Raum ist. Aus der Kommutativität des Diagramms folgt nun die Behauptung.
\end{proof}

Im Sinne dieses Lemmas werden wir im Folgenden die Fundamentalgruppen von $N$, $M_N$, $X$, $X^f$ und $F$ mit $G$ identifizieren. Nun definieren wir Farrells zweites Hindernis. Dazu erinnern wir an die Aktion $\alpha$ von $\ZZ$ auf $G=\pi_1(N)$ durch Konjugation in $\pi_1(M)=G\rtimes_\alpha \ZZ$. Diese induziert eine Aktion von $\ZZ$ auf $\Wh(G)$ durch Gruppenhomomorphismen. Setze $\Wh_\alpha(G):=\Wh(G)/\ZZ$ und
\[\tau(f) := [\tau(M_N,N)]\in\Wh_\alpha(G),\]
also die Restklasse der Torsion des $h$-Kobordismus $M_N$. Es sei darauf hingewiesen, dass dieses Element nicht mit dem von uns in Kapitel \ref{tau} definierten Element $\tau(f)$ übereinstimmt. Dessen Rolle übernimmt ja in diesem Kapitel die Torsion $\tau_2$.  Der folgende Satz fasst die Resultate von Farrel zusammen:

\begin{satz}[Seite 343f.]\label{s1:farrells_satz}
Das Element $\tau(f)$ ist definiert, falls $c(f)=0$, und hängt nur von der Homotopie-Klasse von $f$ ab. Es verschwindet genau dann, wenn ein Splitting $N'$ existiert, sodass $M_{N'}$ ein trivialer Kobordismus ist. Dies wiederum ist äquivalent dazu, dass $f$ homotop zu einem differenzierbaren Faserbündel ist. 
\end{satz}

Die von der Inklusion $j:\,G\lto G\rtimes_\alpha\ZZ$ induzierte Abbildung $j_*:\,\Wh(G)\lto\Wh(G\rtimes_\alpha\ZZ)$ faktorisiert nach \cite{farrell_hsiang}, Seite 210, in einen Monomorphismus $\Wh_\alpha(G)\lto\Wh(G\rtimes_\alpha\ZZ)$, den wir ebenfalls mit $j_*$ bezeichnen. 

\begin{satz}\label{s1:satz1}
Im Falle $\tau_1=0$ gilt:
\[\tau_2 = -j_*\,\tau(f).\]
Insbesondere verschwindet $\tau(f)$ genau dann, wenn auch $\tau_2$ verschwindet.
\end{satz}

\begin{proof}[Beweis]
Wir haben die Torsion der Komposition $x^{-1}\circ\lambda:\,M\lto M^f\lto\Torus(N,T')$ zu berechnen. Diese Abbildung wird auf Ebene der Push-outs vom folgenden kommutativen Diagramm induziert, in dem die unteren beiden Zeilen der Definition von $x:\,\Torus(N,T')\lto M^f$ in Diagramm \eqref{s1:diag3} entsprechen.
\[\xymatrix{
M_N \ar[d]^{\sigma^*\lambda}_{\simeq} && N\amalg N \arinclinv[ll] \ar[rr]^{\id\amalg\id} \ar[d]^{\lambda\vert_N\amalg\lambda\vert_N}_{\simeq} && N \ar[d]^{\lambda\vert_N}_{\simeq}\\
\sigma^*M^f && F\amalg F \arinclinv[ll] \ar[rr]^{\id\amalg\id} && F\\
N\times I \ar[u]_{\varphi} && N\amalg N \arinclinv[ll] \arincl[rr] \ar[u]^{\simeq}_{\lambda\vert_N\amalg(T\circ\lambda\vert_N)} && {\Cyl(\id\amalg T')} \ar[u]^{\simeq}_{\psi}
}\]
Somit lässt sich die gesuchte Torsion durch Anwenden von Lemma \ref{theta:lemma4} aus den Torsionen der jeweiligen vertikalen Abbildungen berechnen.

Für die Abbildung in der rechten Spalte gilt $\tau(\psi^{-1}\circ\lambda\vert_N)=0$, denn nach Konstruktion repräsentieren $\psi$ und $\lambda\vert_N$ dieselbe einfache Struktur auf $F$. Weiter gilt nach Definition von $\tau_1$:
\[j_*\,\tau((\lambda\vert_N)^{-1}\circ T^{-1}\circ\lambda\vert_N)=-\tau_1=0.\]
Wir nennen $r:\,M_N\lto N$ ein Homotopie-Inverses der Inklusion und erhalten nach Lemma \ref{theta:lemma4}:
\[\begin{split}
\tau_2 & =\; j_*\,\tau\bigl[M_N \xrightarrow{\sigma^*\lambda} \sigma^*M^f \xrightarrow{\varphi^{-1}} N\times I\bigr]\\
& = \; j_*\, \tau\bigl[M_N \xrightarrow{r} N = N\times\{0\} \hookrightarrow N\times I\bigr]\\
& = \; j_*\, \tau\bigl[M_N\xrightarrow{r} N\bigr]\\
& = \; -j_*\, \tau(M_N,N)
\end{split}\]
Hieraus folgt die Behauptung.
\end{proof}

\begin{satz}
Sei $M$ eine geschlossene und zusammenhängende differenzierbare Mannigfaltigkeit der Dimension $n\geq 6$, $f:\,M\lto S^1$ eine Abbildung, sodass $f_*:\,\pi_1(M)\lto\pi_1(S^1)$ surjektiv ist. Dann sind die folgenden Bedingungen äquivalent:
\begin{enumerate}\smallsep
\item $F$ hat den Homotopietyp eines endlichen CW-Komplexes, $\tau_1=0$ und $\tau_2=0$,
\item $X$ hat den Homotopietyp eines endlichen CW-Komplexes, $c(f)=0$ und $\tau(f)=0$,
\item $f$ ist homotop zu einem differenzierbaren Faserbündel.
\end{enumerate}
\end{satz}

\begin{proof}[Beweis]
Die Implikation von (i) nach (ii) liefern Proposition \ref{s1:prop1}, Satz \ref{s1:satz2} und Satz \ref{s1:satz1}. Die Implikation von (ii) nach (iii) ist das Resultat von Farrell, vgl.~Satz \ref{s1:farrells_satz}. Die Implikation von (iii) nach (i) schließlich folgt aus den Propositionen \ref{theta:prop_haupt} und \ref{tau:prop_haupt}.
\end{proof}

\backmatter

\bibliography{bibli}


\include{formales}

\end{document}